\documentclass[a4paper,11pt]{article}

\usepackage{authblk}

\usepackage{amssymb}
\let\Emptyset\emptyset

\usepackage{type1cm}        % activate if the above 3 fonts are
\usepackage{makeidx}         % allows index generation
\usepackage{graphicx}        % standard LaTeX graphics tool
\usepackage{multicol}        % used for the two-column index

\usepackage{newtxtext}       % 
\usepackage{newtxmath}       % selects Times Roman as basic font

\makeindex             % used for the subject index
\usepackage{mdframed}

\usepackage{tikz}
\usetikzlibrary{shapes,arrows}
\usetikzlibrary{patterns}

\usepackage{pgfplots}
\pgfplotsset{compat=1.13}
\usepgfplotslibrary{fillbetween}

\definecolor{PrimalOrange}{RGB}{255,153,85}
\definecolor{DualBlue}{RGB}{0,255,255}
\definecolor{SlabGreen}{RGB}{170,212,0}
\definecolor{HSUred}{RGB}{197,0,66}
\definecolor{mountainmeadow}{rgb}{0.19, 0.73, 0.56}
\definecolor{navyblue}{rgb}{0.0, 0.0, 0.5}

\definecolor{Goldenrod}{RGB}{255,185,15}
\definecolor{Eggshell}{RGB}{255,236,139}
\definecolor{LightOlive}{RGB}{204,255,153}

\usepackage{hyperref}
\usepackage{xcolor}
\usepackage{subfig}
\usepackage{graphicx}

\renewcommand{\emptyset}{\Emptyset}

\newcommand{\doi}[1]{\href{https://doi.org/#1}{doi:#1}}

\newcommand{\concentration}{{u}}
\newcommand{\convection}{\mathbf{v}}

\newcommand{\dualz}{{z}}

\newcommand{\phiflow}{\boldsymbol{\varphi}}
\newcommand{\psiflow}{\boldsymbol{\psi}}
\newcommand{\chiflow}{\chi}
\newcommand{\trans}{\textrm{t}}
\newcommand{\flow}{\textrm{f}}
\newcommand{\dualzflow}{\mathbf{z}} % displacement vector
\newtheorem{defi}{Definition}[section]
\newtheorem{theorem}[defi]{Theorem}

\newtheorem{remark}[defi]{Remark}

\newenvironment{mproof}{\paragraph{Proof.}}{\hfill$\blacksquare$}

\numberwithin{equation}{section}
\numberwithin{table}{section}
\numberwithin{figure}{section}

\usepackage{geometry}
\begin{document}

\title{\Large \textbf{A Cost-Efficient Space-Time Adaptive Algorithm for Coupled Flow and Transport}}
\author[M.\ P.\ Bruchh\"auser, M.\ Bause]
{\large \textbf{Marius Paul Bruchh\"auser}\thanks{bruchhaeuser@hsu-hamburg.de ($^\ast$corresponding author)} 
$\,\boldsymbol{\cdot}$
\textbf{Markus Bause}\thanks{bause@hsu-hamburg.de}\\
{\small Helmut Schmidt University, University of the German Federal Armed Forces Hamburg,\\
Faculty of Mechanical and Civil Engineering, Chair of Numerical Mathematics\\ 
Holstenhofweg 85, 22043 Hamburg, Germany}
}
\date{}
\maketitle

\begin{abstract}
% {\bfseries Abstract.} 
\noindent
In this work, a cost-efficient space-time adaptive algorithm based 
on the Dual Weighted Residual (DWR) method is developed and studied for a 
coupled model problem of flow and convection-dominated transport.
Key ingredients are a multirate approach adapted to varying dynamics in time
of the subproblems, weighted and non-weighted error indicators for the transport
and flow problem, respectively, and the concept of space-time slabs based on 
tensor product spaces for the data structure.
In numerical examples the performance of the underlying algorithm is studied
for benchmark problems and applications of practical interest. Moreover, the 
interaction of stabilization and goal-oriented adaptivity is investigated for 
strongly convection-dominated transport.
\end{abstract}


\section{Introduction} 
\label{sec:1:intro}
In this day and age, the efficient numerical approximation of multi-physics and
multi-scale problems is associated with an ever-increasing complexity; cf., e.g.,
\cite{Ahmed2022,Almani2021,Soszynska2021,Jammoul2021,Almani2019,Ge2018,Larson2007}.
Challenges include, inter alia, different characteristic time scales of the subproblems,
different characteristic spatial scales by means of layers and sharp moving fronts,
and thus an efficient handling of the underlying discretization parameters
in space and time likewise.  
In order to reduce complexity, we propose in this work a multirate approach 
using different time step sizes adapted to the dynamics and characteristic scales
of the respective subproblems.
For a general review of multirate methods including a list of references we refer
to \cite{Gander2013,Gupta2016}.

Furthermore, with regard to efficiency reasons it is indispensable that 
additional adaptive mesh refinement strategies in space and time are necessary; 
cf., e.g., \cite{Endtmayer2021,Avijit2022,Larson2007}.
For this purpose, our multirate approach is combined with goal-oriented error 
control based on the Dual Weighted Residual (DWR) method \cite{Becker2001,Bangerth2003}.
In the DWR approach, an a posteriori error representation is derived in terms of
a user-chosen goal functional of physical relevance, where the local residuals 
are weighted by means of approximating an additional so-called dual problem.
Consisten with the underlying goal, local error indicators manage the adaptive
mesh refinement process automatically by marking the respective cells in space 
and time and thus further reduce an aspect regarding the above mentioned complexity.
In spite of all the advantages associated with this goal-oriented approach,
the main criticism is referred to increasing numerical costs by means of solving
the additional dual problem, in particular when using an higher-order finite 
elements approach.
Thus, in order to reduce numerical costs significantly, but obtain adaptive
meshes as efficient as possible at the same time, we present a cost-efficient
space-time adaptive algorithm applied to a model problem of coupled flow and 
transport. 
Here, the adaptive mesh refinements in space and time are based on weighted and
non-weighted local error indicators for the transport and flow problem, respectively. 
The transport problem is represented by a convection-diffusion-reaction equation
involving high dynamic behavior in time, whereas the flow problem is modeled by 
a viscous time-dependent Stokes flow problem.
More precisely, the underlying model problem is described by the following equations,
where the dimensionless Stokes flow system read as
\begin{equation}
\label{eq:1:stokes_problem}
\begin{array}{rcl @{\,\,}l @{\,\,}l @{\,}l}
\partial_t \mathbf{v}
- \nu \Delta \mathbf{v}
+ \nabla p &=& \mathbf{f}
& \textnormal{in} & Q & = \Omega \times I\,,\\[.5ex]
\nabla \cdot \mathbf{v} &=& 0
& \textnormal{in} & Q & = \Omega \times I\,,\\[.5ex]
\mathbf{v} &=& \mathbf{0} %\mathbf{v}_D
& \textnormal{on} & \Sigma
& = \partial\Omega \times I\,,\\[.5ex]
\mathbf{v} &=& \mathbf{v}_0
& \textnormal{on} & \Sigma_0 & = \Omega \times \{ 0 \}\,,
\end{array}
\end{equation}
and the convection-diffusion-reaction transport problem in dimensionless form is 
given by
\begin{equation}
\label{eq:2:transport_problem}
\begin{array}{rcl @{\,\,}l @{\,\,}l @{\,}l}
\partial_{t} u
- \nabla \cdot (\varepsilon \nabla u)
+ \mathbf{v} \cdot \nabla u 
+ \alpha u &=& g
& \textnormal{in} & Q & = \Omega \times I\,,\\[.5ex]
u &=& 0
& \textnormal{on} & \Sigma & = \partial\Omega \times I\,,\\[.5ex]
u &=&  u_0
& \textnormal{on} & \Sigma_0 & = \Omega \times \{ 0 \}\,.
\end{array}
\end{equation}
In \eqref{eq:1:stokes_problem}, \eqref{eq:2:transport_problem}, we denote by 
$Q=\Omega \times I$ the space-time domain, where 
$\Omega \subset \mathbb{R}^{d}$, with $d=2,3$, is a polygonal or polyhedral 
bounded domain with Lipschitz boundary $\partial\Omega$ and $I=(0,T]$, 
$0 < T < \infty$, is a finite time interval. 
Within the flow problem \eqref{eq:1:stokes_problem}, $\mathbf{v}$ is 
the velocity and $p$ is the pressure variable, the parameter $\nu>0$ is 
the dimensionless viscosity, $\mathbf{f}$ is a given volume force, and 
$\mathbf{v}_{0}$ is a given initial condition, respectively. 
Regarding the transport problem \eqref{eq:2:transport_problem}, we assume that
$0 < \varepsilon \leq 1$ is a constant diffusion coefficient, 
$\alpha > 0$ is a non-negative reaction coefficient,
$g$ is a given source of the unknown scalar quantity $u$, and
$u_0$ is a given initial condition, respectively. 
To measure different dynamics in time with regard to the two subproblems, we 
introduce parameter dependent characteristic times $t_{\textrm{transport}}$ and
$t_{\textrm{flow}}$, respectively, defined by
\begin{equation}
\label{eq:3:characteristic-times}
t_{\textrm{flow}} := \frac{L}{V}\,,
\qquad
t_{\textrm{transport}} := \min \Bigg\{\frac{L^2}{\varepsilon}\,;\, \frac{L}{V}\,;\,
\frac{1}{\alpha}\Bigg\}\,,
\end{equation}
where $L$ denotes the characteristic length of the domain $\Omega$
, for instance, its diameter, 
and $V$ denotes a characteristic velocity of the flow field $\mathbf{v} $;
cf.\cite{Bruchhaeuser2022,Gujer2008,Morgenroth2015} for more details.
These characteristic times can be understood as dimensionless time variables 
and serve as indicators for the underlying temporal meshes.

Finally, for the sake of physical realism, the transport problem is supposed to 
be convection-dominated by assuming high P\'{e}clet numbers, cf.~\cite{Burman2014,John2018}.
This poses a further challenge for finding numerical solutions avoiding non-physical 
oscillations or smearing effects close to sharp moving fronts and layers; 
cf. e.g., \cite{John2009}.
The application of stabilization techniques is a typical approach to overcome
these difficulties.
Here, we apply the streamline upwind Petrov--Galerkin (SUPG) method \cite{Hughes1979,Brooks1982}. 
In this context, we investigate the interaction of stabilization and goal-oriented
error control within our multirate approach.

This work is organized as follows. In Sec.~\ref{sec:2:weakdisc} we introduce the
variational formulations of the subproblems as well as our multirate approach
including the space-time discretizations. 
In Sec.~\ref{sec:3:error} we derive a posteriori error representations for both 
the flow and transport problem. 
In Sec.~\ref{sec:4:algorithm} we present the underlying algorithm and explain some
implementational aspects regarding our multirate approach. 
Numerical examples are given in Sec. \ref{sec:5:examples} and in Sec. \ref{sec:6:conclusion}
we summarize with conclusions and give some outlook for future work.

\section{Variational Formulation and Discretization}
\label{sec:2:weakdisc}

In this section, we introduce the variational formulation of the coupled model 
problem and explain some details about the discretization in space and time.

\subsection{Variational Formulation} 
\label{sec:2:1:weak}

The variational formulation of the Stokes flow problem \eqref{eq:1:stokes_problem}
reads as follows: 
\textit{For $\mathbf{f} \in L^2(I;(V^{\prime})^d)$ and 
$\mathbf{v}_0 \in V^d$, find $\mathbf{u}=\{ \mathbf{v}, p\} \in \mathcal{Y}:=\mathcal{Y}_1 \times \mathcal{Y}_2$, 
such that}
\begin{equation}
\label{eq:4:weak_stokes}
B(\mathbf{u})(\boldsymbol{\varphi}) = F(\boldsymbol{\psi}) \quad
\forall \boldsymbol{\varphi}=\{\boldsymbol{\psi}, \chi\} \in \mathcal{Y}\,,
\end{equation}
\textit{where the bilinear form $B(\cdot)(\cdot)$
and the linear form $F(\cdot)$ are defined by}
\begin{equation*}
\begin{array}{rcl}
B(\mathbf{u})(\boldsymbol{\varphi}) &:=&
\displaystyle
\int_{I} \big\{
(\partial_{t} \mathbf{v},\boldsymbol{\psi})
+ b(\mathbf{u})(\boldsymbol{\varphi})
\big\}\; \mathrm{d} t 
+ (\mathbf{v}(0),\boldsymbol{\psi}(0))\,,
\\[1.5ex]
F(\boldsymbol{\psi}) &:=&  
\displaystyle
\int_{I} \big\{
(\mathbf{f}, \boldsymbol{\psi})\big\}\; \mathrm{d} t
+ (\mathbf{v}_0,\boldsymbol{\psi}(0))
\,,
\end{array}
\end{equation*}
\textit{with the inner bilinear form $b(\cdot)(\cdot)$ given by}
\begin{equation}
\label{eq:5:inner_bilinearform_stokes}
b(\mathbf{u})(\boldsymbol{\varphi}):=\nu(\nabla \mathbf{v}, \nabla \boldsymbol{\psi}) 
- (p, \nabla \cdot \boldsymbol{\psi}) + (\nabla \cdot \mathbf{v},\chi)\,.
\end{equation}
Here, $(\cdot, \cdot)$ denotes the inner product of $L^2(\Omega)$ or duality
pairing of $H^{-1}(\Omega)$ with $H^1_0(\Omega)$, respectively, and the 
appearing function spaces are given as follows
\begin{equation}
\label{eq:6:continuous_spaces_stokes}
\begin{array}{rcl}
\mathcal{Y}_1 & := &\{ \mathbf{v} \in L^2(0, T; V^d) \mid
\partial_t \mathbf{v} \in L^2(0, T; (V^\prime)^d) \}\,, \textrm{ with }
V:=\{\mathbf{v}\in H^1_0(\Omega) \mid \nabla\cdot\mathbf{v}=0\}\,,
\\[1.5ex]
\mathcal{Y}_2 & := & \{p \in L^2(0,T;H)\}\,, \textrm{ with }
H=L_0^2(\Omega):=\{q\in L^2(\Omega) \mid \int_\Omega q\; \mathrm{d}x = 0\}\,.
\end{array}
\end{equation}
\begin{remark}
\label{rem:1:boundary}
In some numerical examples, we also consider more general boundary conditions as 
introduced in \eqref{eq:1:stokes_problem}. More precisely, we consider a 
boundary partition $\partial\Omega=\Gamma_{\textnormal{in}} \cup \Gamma_{\textnormal{wall}}
\cup \Gamma_{\textnormal{out}}$ given by
\begin{equation}
\label{eq:7:boundary_stokes}
\mathbf{v}=\mathbf{v}_D \textrm{ on } 
\Gamma_{\textnormal{in}} \times I\,,
\quad
\mathbf{v}=\mathbf{0} \textrm{ on } 
\Gamma_{\textnormal{wall}} \times I\,,
\quad
(\nu\nabla\mathbf{v} - p\mathbf{I})\mathbf{n} = \mathbf{0} \textrm{ on } 
\Gamma_{\textnormal{out}} \times I\,.
\end{equation}
For this configuration, including a so-called "do nothing" outflow condition, 
the space $\mathcal{Y}_2$ has to be modified by $\mathcal{Y}_2 = \{p \in L^2(0,T;L^2(\Omega))\}\,.$ 
For more information about these boundary conditions as well as results 
regarding existence and uniqueness of solutions, we refer to \cite{Ern2021,John2016}. 
\end{remark}
Using the velocity solution $\mathbf{v}$ of \eqref{eq:1:stokes_problem}, the 
variational formulation of the transport problem \eqref{eq:2:transport_problem}
reads as follows: 
\textit{For a given $\mathbf{v} \in \mathcal{Y}_1$ of \eqref{eq:1:stokes_problem}, 
find $u \in \mathcal{X}$ such that}
\begin{equation}
\label{eq:8:weak_transport}
A(u;\mathbf{v})(\varphi) = G(\varphi) \quad \forall \varphi \in \mathcal{X}\,,
\end{equation}
\textit{where the bilinear form $A(\cdot;\cdot)(\cdot)$ and the linear form $G(\cdot)$
are defined by}
\begin{equation*}
\begin{array}{rcl}
A(u;\mathbf{v})(\varphi) & := & 
\displaystyle
\int_{I}\big\{(\partial_t u,\varphi)
+ a(u;\mathbf{v})(\varphi)
\big\} \mathrm{d} t\;
+ (u(0),\varphi(0) )\,,
\\[1.5ex]
G(\varphi) & := & \displaystyle\int_I(g,\varphi)\;\mathrm{d}t + (u_{0},\varphi(0) ) \,,
\end{array}
\end{equation*}
\textit{with the inner bilinear form $a(\cdot;\cdot)(\cdot)$ given by}
\begin{equation}
\label{eq:9:inner_bilinearform_transport}
a(u;\mathbf{v})(\varphi) := (\varepsilon\nabla u, \nabla \varphi) 
+(\mathbf{v}\cdot \nabla u, \varphi) 
+ (\alpha u,\varphi)\,,
\end{equation}
where the function space $\mathcal{X}$ is defined as 
$\mathcal{X} := \{ \mathbf{v} \in L^2(0, T; H^1_0(\Omega)) \mid
\partial_t \mathbf{v} \in L^2(0, T; H^{-1}(\Omega)) \}\,.$
\begin{remark}
\label{rem:2:coupling}
We note that the coupling here is uni-directional via the flow field $\mathbf{v}$
from the flow to the transport problem.
Well-posedness and the existence of a sufficiently regular solution, such that 
all of the arguments and terms used below are well-defined, are tacitly assumed 
without mentioning all technical assumptions about the data and coefficients 
explicitly, cf.~\cite{Roos2008} and \cite{Ern2021}.
\end{remark}

\subsection{Multirate and Discretization in Space and Time} 
\label{sec:2:2:discretization}

Here, we present our multirate in time approach as well as the underlying 
discretization schemes in space and time. 
Since the multirate approach is explained in detail in our previous work 
\cite{Bruchhaeuser2022} and the discretization is rather standard, we keep this 
section short and refer to our former works \cite{Bruchhaeuser2022a,Bruchhaeuser2022,Bause2021} 
for more details.

Assuming a high dynamic behavior in time within the transport problem compared 
to a rather slow behavior of the flow problem, along with with characteristic 
times $t_{\textrm{transport}} \ll t_{\textrm{flow}}$, we are using a finer 
temporal mesh for the transport problem compared to a rather coarse mesh used 
for the flow problem. 
In addition, we allow for adaptive time refinements within both subproblems, 
obtaining temporal meshes being as efficient as possible.
For this multirate decoupling, we divide the time interval $I$ into not 
necessarily equidistant, left-open subintervals $I_n:=(t_{n-1},t_n]\,,$ with 
$n=1,\dots,N\,,$ where $0=:t_0<\dots<t_N:=T$.
Since we take different temporal meshes as a basis for the two subproblems, 
we use indexes $\mathrm{f}$ and $\mathrm{t}$ to distinguish between flow and 
transport here and in the following. 
For simplicity of the implementation, we ensure that each element of the flow set
$\{ t_0^\mathrm{f}, \dots, t_{N^{\mathrm{f}}}^\mathrm{f} \}$
corresponds to an element of the transport set
$\{ t_0^\mathrm{t}, \dots, t_{N^{\mathrm{t}}}^{\mathrm{t}} \}$ such that the 
temporal mesh of the flow problem is not finer than that of the transport problem.
With this in mind, we consider a separation of the global space-time domain
$Q=\Omega \times I$ into a partition of so-called space-time slabs 
$\hat{Q}_n = \Omega \times I_n$.
The time domain of each space-time slab $\hat{Q}_n$ is then discretized using a 
one-dimensional triangulation $\mathcal{T}_{\sigma,n}$ or $\mathcal{T}_{\tau,n}$ 
for the subinterval $I_n^\mathrm{f}$ or $I_n^{\trans}$ of the flow or transport
problem, respectively.
This allows to have more than one cell in time on a slab $\hat{Q}_n$ and a 
different number of cells in time of pairwise different slabs $\hat{Q}_i$ and 
$\hat{Q}_j\,, 1 \le i,j \le N$.
Furthermore, within this context let $\mathcal{F}_\sigma$ and $\mathcal{F}_\tau$ 
be the sets of all interior time points given as
\begin{equation*}
\begin{array}{rcl}
\mathcal{F}_\sigma &:=& ( \{ t_1^\mathrm{f}, \dots, 
t_{N^{\mathrm{f}}}^\mathrm{f} \}
\cup
\{ t^\mathrm{f} \in \partial K_n^\mathrm{f} \mid K_n^\mathrm{f} \in \mathcal{T}_{\sigma,n} \} ) \setminus \{ 0, T \}\,,
\\[1.5ex]
\mathcal{F}_\tau &:=& ( \{ t_1^{\mathrm{t}}, \dots, t_{N^{\mathrm{t}}}^{\mathrm{t}} \}
\cup
\{ t^\mathrm{t} \in \partial K_n^{\mathrm{t}} \mid K_n^{\mathrm{t}} \in \mathcal{T}_{\tau,n} \} ) \setminus \{ 0, T \}\,.
\end{array}
\end{equation*}
The commonly used time step size $\sigma_K$ or $\tau_K$ is here the diameter or 
length of the cell in time $K_n^\mathrm{f}$ of $\mathcal{T}_{\sigma,n}$ or 
$K_n^\mathrm{t}$ of $\mathcal{T}_{\tau,n}$ and the global time discretization 
parameter $\sigma$ or $\tau$ is the maximum time step size $\sigma_K$ or $\tau_K$ 
of all cells in time of all slabs $\hat{Q}_n^\mathrm{f}$ or $\hat{Q}_n^\mathrm{t}$,
respectively.

For the discretization in time, we use a discontinuous Galerkin method dG($r$) 
with an arbitrary polynomial degree $r \ge 0$. 
Then, the time-discrete function spaces for the flow and transport problem, respectively,
are given by
\begin{equation}
\label{eq:10:Def_X_time_spaces}
\begin{array}{r@{\,}c@{\,}l@{\,}}
\mathcal{Y}_{\sigma}^{r} &:= \Big\{ &
\{ \mathbf{v}_{\sigma},p_\sigma \} \in 
L^2(0, T; V^d\times H) \,\,\big|
\mathbf{v}_{\sigma}|_{K_n^\mathrm{f}} \in \mathcal{P}_r(K_n^\mathrm{f}; V^d)\,,
\mathbf{v}_{\sigma}(0)\in L^2(\Omega)\,,\\[1.5ex]
 &~& p_{\sigma}|_{K_n^\mathrm{f}} \in \mathcal{P}_r(K_n^\mathrm{f}; H)\,, 
 K_n^\mathrm{f} \in \mathcal{T}_{\sigma,n}\,,\,\,
 n=1,\dots,N^{\mathrm{f}}
\Big\}\,,
\\[1.5ex]
\mathcal{X}_{\tau}^{r} &:= \Big\{ &
 u_{\tau} \in L^2(0, T; H_0^1(\Omega)) \,\,\big|\,\,
 u_{\tau}|_{K_n^\mathrm{t}} \in \mathcal{P}_r(K_n^\mathrm{t}; H_0^1(\Omega))\,,
 \\[1.5ex]
 &~& K_n^\mathrm{t} \in \mathcal{T}_{\tau,n}\,,\,\,
 n=1,\dots,N^\mathrm{t}\,,\,\,
 u_{\tau}(0)\in L^2(\Omega)
\Big\}\,,
\end{array}
\end{equation}
where $\mathcal{P}_{r}(K_n; V)$ denotes the space of all polynomials
in time up to degree $r \ge 0$ on $K_n$ with values in some function space $V$.
Exemplary for some discontinuous in time function $u_{\tau} \in \mathcal{X}_{\tau}^{r}$,
we define the limits $u_{\tau}(t_F^\pm)$
from above and below as well as their jump at $t_F$ by
\begin{displaymath}
u_\tau(t_F^\pm) := \displaystyle \lim_{t \to t_F \pm 0} u_\tau(t)\,,\quad
[ u_\tau ]_{t_F} := u_\tau(t_F^+) - u_\tau(t_F^-) \,.
\end{displaymath}
Now, the semi-discrete in time scheme of the the flow problem
\eqref{eq:4:weak_stokes} reads as follows:
\textit{Find $\mathbf{u}_{\sigma} = \{\mathbf{v}_{\sigma},p_{\sigma}\} \in \mathcal{Y}_{\sigma}^{r}$ 
such that}
\begin{equation}
\label{eq:11:semi_stokes}
B_{\sigma}(\mathbf{u}_{\sigma})(\phiflow_{\sigma})
=
F_{\sigma}(\psiflow_{\sigma}) 
\quad \forall \phiflow_{\sigma} = \{\psiflow_{\sigma},\chiflow_{\sigma}\} 
\in \mathcal{Y}_{\sigma}^{r}\,,
\end{equation}
\textit{where the bilinar form $B_{\sigma}(\cdot,\cdot)(\cdot,\cdot)$ and the linear form
$F_{\sigma}(\cdot)$ are defined by}
\begin{equation}
\label{eq:12:BsigmaFsigma}
\begin{array}{rcl}
\displaystyle
B_{\sigma}(\mathbf{u}_{\sigma})(\phiflow_{\sigma})
& := & 
\displaystyle
\sum_{n=1}^{N^{\flow}}\sum_{K_{n}^{\flow}\in\mathcal{T}_{\sigma, n}}\int_{K_{n}^{\flow}}
\big\{
(\partial_t \mathbf{v}_{\sigma},\psiflow_{\sigma})
+ b(\mathbf{u}_{\sigma})(\phiflow_{\sigma})
\big\} \mathrm{d} t\; 
\\[1.5ex]
& &
\displaystyle
+ \sum_{t_{\mathcal{F}}\in\mathcal{F}_{\sigma}}
([\mathbf{v}_{\sigma}]_{t_{\mathcal{F}}},\psiflow_{\sigma}(t_{\mathcal{F}}^+) )
+ (\mathbf{v}_{\sigma}(0^{+}),\psiflow_{\sigma}(0^+) )\,,
\\[1.5ex]
F_{\sigma}(\psiflow_{\sigma}) 
& := & 
\displaystyle
\sum_{n=1}^{N^{\flow}}\sum_{K_{n}^{\flow}\in\mathcal{T}_{\sigma, n}}\int_{K_{n}^{\flow}}
(\mathbf{f},\psiflow_{\sigma})\;\mathrm{d}t
+ (\mathbf{v}_{0},\psiflow_{\sigma}(0^+) )\,,
\end{array}
\end{equation}
\textit{with the inner bilinear form $b(\cdot)(\cdot)$ given by 
Eq.~\eqref{eq:5:inner_bilinearform_stokes}}.
The semi-discrete in time scheme of the the transport problem
\eqref{eq:8:weak_transport} reads as follows:
\textit{For a given $\mathbf{v}_{\sigma} \in \mathcal{Y}_{\sigma}^{r}$ 
of \eqref{eq:11:semi_stokes}, find $u_{\tau} \in \mathcal{X}_{\tau}^{r}$ such that}
\begin{equation}
\label{eq:13:semi_transport}
A_{\tau}(u_{\tau};\mathbf{v}_{\sigma})(\varphi_{\tau}) = G_{\tau}(\varphi_{\tau}) \quad 
\forall \varphi_{\tau} \in \mathcal{X}_{\tau}^{r}\,,
\end{equation}
\textit{where the bilinear form $A_{\tau}(\cdot;\cdot)(\cdot)$ and the linear 
form $G_{\tau}(\cdot)$ are defined by}
\begin{equation}
\label{eq:14:AtauGtau}
\begin{array}{rcl}
A_{\tau}(u;\mathbf{v})(\varphi) & := & 
\displaystyle
\sum_{n=1}^{N^{\trans}}\sum_{K_{n}^{\trans}\in\mathcal{T}_{\tau, n}}\int_{K_{n}^{\trans}}
\big\{
(\partial_t u_{\tau},\varphi_{\tau})
+ a(u_{\tau};\mathbf{v}_{\sigma})(\varphi_{\tau})
\big\} \mathrm{d} t\;
\\[1.5ex]
& &
\displaystyle
+ \sum_{t_{\mathcal{F}}\in\mathcal{F}_{\tau}}
([u_{\tau}]_{t_{\mathcal{F}}},\varphi_{\tau}(t_{\mathcal{F}}^+) )
+ (u_{\tau}(0^{+}),\varphi_{\tau}(0^+) )\,,
\\[1.5ex]
G_{\tau}(\varphi_{\tau}) & := & 
\displaystyle
\sum_{n=1}^{N^{\trans}}\sum_{K_{n}^{\trans}\in\mathcal{T}_{\tau, n}}\int_{K_{n}^{\trans}}
(g,\varphi_{\tau})\;\mathrm{d}t + (u_{0},\varphi_{\tau}(0^+) ) \,,
\end{array}
\end{equation}
\textit{with the inner bilinear form $a(\cdot;\cdot)(\cdot)$ given by 
Eq.~\eqref{eq:9:inner_bilinearform_transport} depending on the semi-discrete
flow solution $\mathbf{v}_{\sigma}$}.

For the discretization in space we use standard Lagrange type finite element 
spaces of continuous functions that are piecewise polynomials. 
More precisely, we define the discrete finite element space 
$
V_{h}^{p,n}:=
\big\{v\in C(\overline{\Omega})\mid v_{|K}
\in Q_h^p(K)\,,\forall K\in\mathcal{T}_{h,n}\,,1 \le n \le N
\big\}
$
on a triangulation $\mathcal{T}_{h,n}$ building a decomposition of the domain
$\Omega$ into disjoint elements $K$. Here, the space $Q_h^p(K)$ is defined on 
the reference element with maximum degree $p$ in each variable. 
Now, by replacing the respective spaces in the definition of the semi-discrete 
function spaces $\mathcal{Y}_{\sigma}^{r}$ and $\mathcal{X}_{\tau}^{r}$ by $V_h^{p,n}$, we obtain 
the fully discrete function spaces for the transport and flow problem, respectively, 
\begin{equation}
\label{eq:15:Def_space_spaces}
\begin{array}{r@{\,}c@{\,}l@{\,}}
\mathcal{Y}_{\sigma h}^{r,p_{\mathbf{v}},p_p} := & \Big\{ &
\{ \mathbf{v}_{\sigma h},p_{\sigma h} \} \in 
\mathcal{Y}_{\sigma}^{r} \,\,\big|\,\,
\mathbf{v}_{\sigma h}|_{K_n} \in \mathcal{P}_r(K_n; (H_h^{p_\mathbf{v},n})^d)\,,
\mathbf{v}_{\sigma h}(0)\in (H_h^{p_\mathbf{v},0})^d\,,
\\[1ex]
&~&
p_{\sigma h}|_{K_n} \in \mathcal{P}_r(K_n; L_h^{p_p,n})\,,
K_n \in \mathcal{T}_{\sigma,n}\,,\,\,
 n=1,\dots,N^{\flow}
\Big\}\,,
\\[1.5ex]
\mathcal{X}_{\tau h}^{r,p} := & \Big\{ &
u_{\tau h}\in \mathcal{X}_{\tau}^{r} \,\,\big|\,\,
u_{\tau h}|_{K_n} \in \mathcal{P}_r(K_n;H_h^{p_{u},n})\,,
u_{\tau h}(0) \in H_h^{p_{u},0},
K_n \in \mathcal{T}_{\tau,n}\,,\,\,n=1,\dots,N^{\trans}
\Big\}\,,
\end{array}
\end{equation}
\begin{displaymath}
H_h^{p_{\mathbf{v}},n}:=V_h^{p_{\mathbf{v}},n}\cap H_0^1(\Omega), \quad
L_h^{p_{p},n}:=V_h^{p_{p},n}\cap L_0^2(\Omega), \quad
H_h^{p_{u},n}:=V_h^{p_{u},n}\cap H_0^1(\Omega).
\end{displaymath}
We note that the spatial finite element space $V_h^{p,n}$ is allowed to be 
different on all subintervals $I_n$, which is natural in the context of a 
discontinuous Galerkin approximation of the time variable and allows dynamic 
mesh changes in time.
Due to the conformity of $H_h^{p_{u},n}$, $H_h^{p_{\mathbf{v}},n}$
and $L_h^{p_{p},n}$, we get
$\mathcal{X}_{\tau h}^{r,p}\subseteq \mathcal{X}_{\tau}^{r}$ and
$\mathcal{Y}_{\sigma h}^{r,p_{\mathbf{v}},p_p}\subseteq 
\mathcal{Y}_{\sigma}^{r}$, respectively.
Then, the fully discrete schemes of the flow and transport problem can be easily 
obtained from the semi-discrete in time schemes given by \eqref{eq:11:semi_stokes}
and \eqref{eq:13:semi_transport}, respectively, by simply adding the additional 
index $h$ to the variables and by replacing the respective semi-discrete spaces 
by the above defined fully discrete counterparts.
For the sake of completeness, the fully discrete scheme for the flow problem 
reads as follows:
\textit{Find $\mathbf{u}_{\sigma h} = \{\mathbf{v}_{\sigma h},p_{\sigma h}\} 
\in \mathcal{Y}_{\sigma h}^{r,p_{\mathbf{v}},p_p}$ such that}
\begin{equation}
\label{eq:16:fully_stokes}
B_{\sigma}(\mathbf{u}_{\sigma h})(\phiflow_{\sigma h})
=
F_{\sigma}(\psiflow_{\sigma h}) 
\quad \forall \phiflow_{\sigma h} = \{\psiflow_{\sigma h},\chiflow_{\sigma h}\} 
\in \mathcal{Y}_{\sigma h}^{r,p_{\mathbf{v}},p_p}\,,
\end{equation}
\textit{where the bilinar form $B_{\sigma}(\cdot,\cdot)(\cdot,\cdot)$ and the linear form
$F_{\sigma}(\cdot)$ are defined by \eqref{eq:12:BsigmaFsigma}}.

Finally, since the transport problem is assumed to be convection-dominated, the 
finite element approximation needs to be stabilized in order to avoid spurious 
and non-physical oscillations of the discrete solution arising close to sharp 
fronts and layers. 
Here, we apply the streamline upwind Petrov--Galerkin (SUPG) method introduced 
by Hughes and Brooks \cite{Hughes1979,Brooks1982}.
Then, the stabilized fully discrete scheme for the transport problem reads as follows:
\textit{For a given $\mathbf{v}_{\sigma h} \in \mathcal{Y}_{\sigma h}^{r,p_{\mathbf{v}},p_p}$ 
of \eqref{eq:16:fully_stokes}, find $u_{\tau h} \in \mathcal{X}_{\tau h}^{r,p}$ such that}
\begin{equation}
\label{eq:17:fully_transport}
A_{S}(u_{\tau h};\mathbf{v}_{\sigma h})(\varphi_{\tau h}) = G_{\tau}(\varphi_{\tau h}) \quad 
\forall \varphi_{\tau h} \in \mathcal{X}_{\tau h}^{r,p}\,,
\end{equation}
\textit{where the linear form $G_{\tau}(\cdot)$ is defined in \eqref{eq:14:AtauGtau}
and the stabilized bilinear form $A_S(\cdot;\cdot)(\cdot)$ is given by}
\begin{equation*}
A_{S}(u_{\tau h};\mathbf{v}_{\sigma h})(\varphi_{\tau h}):=
A_{\tau}(u_{\tau h};\mathbf{v}_{\sigma h})(\varphi_{\tau h})
+ S_{A}(u_{\tau h};\mathbf{v}_{\sigma h})(\varphi_{\tau h})\,,
\end{equation*}
\textit{with $A_{\tau}(\cdot;\cdot)(\cdot)$ being defined in \eqref{eq:14:AtauGtau}}.
\textit{Here, $S_A(\cdot;\cdot)(\cdot)$ is the SUPG stabilized bilinear form
obtained by adding weighted residuals.}
In order to keep this work short, we skip here the explicit 
presentation of $S_A$ that can be found, e.g., in
our work \cite[Sec. 2.5]{Bause2021}.

\section{A Posteriori Error Estimation for Coupled Flow and Transport} 
\label{sec:3:error}
In this section, we present DWR-based a posteriori error representations for the 
flow as well as the stabilized transport problem. 
The latter are depending, among other things, on additional so-called coupling 
terms that account for the influence of the error within the flow problem and may 
be interpreted as a modeling error, cf.~\cite{Larson2007}.
The ideas and concepts below are based on the works of Besier and Rannacher \cite{Besier2012} 
as well as Schmich and Vexler \cite{Schmich2008},  
where stabilized Navier-Stokes equations and parabolic problems in general have 
been investigated, respectively.
In order to keep this section short and clear, we restrict ourselves to
the presentation of the main results regarding the seperation of the temporal 
and spatial discretization errors for both subproblems %the flow and transport problem 
and refer to our works \cite{Bruchhaeuser2022a,Bruchhaeuser2022} for a detailed 
derivation and further details.

The following error representation formulas are given in terms of a user-chosen 
goal functional $J(\cdot)$ by using Lagrangian functionals $\mathcal{L}$ within 
a constraint optimization approach, cf.~\cite{Becker2001,Bangerth2003,Besier2012}.
We assume the goal $J(\cdot)$ to be a linear functional, generally  
given as a sum
\begin{equation*}
% \label{eq:?:Jgeneral}
 J(u)=\int_0^T J_1(u(t))\mathrm{d}t 
 + J_2(u(T))\,,
\end{equation*}
where $J_1$ and $J_2$ are three times differentiable functionals 
and each of them may be zero; cf.~\cite{Schmich2008,Besier2012}.
Here, we separate the error representation formulas into temporal and spatial
amounts such that their localized forms can be used as cell-wise error indicators 
for the adaptive mesh refinement process in space and time.

\subsection{An A Posteriori Error Estimator for the Flow Problem} 
\label{sec:3:1:ER_stokes}

To derive the temporal and spatial error representation formulas for the flow problem, 
we first define Lagrangian functionals $\mathcal{L}: \mathcal{Y} \times \mathcal{Y} \rightarrow \mathbb{R}$,
$\mathcal{L}_\sigma: \mathcal{Y}_{\sigma}^{r} \times \mathcal{Y}_{\sigma}^{r}\rightarrow \mathbb{R}$,
and $\mathcal{L}_{\sigma h}:\mathcal{Y}_{\sigma h}^{r,p_{\mathbf{v}},p_p} \times 
\mathcal{Y}_{\sigma h}^{r,p_{\mathbf{v}},p_p} \rightarrow \mathbb{R}$ 
in the following way:
\begin{equation}
\label{eq:18:Lagrangian_stokes}
\begin{array}{r@{\,}c@{\,}l@{\,}}
\mathcal{L}(\mathbf{u},\mathbf{z}) 
& := &
\displaystyle
J(\mathbf{u}) + F(\mathbf{w}) - B(\mathbf{u})(\mathbf{z})\,,
\\[1.5ex]
\mathcal{L}_{\sigma}(\mathbf{u}_{\sigma},\mathbf{z}_{\sigma}) 
& := & 
\displaystyle
J(\mathbf{u}_{\sigma}) + F(\mathbf{w}_{\sigma}) 
- B_{\sigma}(\mathbf{u}_{\sigma})(\mathbf{z}_{\sigma})\,,
\\[1.5ex]
\mathcal{L}_{\sigma h}(\mathbf{u}_{\sigma h},\mathbf{z}_{\sigma h}) 
& := &
\mathcal{L}_{\sigma}(\mathbf{u}_{\sigma h},\mathbf{z}_{\sigma h}) \,.
\end{array}
\end{equation}
The Lagrangian functionals are defined on different discretization levels such 
that the corresponding optimality or stationary conditions can be identified 
with primal and dual problems. While the primal problems correspond to the 
continuous, the semi-discrete in time and the fully discrete problems given by
\eqref{eq:4:weak_stokes}, \eqref{eq:11:semi_stokes} and \eqref{eq:16:fully_stokes},
respectively, the dual problems are kind of auxiliary problems providing
dual variables (Lagrange multipliers) $\mathbf{z}=\{\mathbf{w},q\},
\mathbf{z}_{\sigma}=\{\mathbf{w}_{\sigma},q_{\sigma}\}$ and 
$\mathbf{z}_{\sigma h}=\{\mathbf{w}_{\sigma h},q_{\sigma h}\}$.
These dual solutions are used for weighting the influence of the local residuals on the 
error within the underlying goal quantity. 

Now, using the Lagrangian functionals defined above, we derive the following a 
posteriori error representations for the flow problem in space and time.
%%%%%%%%%%%%%%%%%%%%%%%%%%%%%%%%%%%%%%%%%%%%%%%%%%%%%%%%%%%%%%%%%%%%%%%%%%%%%%%
% Begin THEOREM 3.1 EE Stokes
\begin{theorem}
\label{thm:1:ER_stokes}
Let $\{\mathbf{u},\dualzflow\} \in \mathcal{Y} \times \mathcal{Y}$,
$\{\mathbf{u}_{\sigma},\dualzflow_{\sigma}\} \in \mathcal{Y}_{\sigma}^{r} 
\times \mathcal{Y}_{\sigma}^{r}$, and
$\{\mathbf{u}_{\sigma h},\dualzflow_{\sigma h}\} \in 
\mathcal{Y}_{\sigma h}^{r,p_{\mathbf{v}},p_p} \times \mathcal{Y}_{\sigma h}^{r,p_{\mathbf{v}},p_p}$
be stationary points of
$\mathcal{L}, \mathcal{L}_{\sigma}$, and $\mathcal{L}_{\sigma h}$
on the different levels of discretization, i.e.,
\begin{equation*}
\begin{array}{r@{\;}c@{\;}l@{\;}}
\mathcal{L}^{\prime}(\mathbf{u},\dualzflow)(\delta \mathbf{u}, \delta \dualzflow)
=
\mathcal{L}_{\sigma}^{\prime}(\mathbf{u},\dualzflow)(\delta \mathbf{u}, \delta \dualzflow)
& = & 0 
\quad
\forall \,
\{\delta \mathbf{u}, \delta \dualzflow\} \in 
\mathcal{Y}  \times  \mathcal{Y}\,,
\\[1.5ex]
\mathcal{L}_{\sigma}^{\prime}(\mathbf{u}_{\sigma},\dualzflow_{\sigma})
(\delta \mathbf{u}_{\sigma}, \delta \dualzflow_{\sigma}) 
& = & 0 
\quad
\forall \, \{\delta \mathbf{u}_{\sigma}, \delta \dualzflow_{\sigma}\}
\in \mathcal{Y}_{\sigma}^{r}  \times  \mathcal{Y}_{\sigma}^{r}\,,
\\[1.5ex]
\mathcal{L}_{\sigma h}^{\prime}(\mathbf{u}_{\sigma h},\dualzflow_{\sigma h})
(\delta \mathbf{u}_{\sigma h}, \delta \dualzflow_{\sigma h})  
& = &
\\[1.5ex]
=\mathcal{L}_{\sigma}^{\prime}(\mathbf{u}_{\sigma h},\dualzflow_{\sigma h})
(\delta \mathbf{u}_{\sigma h}, \delta \dualzflow_{\sigma h}) 
& = & 0  
\quad
\forall \, \{\delta \mathbf{u}_{\sigma h}, \delta \dualzflow_{\sigma h}\}
\in \mathcal{Y}_{\sigma h}^{r,p_{\mathbf{v}},p_p}  \times \mathcal{Y}_{\sigma h}^{r,p_{\mathbf{v}},p_p}\,.
\end{array}
\end{equation*}
Then, for the discretization errors in space and time we get the representation 
formulas
\begin{subequations}
\begin{equation}
\label{eq:19:ER_stokes_time}
\begin{array}{r@{\;}c@{\;}l@{\;}}
J(\mathbf{u})-J(\mathbf{u}_{\sigma})  & = &
\frac{1}{2}\rho_{\sigma}(\mathbf{u}_{\sigma})
(\dualzflow-\tilde{\dualzflow}_{\sigma})
+
\frac{1}{2}
\rho_{\sigma}^{\ast}
(\mathbf{u}_{\sigma},\dualzflow_{\sigma})
(\mathbf{u}-\tilde{\mathbf{u}}_{\sigma})
+ \mathcal{R}_{\sigma}\,,
% \\[1.5ex]
\end{array}
\end{equation}
\begin{equation}
\label{eq:19:ER_stokes_space}
\begin{array}{r@{\;}c@{\;}l@{\;}}
J(\mathbf{u}_{\sigma})-J(\mathbf{u}_{\sigma h}) & = &
\frac{1}{2}\rho_{\sigma}(\mathbf{u}_{\sigma h})
(\dualzflow_{\sigma}-\tilde{\dualzflow}_{\sigma h})
+ \frac{1}{2}
\rho_{\sigma}^{\ast}
(\mathbf{u}_{\sigma h},\dualzflow_{\sigma h})
(\mathbf{u}_{\sigma}-\tilde{\mathbf{u}}_{\sigma h})
+ \mathcal{R}_{\sigma h}\,.
\end{array}
\end{equation}
\end{subequations}
Here, 
$\{\tilde{\mathbf{u}}_{\sigma},\tilde{\dualzflow}_{\sigma}\}
\in \mathcal{Y}_{\sigma}^{r} \times \mathcal{Y}_{\sigma}^{r}$, 
and 
$\{\tilde{\mathbf{u}}_{\sigma h},\tilde{\dualzflow}_{\sigma h}\} 
\in \mathcal{Y}_{\sigma h}^{r,p_{\mathbf{v}},p_p} \times \mathcal{Y}_{\sigma h}^{r,p_{\mathbf{v}},p_p}$
can be chosen arbitrarily and the remainder terms $\mathcal{R}_{\sigma}, 
\mathcal{R}_{\sigma h}$ are of higher-order with respect to the errors 
$\mathbf{u}-\mathbf{u}_{\sigma},\mathbf{z}-\mathbf{z}_{\sigma}$
$\mathbf{u}_{\sigma}-\mathbf{u}_{\sigma h},\mathbf{z}_{\sigma}-\mathbf{z}_{\sigma h}$,
respectively.
Furthermore, $\rho_{\sigma}$ and $\rho_{\sigma}^{\ast}$ denote the primal and 
dual residuals based on the semi-discrete in time schemes, where their explicit
definitions are given in the Appendix.
\end{theorem}
%%%%%%%%%%%%%%%%%%%%%%%%%%%%%%%%%%%%%%%%%%%%%%%%%%%%%%%%%%%%%%%%%%%%%%%%%%%%%%%
% End THEOREM 3.1 EE Stokes
\begin{mproof}
A detailed proof can be found in \cite[Ch.~5]{Bruchhaeuser2022a}, which is  
based on the idea given for parabolic problems in general that can be found 
in \cite[Thm.~3.2]{Schmich2008}. 
Basically, we are using a general result given in \cite[Prop.~3.1]{Schmich2008},
first derived in \cite{Becker2001}, with the following settings:
\begin{equation*}
\begin{array}{r@{\;}c@{\;}l@{\;} @{\,\,}l @{\,\,}l @{\,}l @{\,\,}l}
L & = & \mathcal{L}_\sigma\,,\;\;
Y   = (\mathcal{Y}+\mathcal{Y}_{\sigma}^{r}) \times (\mathcal{Y}+\mathcal{Y}_{\sigma}^{r})\,,\;\;
& Y_0 & = \mathcal{Y}_{\sigma}^{r} \times \mathcal{Y}_{\sigma}^{r}\,,
& \textnormal{ for \eqref{eq:19:ER_stokes_time}: }
\\[1.5ex]
L & = & \mathcal{L}_\sigma\,,\;\;
Y  = \mathcal{Y}_{\sigma}^{r} \times \mathcal{Y}_{\sigma}^{r}\,,\;\;
& Y_0 & = \mathcal{Y}_{\sigma h}^{r,p_{\mathbf{v}},p_p} \times 
\mathcal{Y}_{\sigma h}^{r,p_{\mathbf{v}},p_p}\,,
& \textnormal{ for \eqref{eq:19:ER_stokes_space}: }
\end{array}
\end{equation*}
where $\mathcal{L}_{\sigma}$ is the Lagrangian functional given by
\eqref{eq:18:Lagrangian_stokes} and $Y, Y_0$ are function spaces defined in 
\cite[Prop.~3.1]{Schmich2008}. 
\end{mproof}

\subsection{An A Posteriori Error Estimator for the Transport Problem} 
\label{sec:3:2:ER_transport}

To derive the temporal and spatial error representation formulas for the transport
problem, we define Lagrangian functionals $\mathcal{L}: \mathcal{X} \times \mathcal{X} \rightarrow \mathbb{R}$,
$\mathcal{L}_\tau: \mathcal{X}_{\tau}^{r} \times \mathcal{X}_{\tau}^{r}\rightarrow \mathbb{R}$,
and $\mathcal{L}_{\tau h}:\mathcal{X}_{\tau h}^{r,p} \times \mathcal{X}_{\tau h}^{r,p} \rightarrow \mathbb{R}$ 
by means of:
\begin{equation}
\label{eq:20:Lagrangian_transport}
\begin{array}{r@{\,}c@{\,}l@{\,}}
\mathcal{L}(u,\dualz;\convection) & := & J(u)
+ G(\dualz)
- A(u; \convection)(\dualz)\,,
\\[1.ex]
\mathcal{L}_{\tau}(u_\tau,\dualz_\tau;\convection_{\sigma}) & := &
J(u_{\tau}) + G_\tau(\dualz_{\tau})
- A_{\tau}(u_{\tau}; \convection_{\sigma})(\dualz_{\tau})\,,
\\[1.ex]
\mathcal{L}_{\tau h}(u_{\tau h},\dualz_{\tau h};\convection_{\sigma h}) 
& := &
J(u_{\tau h})
+ G_\tau (\dualz_{\tau h})
- A_S(u_{\tau h}; \convection_{\sigma h})(\dualz_{\tau h})\,.
\end{array}
\end{equation}
Using these Lagrangian functionals, we derive the following a 
posteriori error representations for the transport problem in space and time.
Compared to the result given in Thm.~\ref{thm:1:ER_stokes}, additional 
non-vanishing Galerkin orthogonality terms caused by the uni-directional 
coupling occur here, cf.~Rem.~\ref{rem:2:coupling}.
%%%%%%%%%%%%%%%%%%%%%%%%%%%%%%%%%%%%%%%%%%%%%%%%%%%%%%%%%%%%%%%%%%%%%%%%%%%%%%%
% Begin THEOREM 3.2 EE Transport
\begin{theorem}
\label{thm:2:ER_transport}
Let $\{\concentration,\dualz\}\in X \times X$,
$\{\concentration_{\tau},\dualz_{\tau}\}
\in
X_{\tau}^{r} \times X_{\tau}^{r}$,
and
$\{\concentration_{\tau h},\dualz_{\tau h}\}
\in X_{\tau h}^{r,p} \times X_{\tau h}^{r,p}$
be stationary points of
$\mathcal{L}, \mathcal{L}_{\tau}$, and $\mathcal{L}_{\tau h}$
on the different levels of discretization, i.e.,
\begin{equation*}
\begin{array}{r@{\;}c@{\;}l@{\;}}
\mathcal{L}^{\prime}(\concentration,\dualz;\convection)(\delta \concentration, \delta \dualz)
& = & 0 \quad
\forall \{\delta \concentration,\delta \dualz\}\in \mathcal{X} \times \mathcal{X}\,,
\\[1ex]
\mathcal{L}_{\tau}^{\prime}(\concentration_{\tau},\dualz_{\tau};\convection_{\sigma})
(\delta \concentration_{\tau}, \delta \dualz_{\tau})
& = & 0
\quad \forall \{\delta \concentration_{\tau},\delta \dualz_{\tau}\}
\in \mathcal{X}_{\tau}^{r} \times \mathcal{X}_{\tau}^{r}\,,
\\[1ex]
\mathcal{L}_{\tau h}^{\prime}(\concentration_{\tau h},\dualz_{\tau h};\convection_{\sigma h})
(\delta \concentration_{\tau h}, \delta \dualz_{\tau h})
& = & 0
\quad \forall \{\delta \concentration_{\tau h},\delta \dualz_{\tau h}\}
\in \mathcal{X}_{\tau h}^{r,p} \times \mathcal{X}_{\tau h}^{r,p}\,.
\end{array}
\end{equation*}
Then, for the discretization errors in space and time we get the representation 
formulas
\begin{subequations}
\label{eq:21:ER_transport}
\begin{equation}
\label{eq:21:ER_transport_time}
\begin{array}{r@{\;}c@{\;}l@{\;}}
J(\concentration)-J(\concentration_{\tau}) & = &
\frac{1}{2}\rho(\concentration_{\tau};\convection)(\dualz-\tilde{\dualz}_{\tau})
+ \frac{1}{2}\rho^{\ast}(\concentration_{\tau},\dualz_{\tau};\convection)
(\concentration-\tilde{\concentration}_{\tau})
\\[1.5ex]
& & 
+ \frac{1}{2} \mathcal{D}_{\tau}^{\prime}(\concentration_{\tau},\dualz_{\tau})
(\tilde{\concentration}_{\tau}-\concentration_{\tau},\tilde{\dualz}_{\tau}-\dualz_{\tau})
\\[1.5ex]
& & 
+ \mathcal{D}_{\tau}(\concentration_{\tau},\dualz_{\tau})
+ \mathcal{R}_{\tau}\,,
% \\[1.5ex]
% \hspace{-0.3cm}
\end{array}
\end{equation}
\begin{equation}
\label{eq:21:ER_transport_space}
\begin{array}{r@{\;}c@{\;}l@{\;}}
J(\concentration_{\tau})-J(\concentration_{\tau h}) & = &
\frac{1}{2}\rho_{\tau}(\concentration_{\tau h};\convection_{\sigma})(\dualz_{\tau}-\tilde{\dualz}_{\tau h})
+ \frac{1}{2}
\rho_{\tau}^{\ast}(\concentration_{\tau h},\dualz_{\tau h};\convection_{\sigma})
(\concentration_{\tau}-\tilde{\concentration}_{\tau h})
\\[1.5ex]
& & 
+ \frac{1}{2} \mathcal{D}_{\tau h}^{\prime}(\concentration_{\tau h},\dualz_{\tau h})
(\tilde{\concentration}_{\tau h}-\concentration_{\tau h},\tilde{\dualz}_{\tau h}-\dualz_{\tau h})
\\[1.5ex]
& & 
+ \mathcal{D}_{\tau h}(\concentration_{\tau h},\dualz_{\tau h}) + \mathcal{R}_{h}\,,
\end{array}
\end{equation}
\end{subequations}
with non-vanishing Galerkin orthogonality terms $\mathcal{D}_{\tau}(\cdot,\cdot)$ 
and $\mathcal{D}_{\tau h}(\cdot,\cdot)$, given by
\begin{equation}
\label{eq:22:Def_D_tau_h}
\begin{array}{r@{\;}c@{\;}l@{\;}}
\mathcal{D}_{\tau}(\concentration_{\tau},\dualz_{\tau}) 
&=&
\displaystyle\sum_{t_F \in \mathcal{F}_\tau}
(
[\concentration_\tau]_{t_F},\dualz_{\tau}(t_F^+))
\displaystyle
-\sum_{n=1}^{N^{\trans}}\sum_{K_n\in\mathcal{T}_{\tau,n}}\int_{K_n}\big(
(\convection-\convection_{\sigma})\cdot\nabla \concentration_\tau,\dualz_\tau
\big)\;\mathrm{d}t\,,
\\
\mathcal{D}_{\tau h}(\concentration_{\tau h},\dualz_{\tau h}) &=&
\displaystyle
S_A(\concentration_{\tau h}; \convection_{\sigma h})(\dualz_{\tau h})
\displaystyle
-\sum_{n=1}^{N^{\trans}}\sum_{K_n\in\mathcal{T}_{\tau,n}}\int_{K_n}\big(
(\convection_{\sigma}-\convection_{\sigma h})\cdot\nabla \concentration_{\tau h},
\dualz_{\tau h}
\big)\;\mathrm{d}t\,,
\end{array}
\end{equation}
where $\mathcal{D}^\prime_{\tau}(\cdot,\cdot)(\cdot,\cdot)$ and 
$\mathcal{D}^\prime_{\tau h}(\cdot,\cdot)(\cdot,\cdot)$ denote the G\^{a}teaux 
derivatives with respect to the first and second argument.
Here, $\{\tilde{\concentration}_{\tau},\tilde{\dualz}_{\tau}\}\in X_{\tau}^{r}
\times X_{\tau}^{r}$, and
$\{\tilde{\concentration}_{\tau h},\tilde{\dualz}_{\tau h}\} \in
X_{\tau h}^{r,p} \times X_{\tau h}^{r,p}$
can be chosen arbitrarily and the remainder terms $\mathcal{R}_{\tau},
\mathcal{R}_{h}$ are of higher-order with respect to the errors
$\concentration-\concentration_\tau,\dualz-\dualz_\tau$ and
$\concentration_{\tau}-\concentration_{\tau h},\dualz_{\tau}-\dualz_{\tau h}$,
respectively.
Furthermore, The explicit presentations of the primal and dual residuals based 
on the continuous and semi-discrete schemes $\rho, \rho^{\ast}$ and $\rho_\tau,\rho_\tau^\ast$,
respectively, are given in the Appendix.
\end{theorem}
%%%%%%%%%%%%%%%%%%%%%%%%%%%%%%%%%%%%%%%%%%%%%%%%%%%%%%%%%%%%%%%%%%%%%%%%%%%%%%%
% End THEOREM 3.2 EE Transport
\begin{mproof}
A detailed proof can be found in \cite[Ch.~5]{Bruchhaeuser2022a}, which is  
based on the idea given 
in \cite[Thm.~5.2]{Besier2012} stated for the nonstationary Navier-Stokes 
equations stabilized by local projection stabilization.
Basically, we are using a general result given in \cite[Lemma~5.1]{Besier2012} with 
the following settings:
\begin{equation*}
\begin{array}{r@{\;}c@{\;}l@{\;} @{\,\,}l @{\,\,}l @{\,}l @{\,\,}l @{\,\,}l @{\,}l @{\,\,}l @{\,\,}l @{\,}l @{\,\,}l @{\,\,}l @{\,}l @{\,\,}l}
L & = & \mathcal{L}\,,\;\;
& \tilde{L} & = \mathcal{L}_{\tau}\,,\;\;
& Y_1 & = \mathcal{X} \times \mathcal{X}\,,\;\;
& Y_2 & = \mathcal{X}_{\tau}^{r} \times \mathcal{X}_{\tau}^{r}\,,\;\;
& Y & = Y_1+Y_2
& \textnormal{ for \eqref{eq:21:ER_transport_time}: }
\\[1.5ex]
L & = & \mathcal{L}_\tau\,,\;\;
& \tilde{L} & = \mathcal{L}_{\tau h}\,,\;\;
& Y_1 & = \mathcal{X}_{\tau}^{r} \times \mathcal{X}_{\tau}^{r}\,,\;\;
& Y_2 & = \mathcal{X}_{\tau h}^{r,p} \times \mathcal{X}_{\tau h}^{r,p}\,,\;\;
& Y & = Y_1
& \textnormal{ for \eqref{eq:21:ER_transport_space}: }
\end{array}
\end{equation*}
where $\mathcal{L}, \mathcal{L}_{\tau}$ and $\mathcal{L}_{\tau h}$ are the 
Lagrangian functionals given by \eqref{eq:20:Lagrangian_transport} and 
$Y, Y_1$ and $Y_2$ are function spaces defined in 
\cite[Lemma~5.1]{Besier2012}. 
\end{mproof}

\section{Algorithm and Practical Aspects}
\label{sec:4:algorithm}

In this section, we illustrate some useful aspects for the practical implementation 
of the DWR-based error estimators derived in the section before.
We introduce localized forms of the error representations called error indicators
and explain how to compute them.
Moreover, we give insight into some implementational aspects with regard to our 
multirate approach.
Finally, we present the underlying cost-efficient adaptive space-time algorithm 
for the coupled flow and transport problem.
\subsection{Error Indicators and Approximation of Weights}
\label{sec:4:1:EIandWeights}
With regard to most application scenarios of the underlying model problem
of coupled flow and transport, for instance oil reservoir simulations or reactive
transport and degradation in the subsurface, the primary focus is to control the
transport problem under the condition that the influence of the error in the flow
problem stays comparatively small.
With this in mind, we focus here on creating most efficient, adaptively refined 
meshes in space and time for the transport problem using so-called weighted 
error indicators derived by means of the DWR approach.
Furthermore, as proclaimed at the beginning, we try to minimize numerical costs 
as far as possible. Thus, to reduce these costs significantly, but at the same 
time obtaining adaptive meshes as efficient as possible, we use here non-weighted, 
so-called auxiliary error indicators for the flow problem which avoid an explicit 
computation of the dual problem.

To use these error indicators locally within the adaptive mesh refinement 
process in space and time, we have to consider elementwise contributions of the
error representation formulas derived in Thm.~\ref{thm:1:ER_stokes} and 
Thm.~\ref{thm:2:ER_transport}, respectively
More precisely, the local error indicators $\eta_{\tau}^{\trans}$ and $\eta_{h}^{\trans}$ 
for the transport problem are obtained by neglecting the higher-order remainder 
terms $\mathcal{R}_{\tau}$ and $\mathcal{R}_{\tau h}$ in \eqref{eq:21:ER_transport}
and splitting the resulting quantities into elementwise contributions by means of 
the classical approach using integration by parts on every single mesh element,
cf.~\cite{Becker2001}.
\begin{equation}
\label{eq:24:eta_transport}
\begin{array}{r@{\;}c@{\;}l@{\;}}
J(u)-J(u_{\tau}) & \doteq &
\eta_{\tau}^{\trans}=\displaystyle \sum_{n=1}^{N^{\trans}} \eta_{\tau}^{\trans,n}
=\displaystyle \sum_{n=1}^{N^{\trans}}\sum_{K_n^{\trans}\in\mathcal{T}_{\tau,n}}\eta_{\tau,K_n^{\trans}}^{\trans}\,,
\\[1.5ex]
J(u_{\tau})-J(u_{\tau h}) & \doteq &
\eta_{h}^{\trans}=\displaystyle \sum_{n=1}^{N^{\trans}} \eta_{h}^{\trans,n}
=\displaystyle \sum_{n=1}^{N^{\trans}}
\sum_{K_n^{\trans}\in\mathcal{T}_{\tau,n}}
\sum_{K^{\trans}\in\mathcal{T}_{h}^{\trans,n}} \eta_{h,K^{\trans}}^{\trans,n}\,.
\end{array}
\end{equation}
Here, $\eta_{\tau,K_n^{\trans}}^{\trans}$ and $\eta_{h,K^{\trans}}^{\trans,n}$
denote the elementwise contributions on a temporal and spatial mesh cell 
$K_n^{\trans}$ and $K^{\trans}$, respectively, where the contribution on a single
temporal cell is computed by collecting the contributions of all spatial cells
of the spatial triangulation $\mathcal{T}_{h}^{\trans,n}$ corresponding to this 
temporal cell, more precisely $\eta_{\tau,K_n^{\trans}}^{\trans}:=
\sum_{K^{\trans}\in\mathcal{T}_{h}^{\trans,n}} \eta_{h,K^{\trans}}^{\trans,n}\,.$
For an explicit presentation of these elementwise indicators and further details,
we refer to \cite[Ch.~5.3.2.2]{Bruchhaeuser2022a}.
The non-weighted, auxiliary error indicators $\tilde{\eta}_{\sigma h}^{\flow}=
\tilde{\eta}_{\sigma}^{\flow}+\tilde{\eta}_{h}^{\flow}$ are given by
\begin{equation}
\label{eq:24:auxiliary_eta_stokes}
\begin{array}{r@{\;}c@{\;}l@{\;}}
J(\mathbf{u})-J(\mathbf{u}_{\sigma}) & \approx &
\tilde{\eta}_{\sigma}^{\flow}=\displaystyle \sum_{n=1}^{N^{\flow}} \tilde{\eta}_{\sigma}^{\flow,n}
=\displaystyle \sum_{n=1}^{N^{\flow}}\sum_{K_n^{\flow}\in\mathcal{T}_{\sigma,n}}\tilde{\eta}_{\sigma,K_n^{\flow}}^{\flow}\,,
\\[1.5ex]
J(\mathbf{u}_{\sigma})-J(\mathbf{u}_{\sigma h}) & \approx &
\tilde{\eta}_{h}^{\flow}=\displaystyle \sum_{n=1}^{N^{\flow}} \tilde{\eta}_{h}^{\flow,n}
=\displaystyle \sum_{n=1}^{N^{\flow}}
\sum_{K_n^{\flow}\in\mathcal{T}_{\sigma,n}}
\sum_{K^{\flow}\in\mathcal{T}_{h}^{\flow,n}} \tilde{\eta}_{h,K^{\flow}}^{\flow,n}\,,
\end{array}
\end{equation}
where, in contrast to the classical DWR-based derived indicators above, the 
local errror indicators $\tilde{\eta}_{h,K^{\flow}}^{\flow,n}$ are based on the 
so-called Kelly Error Estimator, cf.~\cite{Kelly1983} as well as the reference 
documentation of the \texttt{deal.II} library \cite{dealiiReference93} for more 
details.

To compute these error indicators, we replace all unknown solutions by the 
approximated fully discrete solutions $u_{\tau h} \in \mathcal{X}_{\tau h}^{r,p}$,
$z_{\tau h} \in \mathcal{X}_{\tau h}^{r,q}$, with $p < q$, and 
$\mathbf{u}_{\sigma h}=\{\mathbf{v}_{\sigma h},p_{\sigma h}\} \in 
\mathcal{Y}_{\sigma h}^{0,p_{\mathbf{v}},p_{p}}$, $p_{p}+1=p_{\mathbf{v}} \geq 2$.
More precisely, we restrict the solution $\{\mathbf{v}_{\sigma h},p_{\sigma h}\}$
of the flow problem on each $I_n^{\flow}$ to a piecewise constant discontinuous
Galerkin (dG($0$)) time approximation. This is due to simplicity reasons for
the implementation of the solution transfer of the flow field $\mathbf{v}_{\sigma h}$
required within the transport problem \eqref{eq:17:fully_transport}, which is 
described in detail in the following section.
Moreover, the temporal and spatial weights as well as the flow field differences 
arising within the DWR-based error indicators \eqref{eq:24:eta_transport} of the 
transport problem are approximated in the following way:
\begin{itemize}
\itemsep1.5ex

\item Approximate the temporal weights $u-\tilde{u}_{\tau}, z-\tilde{z}_{\tau}$
by means of a higher-order reconstruction using Gauss-Lobatto quadrature points,
exemplary given by
\begin{displaymath}
\begin{array}{r@{\,}c@{\,}l}
u-\tilde{u}_{\tau} & \approx & \mathrm{E}_{\tau}^{(r+1)}u_{\tau h}-u_{\tau h}\,,
\end{array}
\end{displaymath}
using an reconstruction in time operator $\mathrm{E}_{\tau}^{(r+1)}$ 
thats acts on a time cell $K_n$ of length $\tau_K$ and lifts the solution to a piecewise 
polynomial of degree ($r$+$1$) in time, cf. Fig.~\ref{fig:1:hoExtrapolationGL}.
%
%%%%%%%%%%%%%%%%%%%%%%%%%%%%%%%%%%%%%%%%%%%%%%%%%%%%%%%%%%%%%%%%%%%%%%%%%%%%%%%%
%% Begin Figure ho Extrapolation Gauss Lobatto
\begin{figure}[hbt!]
\centering
\vskip1.5ex
\resizebox{0.95\linewidth}{!}{%
\begin{tikzpicture}
\tikzstyle{ns1}=[line width=0.4]
\tikzstyle{ns2}=[line width=0.15,opacity=.5]
\tikzstyle{ns3}=[line width=0.6]
\tikzstyle{ns4}=[line width=0.9]
%%%%%%%%%%%%%%%%%%%%%%%%%%%%%%%%%%%%%%%%%%%%%%%%%%%%%%%%%%%%%%%%%%%%%%%%%%%%%%%%
%% dG(0)
% Timescale
\draw[->,ns3] (-8.5,-0.5) -- (-3.5,-0.5);
\draw[ns1] (-8,-0.65) -- (-8,-0.35);
\draw[ns1] (-6,-0.65) -- (-6,-0.35);
\draw[ns1] (-4,-0.65) -- (-4,-0.35);
\node[] at (-8,-1.0) {$t_{n-1}$};
\node[] at (-7,-0.2) {$K_{n_{\ell}}$};
\node[] at (-6,-1.0) {$t_{\mathcal{F}, n_{\ell}}$};
\node[] at (-5,-0.2) {$K_n$};
\node[] at (-4,-1.0) {$t_{n}$};
% Values
\draw[domain=-8:-6] plot (\x, {0.8});
\draw[fill=black](-7,0.8)circle(1.5pt);
\draw[HSUred,fill=HSUred](-6,0.8)circle(1.5pt);
\draw[domain=-6:-4] plot (\x, {1.5});
\draw[fill=black](-5,1.5)circle(1.5pt);
\draw[HSUred,fill=HSUred](-4,1.5)circle(1.5pt);
\draw[HSUred,ns4,densely dashed,domain=-6:-4] plot (\x, {0.35*\x+2.9});
% Functions
\node[] at (-3.5,1.5) {$v_{\tau,n}^{0}$};
\node[] at (-3.2,0.8) {\textcolor{HSUred}{$\operatorname{E}_{\tau}^{1}v_{\tau,n}^{0}$}};
%%%%%%%%%%%%%%%%%%%%%%%%%%%%%%%%%%%%%%%%%%%%%%%%%%%%%%%%%%%%%%%%%%%%%%%%%%%%%%%%
%% dG(1)
% Timescale
\draw[->,ns3] (-2.5,-0.5) -- (2.5,-0.5);
\draw[ns1] (-2,-0.65) -- (-2,-0.35);
\draw[ns1] (0,-0.65) -- (0,-0.35);
\draw[ns1] (2,-0.65) -- (2,-0.35);
\node[] at (-2,-1) {$t_{n-1}$};
\node[] at (-1,-0.2) {$K_{n_{\ell}}$};
\node[] at (0,-1) {$t_{\mathcal{F}, n_{\ell}}$};
\node[] at (1,-0.2) {$K_n$};
\node[] at (2,-1) {$t_{n}$};
% Values
\draw[domain=-2:0] plot (\x, {-.4*\x+.8});
\draw[fill=black](-1.5773,1.43)circle(1.5pt);
\draw[fill=black](-0.4227,0.97)circle(1.5pt);
\draw[HSUred,fill=HSUred](0,0.8)circle(1.5pt);
\draw[domain=0:2] plot (\x, {.4*\x+1.5});
\draw[fill=black](0.4227,1.66)circle(1.5pt);
\draw[HSUred,fill=HSUred](1,1.9)circle(1.5pt);
\draw[fill=black](1.5773,2.13)circle(1.5pt);
\draw[HSUred,fill=HSUred](2,2.3)circle(1.5pt);
\draw[HSUred,ns4,densely dashed,domain=0:2] plot (\x, {-0.35*\x*\x+1.45*\x+0.8});
% Functions
\node[] at (2.5,2.3) {$v_{\tau,n}^{1}$};
\node[] at (2.8,1.6) {\textcolor{HSUred}{$\operatorname{E}_{\tau}^{2}v_{\tau,n}^{1}$}};
%%%%%%%%%%%%%%%%%%%%%%%%%%%%%%%%%%%%%%%%%%%%%%%%%%%%%%%%%%%%%%%%%%%%%%%%%%%%%%%%
%% dG(2)
% Timescale
\draw[->,ns3] (3.5,-0.5) -- (8.5,-0.5);
\draw[ns1] (4,-0.65) -- (4,-0.35);
\draw[ns1] (6,-0.65) -- (6,-0.35);
\draw[ns1] (8,-0.65) -- (8,-0.35);
\node[] at (4,-1) {$t_{n-1}$};
\node[] at (5,-0.2) {$K_{n_{\ell}}$};
\node[] at (6,-1) {$t_{\mathcal{F}, n_{\ell}}$};
\node[] at (7,-0.2) {$K_n$};
\node[] at (8,-1) {$t_{n}$};
% Values
\draw[domain=4:6] plot (\x, {1.0782*\x*\x-10.7817*\x+27.7542});
\draw[fill=black](4.254,1.4)circle(1.5pt);
\draw[fill=black](5,0.8)circle(1.5pt);
\draw[fill=black](5.746,1.4)circle(1.5pt);
\draw[HSUred,fill=HSUred](6,1.8792)circle(1.5pt);
\draw[domain=6:8] plot (\x, {-1.0782*\x*\x+15.0943*\x-49.2302});
\draw[fill=black](6.254,3.0)circle(1.5pt);
\draw[HSUred,fill=HSUred](6.5528,3.3827)circle(1.5pt);
\draw[fill=black](7,3.6)circle(1.5pt);
\draw[HSUred,fill=HSUred](7.4472,3.3822)circle(1.5pt);
\draw[fill=black](7.746,3.0)circle(1.5pt);
\draw[HSUred,fill=HSUred](8,2.5194)circle(1.5pt);
\draw[HSUred,ns4,densely dashed,domain=6:8] plot (\x, {0.4001*\x*\x*\x-9.8808*\x*\x+79.4344*\x-205.4423});
% Functions
\node[] at (8.5,2.5) {$v_{\tau,n}^{2}$};
\node[] at (8.8,1.8) {\textcolor{HSUred}{$\operatorname{E}_{\tau}^{3}v_{\tau,n}^{2}$}};
% dots
\node[] at (9.5,-0.5) {$\dots$};
\end{tikzpicture}
}
\caption{Reconstruction of a discontinuous constant (left), linear (middle) and 
quadratic (right) in time function on an exemplary cell in time $K_n$ using 
Gauss-Lobatto quadrature points.}
\label{fig:1:hoExtrapolationGL}
\end{figure}
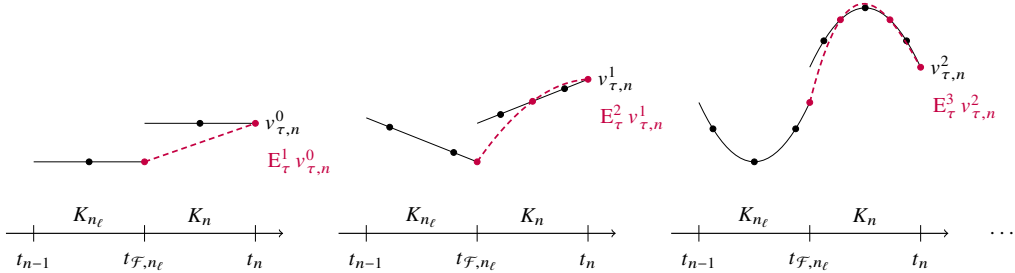
%%%%%%%%%%%%%%%%%%%%%%%%%%%%%%%%%%%%%%%%%%%%%%%%%%%%%%%%%%%%%%%%%%%%%%%%%%%%%%%%
%% End Figure ho Extrapolation

We point out that within this approximation strategy the respective dual problem 
is solved in the same finite element space $\mathcal{X}_{\tau}^{r}$ as used for
the primal problem.
This approximation technique is done for the purpose to further reduce numerical 
costs solving the dual transport problem compared to a higher-order finite 
element approximation that is used, for instance, in \cite[Sec. 4]{Bause2021}.
\item Approximate the spatial weights $u_{\tau}-\tilde{u}_{\tau h}$ and 
$z_{\tau}-\tilde{z}_{\tau h}$ by means of a patch-wise higher-order interpolation
and a higher-order finite elements approach, respectively, given by
\begin{displaymath}
\begin{array}{r@{\,}c@{\,}l}
u_{\tau}-\tilde{u}_{\tau h} & \approx & \mathrm{I}_{2h}^{(2p)}u_{\tau h}-u_{\tau h}\,,
\\[1.0ex]
z_{\tau}-\tilde{z}_{\tau h} & \approx & z_{\tau h}-\mathrm{R}_{h}^{p}z_{\tau h}\,,
\end{array}
\end{displaymath}
using an interpolation in space operator $\mathrm{I}_{2h}^{(2p)}$ and an
restriction in space operator $\mathrm{R}_{h}^{p}$ that are described in 
detail in our work \cite[Sec. 4]{Bause2021}.
\item Approximate the temporal flow field difference $\mathbf{v}-\mathbf{v}_{\sigma}$ 
by means of a higher-order extrapolation using Gauss-Lobatto quadrature points 
given by
\begin{displaymath}
\mathbf{v}-\mathbf{v}_{\sigma} \approx 
\boldsymbol{\operatorname{E}}_{\tau}^{(r+1)}\mathbf{v}_{\sigma}-\mathbf{v}_{\sigma}\,,
\end{displaymath}
with $\boldsymbol{\operatorname{E}}_{\tau}^{(r+1)}$ acting componentwise like 
$\operatorname{E}_{\tau}^{(r+1)}$.
\item 
Approximate the spatial flow field difference $\mathbf{v}_{\sigma}-\mathbf{v}_{\sigma h}$ 
by means of a patch-wise higher-order interpolation given by
\begin{displaymath}
\mathbf{v}_{\sigma}-\mathbf{v}_{\sigma h} \approx 
\boldsymbol{\operatorname{I}}_{2h}^{(2p)}\mathbf{v}_{\sigma h}-\mathbf{v}_{\sigma h}\,,
\end{displaymath}
with $\boldsymbol{\operatorname{I}}_{2h}^{(2p)}$ acting componentwise like 
$\operatorname{I}_{2h}^{(2p)}$.
Note that the fully discrete solution $\mathbf{v}_{\sigma h}$ has to be 
transferred to the spatial mesh used for the transport problem. We call this 
process a solution mesh transfer that will be explained in greater detail in 
Sec.~\ref{sec:4:2:MultirateAspects}.
\end{itemize}

\subsection{Implementation of Multirate Aspects}
\label{sec:4:2:MultirateAspects}

Here, we present some implementational aspects with regard to our multirate
approach for coupled flow and transport problems.
With regard to a practical realization of this approach, the following 
three aspects are of particular importance, being specified subsequently: % and will be described below:
\begin{itemize}
\item Initialization of spatial and temporal meshes for both flow and transport
 problem: Space-time slabs.
\item Interaction of spatial and temporal meshes between flow and transport 
 problem: Solution mesh transfer.
\item Realization of adaptive mesh refinement in space and time: Involvement of slabs.
\end{itemize}
In order to keep this work short, we only address the main parts of these aspects
and refer to \cite{Bruchhaeuser2022,Bruchhaeuser2022a} for a detailed explanation.
\subsubsection*{Initialization of Spatial and Temporal Meshes for Flow and Transport}
For an adaptive numerical approximation of the coupled flow and transport problem
\eqref{eq:1:stokes_problem}, \eqref{eq:2:transport_problem}, the space-time domain
$Q = \Omega \times I$ is divided into non-overlapping space-time slabs 
$Q_{n}^{\flow}=\mathcal{T}_{h,n}^{\flow}\times \mathcal{T}_{\sigma,n},
n=1,\dots,N^{\flow}\,,$ as well as $Q_{n}^{\trans}=\mathcal{T}_{h,n}^{\trans}
\times \mathcal{T}_{\tau,n}, n=1,\dots,N^{\trans}\,,$ with $N^{\flow}\leq N^{\trans}\,,$ 
for the flow and transport problem, respectively.
On such a slab, a tensor-product of a $d$-dimensional, $d=1,2,3\,,$ spatial 
finite element space with a one-dimensional temporal finite element space is 
implemented.
An exemplary illustration of such slabs is given by Fig.~\ref{fig:4:SolutionTransferStokes}. 
The temporal finite element space is based on a discontinuous Galerkin dG($r$) 
method of arbitrary order $r\,, r \geq 0$, whereas the spatial finite element 
space is based on a continuous Galerkin cG($p$) method of arbitrary order 
$p\,, p \geq 1$. 
Thus, we are using here and in the following the notation cG($p$)-dG($r$)
method.

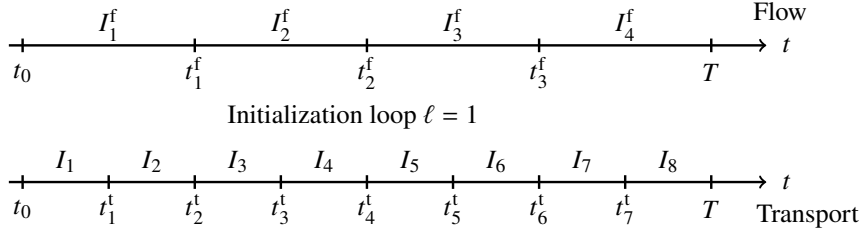
\begin{figure}[hbt!]
\centering
\resizebox{0.8\linewidth}{!}{%
\begin{tikzpicture}
\tikzstyle{ns1}=[line width=1.]

% temporal double axis with captions
\node at (11.6,-1.5) {\normalsize Transport};
\node at (11.3,1.0) {\normalsize $t$};
\node at (11.3,-1.0) {\normalsize $t$};
\draw[->,ns1] (0,1.0) -- (11,1.0);

\draw[->,ns1] (0,-1) -- (11,-1);
\node at (11.2,1.5) {Flow};

\node at (5.0,-.05) {Initialization loop $\ell=1$};

% time points: Transport
\node at (0.2, -1.4) {\normalsize $t_0$};
\draw[ns1] (0.2, -1.1) -- (0.2, -0.9);

\node at (0.825, -.7) {\normalsize $I_1$};

\node at (1.45, -1.4) {\normalsize $t_1^\trans$};
\draw[ns1] (1.45, -1.1) -- (1.45, -0.9);

\node at (2.075, -.7) {\normalsize $I_2$};

\node at (2.7, -1.4) {\normalsize $t_2^\trans$};
\draw[ns1] (2.7, -1.1) -- (2.7, -0.9);

\node at (3.325, -.7) {\normalsize $I_3$};

\node at (3.95, -1.4) {\normalsize $t_3^\trans$};
\draw[ns1] (3.95, -1.1) -- (3.95, -0.9);

\node at (4.575, -.7) {\normalsize $I_4$};

\node at (5.2, -1.4) {\normalsize $t_4^\trans$};
\draw[ns1] (5.2, -1.1) -- (5.2, -0.9);

\node at (5.825, -.7) {\normalsize $I_5$};

\node at (6.45, -1.4) {\normalsize $t_5^\trans$};
\draw[ns1] (6.45, -1.1) -- (6.45, -0.9);

\node at (7.075, -.7) {\normalsize $I_6$};

\node at (7.7, -1.4) {\normalsize $t_6^\trans$};
\draw[ns1] (7.7, -1.1) -- (7.7, -0.9);

\node at (8.325, -.7) {\normalsize $I_7$};

\node at (8.95, -1.4) {\normalsize $t_7^\trans$};
\draw[ns1] (8.95, -1.1) -- (8.95, -0.9);

\node at (9.575, -.7) {\normalsize $I_8$};

\node at (10.2, -1.4) {\normalsize $T$};
\draw[ns1] (10.2, -1.1) -- (10.2, -0.9);

% time points: Stokes flow
\draw[ns1] (0.2, 0.9) -- (0.2, 1.1);
\node at (0.2, 0.6) {\normalsize $t_0$};

\node at (1.45, 1.3) {\normalsize $I_1^{\flow}$};

\draw[ns1] (2.7, 0.9) -- (2.7, 1.1);
\node at (2.7, 0.6) {\normalsize $t_1^\flow$};

\node at (3.95, 1.3) {\normalsize $I_2^{\flow}$};

\draw[ns1] (5.2, 0.9) -- (5.2, 1.1);
\node at (5.2, 0.6) {\normalsize $t_2^\flow$};

\node at (6.45, 1.3) {\normalsize $I_3^{\flow}$};

\draw[ns1] (7.7, 0.9) -- (7.7, 1.1);
\node at (7.7, 0.6) {\normalsize $t_3^\flow$};

\node at (8.95, 1.3) {\normalsize $I_4^{\flow}$};

\draw[ns1] (10.2, 0.9) -- (10.2, 1.1);
\node at (10.2, 0.6) {\normalsize $T$};

\end{tikzpicture}}

\caption{Exemplary initialization of different temporal meshes for flow and
transport.}
\label{fig:2:InitializationTemporalMeshes}
\end{figure}

As mentioned at the beginning of this work, the behavior of the underlying flow
and transport problem is rather contrary with regard to the processes that take 
place in time.
Due to these contrary behaviors, measured by means of the introduced characteristic times
\eqref{eq:3:characteristic-times}, the flow and transport problem are initialized
independently on different time scales fulfilling the following conditions, 
cf.~Fig.~\ref{fig:2:InitializationTemporalMeshes}:
\begin{itemize}
\item The temporal mesh of the flow problem is coarser or equal to that of 
the transport problem.
\item The endpoints of the temporal mesh of the flow problem must match with 
endpoints of the temporal mesh of the transport problem.
\end{itemize}
% 
%%%%%%%%%%%%%%%%%%%%%%%%%%%%%%%%%%%%%%%%%%%%%%%%%%%%%%%%%%%%%%%%%%%%%%%%%%%%%%%%
%% Begin Figure 1 Different Time Scales
%%%%%%%%%%%%%%%%%%%%%%%%%%%%%%%%%%%%%%%%%%%%%%%%%%%%%%%%%%%%%%%%%%%%%%%%%%%%%%%%
%% End Figure 1 Different Time Scales
%

\begin{figure}[h!]
\centering
% \vskip1.5ex
\resizebox{0.5\linewidth}{!}{%
\begin{tikzpicture}
% T_h
\draw (0,0) -- (4,0) -- (4,4) -- (0,4) -- (0,0);
\draw (0,2) -- (4,2);
\draw (2,0) -- (2,4);
\draw (2,1) -- (4,1);
\draw (3,0) -- (3,2);
\draw[fill=black](2,1)circle(2pt);
\draw[fill=black](3,2)circle(2pt);
\draw (2,-0.5) node[]{Flow $\mathcal{T}_{h,n}^{\flow}$};
% T_2h
\draw (5,0) -- (9,0) -- (9,4) -- (5,4) -- (5,0);
\draw (5,1) -- (9,1);
\draw (5,2) -- (9,2);
\draw (5,3) -- (9,3);
\draw (6,0) -- (6,4);
\draw (7,0) -- (7,4);
\draw (8,0) -- (8,4);
\draw (7.5,0) -- (7.5,2);
\draw (8.5,0) -- (8.5,2);
\draw (7,0.5) -- (9,0.5);
\draw (7,1.5) -- (9,1.5);
\draw[fill=black](7,0.5)circle(2pt);
\draw[fill=black](7,1.5)circle(2pt);
\draw[fill=black](7.5,2)circle(2pt);
\draw[fill=black](8.5,2)circle(2pt);
\draw (7,-0.5) node[]{Transport $\mathcal{T}_{h,n}^{\trans}$};
\end{tikzpicture}
}
\caption{Exemplary initialization of different spatial meshes for the flow
($\mathcal{T}_{h,n}^{\flow}$) and transport ($\mathcal{T}_{h,n}^{\trans}$) problem
organized in a patch-wise manner.}
\label{fig:3:InitializationSpatialMeshes}
\end{figure}
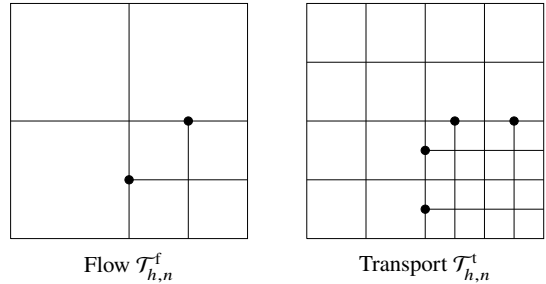

With regard to the underlying spatial triangulations, we state the following.
We assume the spatial triangulations to be regular and organized in a patch-wise
manner, but allowing hanging nodes with regard to adaptive refinements, 
cf.~Fig.~\ref{fig:3:InitializationSpatialMeshes}.
We point out that the global conformity of the finite element approach is 
preserved since the unknowns at such hanging nodes are eliminated by 
interpolation between the neighboring 'regular' nodes; 
cf.~\cite{Carey1984,Bangerth2003}.
For the sake of implementational simplicity, we allow the spatial meshes to change 
between two consecutive slabs, but to be equal on all degrees of freedom in time 
used within one slab, cf.~Fig.~\ref{fig:4:SolutionTransferStokes}.
This approach is referred to as the concept of dynamic meshes, cf., e.g., \cite{Schmich2008}.
Thus, for the initialization of the spatial meshes for the flow and transport
problem, we assume the following:
\begin{itemize}
\item The spatial mesh of the flow problem is coarser or equal to that of the 
transport problem.
\item If the flow spatial mesh is coarser, the transport spatial mesh 
has to be emerged from the flow mesh in the sense of a patch-wise manner, i.e.
$\mathcal{T}_{h,n}^{\trans}$ is obtained by uniform refinement of the coarser 
decomposition $\mathcal{T}_{h,n}^{\flow}$ such that it is always possible to 
combine four ($d=2$) or eight ($d=3$) adjacent elements of $\mathcal{T}_{h,n}^{\trans}$ 
to obtain one element of $\mathcal{T}_{h,n}^{\flow}$, cf.~Fig.~\ref{fig:3:InitializationSpatialMeshes}.
\end{itemize}
% 
%%%%%%%%%%%%%%%%%%%%%%%%%%%%%%%%%%%%%%%%%%%%%%%%%%%%%%%%%%%%%%%%%%%%%%%%%%%%%%%%
%% Begin Figure Initialization Spatial Mesh
%%%%%%%%%%%%%%%%%%%%%%%%%%%%%%%%%%%%%%%%%%%%%%%%%%%%%%%%%%%%%%%%%%%%%%%%%%%%%%%%
%% End Figure Initialization Spatial Mesh

\subsubsection*{Interaction of Spatial and Temporal Meshes between Transport 
and Flow}

As outlined in Rem.~\ref{rem:2:coupling}, the coupling within the underlying
model problem \eqref{eq:1:stokes_problem}, \eqref{eq:2:transport_problem} is 
given via the convection variable $\mathbf{v}$. 
Thus, we need the fully discrete flow field solution $\mathbf{v}_{\sigma h}$ for the
numerical approximation of the stabilized primal and dual transport problems given
by \eqref{eq:17:fully_transport} and \eqref{eq:?:DualProblems}, respectively.
More precisely, the fully discrete flow field solution $\mathbf{v}_{\sigma h}$ has to
be transferred to the respective transport meshes in an appropriate way.
This solution mesh transfer is handled in a similar fashion as used for the 
transfer of a solution between two consecutive slabs with different underlying 
spatial meshes and will be specified in the following. 
%
%%%%%%%%%%%%%%%%%%%%%%%%%%%%%%%%%%%%%%%%%%%%%%%%%%%%%%%%%%%%%%%%%%%%%%%%%%%%%%%%
%% Begin Figure 3 Solution Transfer Stokes
\begin{figure}[hbt!]
\centering
\resizebox{.8\linewidth}{!}{%
\begin{tikzpicture}
\tikzstyle{ns1}=[line width=1.]
\tikzstyle{ns2}=[line width=1.,opacity=.5]
\tikzstyle{ns3}=[line width=1.,opacity=.5,dashed]
\tikzstyle{ns4}=[line width=0.9]

%%%%%%%%%%%%%%%%%%%%%%%%%%%%%%%%%%%%%%%%%%%%%%%%%%%%%%%%%%%%%%%%%%%%%%%%%%%%%%%%
% STOKES 
% Q_n
% Q_n label
\draw[very thick,HSUred,decorate,decoration={brace,amplitude=14pt}]
(-4.2,19) -- (17.8,19) node[midway, above,yshift=14pt,]{\LARGE $Q_{n}^{\flow}$};

% Front
\draw[ns1] (-4.2,15) -- (-1.5,14);
\draw[ns1] (-1.5,14) -- (-1.5,18);
\draw[ns1] (-1.5,18) -- (-4.2,19);
\draw[ns1] (-4.2,19) -- (-4.2,15);

% Connection Lines Front<-->Back
\draw[ns1] (-1.5,14) -- (20.5,14);
\draw[ns1, dashed] (-4.2,15) -- (17.8,15);
\draw[ns1] (-1.5,18) -- (20.5,18);
\draw[ns1] (-4.2,19) -- (17.8,19);

% Middle
\draw[ns1] (6.8,15) -- (9.5,14);
\draw[ns1] (9.5,14) -- (9.5,18);
\draw[ns1] (9.5,18) -- (6.8,19);
\draw[ns1] (6.8,19) -- (6.8,15);
% Inner mesh middle
\draw[ns1] (8.15,14.5) -- (8.15,18.5);
\draw[ns1] (6.8,17) -- (9.5,16);
\draw[ns1] (8.825,14.25) -- (8.825,16.25);
\draw[ns1] (8.15,15.5) -- (9.5,15);

% Back
\draw[ns1,dashed] (17.8,15) -- (20.5,14);
\draw[ns1] (20.5,14) -- (20.5,18);
\draw[ns1] (20.5,18) -- (17.8,19);
\draw[ns1,dashed] (17.8,19) -- (17.8,15);

% t axis
\node[] at (21,12) {\LARGE $t$};
\draw[->,ns1] (-4.2,12) -- (20.5,12);
\node[] at (-2.85,11.3) {\LARGE $t_{n-1}^{\flow}$};
\draw[HSUred,fill=HSUred](-2.85,14.5)circle(3.5pt);
\draw[->,ns1, HSUred, dotted] (-2.85,14.5) -- (-2.85,12.2);
\draw[HSUred,fill=HSUred](-2.85,12)circle(3.5pt);
\node[] at (8.15,11.3) {\LARGE $t_{\textnormal{DoF}}^{\flow,1}$};
\draw[HSUred,fill=HSUred](8.15,14.5)circle(3.5pt);
\draw[->,ns1, HSUred, dotted] (8.15,14.5) -- (8.15,12.2);
\draw[HSUred,fill=HSUred](8.15,12)circle(3.5pt);
\node[] at (19.15,11.3) {\LARGE $t_{n}^{\flow}$};
\draw[HSUred,fill=HSUred](19.15,14.5)circle(3.5pt);
\draw[->,ns1, HSUred, dotted] (19.15,14.5) -- (19.15,12.2);
\draw[HSUred,fill=HSUred](19.15,12)circle(3.5pt);

% Arrows Interpolation
%Qn
\draw[->,ns1, HSUred, dotted] (8.15,14.5) -- (-1.25,0.65);
\node[HSUred] at (3.45,8.5) {\LARGE $\operatorname{I}_{h}^{n}$};
\draw[->,ns1, HSUred, dotted] (8.15,14.5) -- (2.7,0.65);
\node[HSUred] at (5.25,8.5) {\LARGE $\operatorname{I}_{h}^{n}$};
\draw[->,ns1, HSUred, dotted] (8.15,14.5) -- (6.66,0.68);
\node[HSUred] at (7.1,8.5) {\LARGE $\operatorname{I}_{h}^{n}$};
%Qn+1
\draw[->,ns1, HSUred, dotted] (8.15,14.5) -- (9.64,0.68);
\node[HSUred] at (9.55,8.5) {\LARGE $\operatorname{I}_{h}^{n+1}$};
\draw[->,ns1, HSUred, dotted] (8.15,14.5) -- (13.6,0.65);
\node[HSUred] at (11.3,8.5) {\LARGE $\operatorname{I}_{h}^{n+1}$};
\draw[->,ns1, HSUred, dotted] (8.15,14.5) -- (17.55,0.65);
\node[HSUred] at (13.1,8.5) {\LARGE $\operatorname{I}_{h}^{n+1}$};

%%%%%%%%%%%%%%%%%%%%%%%%%%%%%%%%%%%%%%%%%%%%%%%%%%%%%%%%%%%%%%%%%%%%%%%%%%%%%%%%
% TRANSPORT
% Q_n
% 
% First adaptive spatial mesh (on t^1_DoF of Q_n)
%
% Mesh boundary
\draw[ns1] (-2.7,1) -- (0,0);
\draw[ns1] (0,0) -- (0,4);
\draw[ns1] (0,4) -- (-2.7,5);
\draw[ns1] (-2.7,5) -- (-2.7,1);
%
% Inner Mesh
\draw[ns1] (-2.7,3) -- (0,2);
\draw[ns1] (-1.35,0.5) -- (-1.35,4.5);
\draw[ns1] (-2.7,2) -- (-1.35,1.5);
\draw[ns1] (-2.025,0.75) -- (-2.025,2.75);
\draw[ns1] (-2.025,2.25) -- (-1.35,2.);
\draw[ns1] (-1.6875,1.625) -- (-1.6875,2.625);
\draw[ns1] (-2.025,2.) -- (-1.6875,1.875);
\draw[ns1] (-1.8562,1.7) -- (-1.8562,2.2);
\draw[ns1] (-1.35,3.5) -- (0,3.);
\draw[ns1] (-0.65,2.25) -- (-0.65,4.25);

%\visible<2->{
\node[] at (-1.35,-2.5) {\LARGE $t_{\textnormal{DoF}}^{\trans, 1}$};
\draw[HSUred,fill=HSUred](-1.35,0.5)circle(3.5pt);
\draw[->,ns1, HSUred, dotted] (-1.35,0.5) -- (-1.35,-1.8);
\draw[HSUred,fill=HSUred](-1.35,-2)circle(3.5pt);

\node[] at (-2.85,-2.5) {\LARGE $t_{n-1}^{\trans}$};
\draw[HSUred,fill=HSUred](-2.85,0.5)circle(3.5pt);
\draw[->,ns1, HSUred, dotted] (-2.85,0.5) -- (-2.85,-1.8);
\draw[HSUred,fill=HSUred](-2.85,-2)circle(3.5pt);

\node[] at (8.15,-2.5) {\LARGE $t_{n}^{\trans}$};
\draw[HSUred,fill=HSUred](8.15,0.5)circle(3.5pt);
\draw[->,ns1, HSUred, dotted] (8.15,0.5) -- (8.15,-1.8);
\draw[HSUred,fill=HSUred](8.15,-2)circle(3.5pt);

% Q_n label
\draw[very thick,HSUred,decorate,decoration={brace,amplitude=14pt}]
(-4.2,5) -- (6.8,5) node[midway, above,yshift=14pt,]{\LARGE $Q_{n}^{\trans}$};

% Front extension t_n-1
\draw[ns1] (-4.2,1) -- (-1.5,0);
\draw[ns1] (-1.5,0) -- (-1.5,4);
\draw[ns1] (-1.5,4) -- (-4.2,5);
\draw[ns1] (-4.2,5) -- (-4.2,1);

\draw[ns1] (-4.2,5) -- (-2.7,5);
\draw[ns1] (-1.5,4) -- (0,4);
\draw[ns1, dashed] (-4.2,1) -- (-2.7,1);
\draw[ns1] (-1.5,0) -- (0,0);

\draw[ns1] (0,0) -- (8,0);
\draw[ns1] (8,4) -- (0,4);
\draw[ns1] (-2.7,5) -- (5.3,5); %Frontedge
\draw[ns1, dashed] (-2.7,1) -- (5.3,1); % Backedge

% Middle extension t_n
\draw[ns1,dashed] (6.8,1) -- (9.5,0);
\draw[ns1] (9.5,0) -- (9.5,4);
\draw[ns1] (9.5,4) -- (6.8,5);
\draw[ns1,dashed] (6.8,5) -- (6.8,1);

\draw[ns1] (5.3,5) -- (6.8,5);
\draw[ns1] (8,4) -- (9.5,4);
\draw[ns1, dashed] (5.3,1) -- (6.8,1);
\draw[ns1] (8,0) -- (9.5,0);
%}

%\visible<3->{
\node[] at (2.65,-2.5) {\LARGE $t_{\textnormal{DoF}}^{\trans, 2}$};
\draw[HSUred,fill=HSUred](2.65,0.5)circle(3.5pt);
\draw[->,ns1, HSUred, dotted] (2.65,0.5) -- (2.65,-1.8);
\draw[HSUred,fill=HSUred](2.65,-2)circle(3.5pt);

% Second adaptive spatial mesh (on t^2_DoF of Q_n)
%
% boundary
\draw[ns2] (1.3,1) -- (4,0);
\draw[ns2] (4,0) -- (4,4);
\draw[ns2] (4,4) -- (1.3,5);
\draw[ns2] (1.3,5) -- (1.3,1);
% Inner Mesh
\draw[ns2] (1.3,3) -- (4,2);
\draw[ns2] (2.65,0.5) -- (2.65,4.5);
\draw[ns2] (1.3,2) -- (2.65,1.5);
\draw[ns2] (1.975,0.75) -- (1.975,2.75);
\draw[ns2] (1.975,2.25) -- (2.65,2.);
\draw[ns2] (2.3125,1.625) -- (2.3125,2.625);
\draw[ns2] (1.975,2.) -- (2.3125,1.875);
\draw[ns2] (2.1438,1.7) -- (2.1438,2.2);
\draw[ns2] (2.65,3.5) -- (4,3.);
\draw[ns2] (3.35,2.25) -- (3.35,4.25);
%}

%\visible<3->{
\node[] at (6.65,-2.5) {\LARGE $t_{\textnormal{DoF}}^{\trans, 3}$};
\draw[HSUred,fill=HSUred](6.65,0.5)circle(3.5pt);
\draw[->,ns1, HSUred, dotted] (6.65,0.5) -- (6.65,-1.8);
\draw[HSUred,fill=HSUred](6.65,-2)circle(3.5pt);

% Third adaptive spatial mesh (on t^3_DoF of Q_n)
%boundary
\draw[ns2] (8,0) -- (8,4);
\draw[ns2] (0,4) -- (0,0);
\draw[ns2] (5.3,1) -- (8,0);
\draw[ns2] (8,0) -- (8,4);
\draw[ns2] (8,4) -- (5.3,5);
\draw[ns2] (5.3,5) -- (5.3,1);
% Inner Mesh
\draw[ns2] (5.3,3) -- (8,2);
\draw[ns2] (6.65,0.5) -- (6.65,4.5);
\draw[ns2] (5.3,2) -- (6.65,1.5);
\draw[ns2] (5.975,0.75) -- (5.975,2.75);
\draw[ns2] (5.975,2.25) -- (6.65,2.);
\draw[ns2] (6.3125,1.625) -- (6.3125,2.625);
\draw[ns2] (5.975,2.) -- (6.3125,1.875);
\draw[ns2] (6.1438,1.7) -- (6.1438,2.2);
\draw[ns2] (6.65,3.5) -- (8,3.);
\draw[ns2] (7.35,2.25) -- (7.35,4.25);
%}

%%%%%%%%%%%%%%%%%%%%%%%%%%%%%%%%%%%%%%%%%%%%%%%%%%%%%%%%%%%%%%%%%%%%%%%%%%%%%%%%
%
% Q_n+1

%\visible<4->{ % Q_n+1

\node[] at (9.65,-2.5) {\LARGE $t_{\textnormal{DoF}}^{\trans, 1}$};
\draw[HSUred,fill=HSUred](9.65,0.5)circle(3.5pt);
\draw[->,ns1, HSUred, dotted] (9.65,0.5) -- (9.65,-1.8);
\draw[HSUred,fill=HSUred](9.65,-2)circle(3.5pt);

% First adaptive spatial mesh (on t^1_DoF of Q_n+1)
%
% Mesh boundary
\draw[ns1] (8.3,1) -- (11,0);
\draw[ns1] (11,0) -- (11,4);
\draw[ns1] (11,4) -- (8.3,5);
\draw[ns1] (8.3,5) -- (8.3,1);
% Inner Mesh
\draw[ns1] (8.3,3) -- (11,2);
\draw[ns1] (9.65,0.5) -- (9.65,4.5);
\draw[ns1] (9.65,1.5) -- (11,1.0);
\draw[ns1] (10.35,0.25) -- (10.35,2.25);
\draw[ns1] (9.65,2.0) -- (10.325,1.75);
\draw[ns1] (10,1.375) -- (10,2.375);
\draw[ns1] (10.325,1.75) -- (11,1.5);
\draw[ns1] (10.675,1.125) -- (10.675,2.125);
\draw[ns1] (10,2.125) -- (10.35,2.0);
\draw[ns1] (10.175,1.8125) -- (10.175,2.3125);
\draw[ns1] (8.3,2) -- (9.65,1.5);
\draw[ns1] (8.975,0.75) -- (8.975,2.75);

\node[] at (19.15,-2.5) {\LARGE $t_{n+1}^{\trans}$};
\draw[HSUred,fill=HSUred](19.15,0.5)circle(3.5pt);
\draw[->,ns1, HSUred, dotted] (19.15,0.5) -- (19.15,-1.8);
\draw[HSUred,fill=HSUred](19.15,-2)circle(3.5pt);

\draw[very thick,HSUred,decorate,decoration={brace,amplitude=14pt}]
(6.8,5) -- (17.8,5) node[midway, above,yshift=14pt,]{\LARGE $Q_{n+1}^{\trans}$};

% Middle extension t_n
\draw[ns1] (6.8,5) -- (8.3,5);
\draw[ns1] (9.5,4) -- (11,4);
\draw[ns1,dashed] (6.8,1) -- (8.3,1);
\draw[ns1] (9.5,0) -- (11,0);

\draw[ns1] (11,0) -- (19,0);
\draw[ns1] (19,4) -- (11,4);
\draw[ns1] (8.3,5) -- (16.3,5); % Frontedge
\draw[ns1, dashed] (8.3,1) -- (16.3,1); % Backedge

% back extension t_n+1
\draw[ns1] (16.3,5) -- (17.8,5);
\draw[ns1] (19,4) -- (20.5,4);
\draw[ns1,dashed] (16.3,1) -- (17.8,1);
\draw[ns1] (19,0) -- (20.5,0);

\draw[ns1,dashed] (17.8,1) -- (20.5,0);
\draw[ns1] (20.5,0) -- (20.5,4);
\draw[ns1] (20.5,4) -- (17.8,5);
\draw[ns1,dashed] (17.8,5) -- (17.8,1);

\node[] at (13.65,-2.5) {\LARGE $t_{\textnormal{DoF}}^{\trans, 2}$};
\draw[HSUred,fill=HSUred](13.65,0.5)circle(3.5pt);
\draw[->,ns1, HSUred, dotted] (13.65,0.5) -- (13.65,-1.8);
\draw[HSUred,fill=HSUred](13.65,-2)circle(3.5pt);

% Second adaptive spatial mesh (on t^2_DoF of Q_n+1)
%
% Mesh boundary
\draw[ns2] (12.3,1) -- (15,0);
\draw[ns2] (15,0) -- (15,4);
\draw[ns2] (15,4) -- (12.3,5);
\draw[ns2] (12.3,5) -- (12.3,1);
% Inner Mesh
\draw[ns2] (12.3,3) -- (15,2);
\draw[ns2] (13.65,0.5) -- (13.65,4.5);
\draw[ns2] (13.65,1.5) -- (15,1.0);
\draw[ns2] (14.35,0.25) -- (14.35,2.25);
\draw[ns2] (13.65,2.0) -- (14.325,1.75);
\draw[ns2] (14,1.375) -- (14,2.375);
\draw[ns2] (14.325,1.75) -- (15,1.5);
\draw[ns2] (14.675,1.125) -- (14.675,2.125);
\draw[ns2] (14,2.125) -- (14.35,2.0);
\draw[ns2] (14.175,1.8125) -- (14.175,2.3125);
\draw[ns2] (12.3,2) -- (13.65,1.5);
\draw[ns2] (12.975,0.75) -- (12.975,2.75);

\node[] at (17.65,-2.5) {\LARGE $t_{\textnormal{DoF}}^{\trans, 3}$};
\draw[HSUred,fill=HSUred](17.65,0.5)circle(3.5pt);
\draw[->,ns1, HSUred, dotted] (17.65,0.5) -- (17.65,-1.8);
\draw[HSUred,fill=HSUred](17.65,-2)circle(3.5pt);

% Third adaptive spatial mesh (on t^3_DoF of Q_n+1)
%
% Mesh boundary
\draw[ns2] (16.3,1) -- (19,0);
\draw[ns2] (19,0) -- (19,4);
\draw[ns2] (19,4) -- (16.3,5);
\draw[ns2] (16.3,5) -- (16.3,1);
% Inner Mesh
\draw[ns2] (16.3,3) -- (19,2);
\draw[ns2] (17.65,0.5) -- (17.65,4.5);
\draw[ns2] (17.65,1.5) -- (19,1.0);
\draw[ns2] (18.35,0.25) -- (18.35,2.25);
\draw[ns2] (17.65,2.0) -- (18.325,1.75);
\draw[ns2] (18,1.375) -- (18,2.375);
\draw[ns2] (18.325,1.75) -- (19,1.5);
\draw[ns2] (18.675,1.125) -- (18.675,2.125);
\draw[ns2] (18,2.125) -- (18.35,2.0);
\draw[ns2] (18.175,1.8125) -- (18.175,2.3125);
\draw[ns2] (16.3,2) -- (17.65,1.5);
\draw[ns2] (16.975,0.75) -- (16.975,2.75);

%} % visible Q_n+1

% t axis
\node[] at (21,-2) {\LARGE $t$};
\draw[->,ns1] (-4.2,-2) -- (20.5,-2);
\end{tikzpicture}
}

\caption{Exemplary solution mesh transfer from a flow slab $Q_n^{\flow}$
(discontinuous Galerkin dG($0$) time discretization generated 
with one Gaussian quadrature point) to transport slabs $Q_n^{\trans}$ and
$Q_{n+1}^{\trans}$ (discontinuous Galerkin dG($2$) time 
discretization generated with three Gaussian quadrature points), respectively.
Here, each of the illustrated slabs consists of one cell in time and an 
independent and adaptively refined spatial triangulation.
}
\label{fig:4:SolutionTransferStokes}
\end{figure}
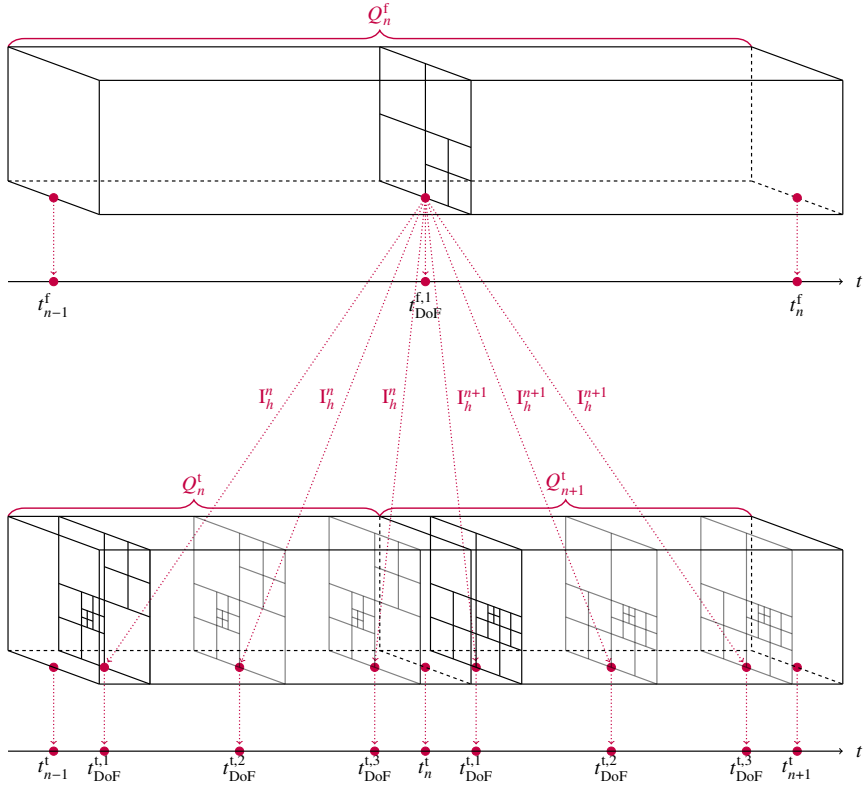
%%%%%%%%%%%%%%%%%%%%%%%%%%%%%%%%%%%%%%%%%%%%%%%%%%%%%%%%%%%%%%%%%%%%%%%%%%%%%%%%
%% End Figure 3 Solution Transfer Stokes

For the sake of simplicity, we approximate the solution $\{\mathbf{v}, p \}$ 
of the flow problem on each $I_n^\flow$ by means of a globally 
piecewise constant discontinuous Galerkin (dG($0$)) time approximation.
Thus, in accordance with the above described conditions for the temporal 
discretizations of both subproblems, we simply have to guarantee the correct 
choice of the corresponding flow slab when solving on the current transport slab 
and avoid an additional evaluation of the flow solution at the temporal degrees 
of freedom within the slab.
%
% For example, with regard to Fig.~\ref{fig:2:InitializationTemporalMeshes},
% for the approximation of the transport solution on $I_4^\trans$ (corresponding 
% to slab $Q_4^\trans$) the flow solution on $I_2^{\flow}$ (corresponding to slab 
% $Q_2^{\flow}$) is relevant.

The flow solution transfer with regard to the spatial meshes is handled by means 
of introducing a temporary additional flow triangulation build as a copy of the 
transport triangulation. 
Then, the flow field solution is interpolated to this temporary triangulation using
an interpolation operator $\operatorname{I}_{h}^{n}$ onto the primal or dual finite
element space used on the current slab $Q_n^\trans$, handled by a precasted function
called \texttt{interpolate\_to\_different\_mesh()} within the \texttt{deal.II} library 
\cite{dealiiReference93}.
An exemplary solution transfer of the flow field solution to the transport 
meshes is illustrated in Fig.~\ref{fig:4:SolutionTransferStokes}.
Of course, this approach entails an additional interpolation error. However, in 
our numerical examples in Sec.~\ref{sec:5:examples}, we observe that this impact 
as well as the restriction to a dG($0$) approximation in time for the flow 
problem is negligibly to obtain quantitatively good results with regard to the 
transport solution.

\subsubsection*{Realization of Adaptive Mesh Refinement in Space and Time}

A further important aspect when dealing with goal-oriented error control of 
coupled problems is the treatment of adaptive mesh refinement within the single
subproblems as well as the question of the consequences for the respective 
other meshes due to these refinements.
The adaptive refinement is performed cell-wise in space and time based on the 
local error indicators $\tilde{\eta}_{\sigma}^{\flow}, \tilde{\eta}_{h}^{\flow}$
and $\eta_{\tau}^{\trans}, \eta_{h}^{\trans}$ for the flow and transport problem, 
respectively.
With regard to a single subproblem, the adaptive mesh refinement process 
including the involvement of additional space-time slabs is handled in the 
following way, cf.~\cite{Bruchhaeuser2022a} for more details:
\begin{itemize}
\item Store the space-time slabs within a \texttt{std::list} object; 
cf.~\cite{Koecher2019} for more details about this list approach.
\item Execute at first the spatial refinement and coarsening of the underlying
triangulations on all slabs.
\item Refine in time by involving new created slabs by copying the just refined 
spatial triangulation of a slab that is marked for refinement.
\end{itemize}
In order to respond the question above regarding the resulting consequences for 
the meshes of the other subproblem within the adaptive refinement process, 
we assume the following:
\begin{itemize}
\item Refine the transport meshes after the flow meshes.
\item As in the case of the initialization of the temporal meshes, the endpoints 
of the temporal mesh of the flow problem must match with endpoints of the 
temporal mesh of the transport problem.
\end{itemize}

\subsection{Algorithm}
\label{sec:4:3:algorithm}

Finally, we present our cost-efficient space-time adaptive algorithm for coupled 
flow and transport problems. 
Here, we are primary interested to control the transport problem under the 
condition that the influence of the error in the flow problem stays small
by reducing the incidental numerical costs at the same time. 
This is enabled by using auxiliary, non-weighted error indicators 
\eqref{eq:24:auxiliary_eta_stokes} in the flow problem avoiding an explicit 
computation of the dual flow problem.
Moreover, the temporal weights arising within the DWR-based error indicators in
the transport problem are approximated by a higher-order reconstruction 
approach using the same polynomial degree $r$ within the discontinuous Galerkin 
dG($r$) time discretization for the primal and dual transport problem.

\noindent\rule{\textwidth}{1pt}
  \begin{center}
   \textbf{Algorithm: DWR cost-efficient multirate space-time adaptivity}
  \end{center}
\vspace{-0.3cm}
\noindent\rule{\textwidth}{0.5pt}
\textbf{Initialization:}
Generate the initial space-time slabs $Q_{n}^{\trans,1}=\mathcal{T}_{h,n}^{\trans,1}\times 
\mathcal{T}_{\tau,n}^{\trans,1}\,,n=1,\dots,N^{\trans,1}\,,$  as well as
$Q_{n}^{\flow,1}=\mathcal{T}_{h,n}^{\flow,1}\times 
\mathcal{T}_{\sigma,n}^{\flow,1}\,,n=1,\dots,N^{\flow,1}\,,$ 
$N^{\flow,1} \leq N^{\trans,1}\,,$ for the transport and Stokes flow problem, 
respectively, where we restrict $\mathcal{T}_{\tau,n}^{\trans,1},
\mathcal{T}_{\sigma,n}^{\flow,1}$ to consist of only one cell in time
on each slab. Set \textbf{DWR-loop} $\ell=1,\dots$:

\noindent\rule{\textwidth}{0.5pt}
\begin{enumerate}
\item %
  \textbf{Find the solution} $\mathbf{u}_{\sigma h}^{\ell}=\{\mathbf{v}_{\sigma h}^{\ell},p_{\sigma h}^{\ell}\} \in 
\mathcal{Y}_{\sigma h}^{0,p_{\mathbf{v}},p_{p}}$, $p_{p}+1=p_{\mathbf{v}} \geq 2$
  of the flow problem \eqref{eq:16:fully_stokes}.

\item \textbf{Find the primal solution}
  $\concentration_{\tau h}^{\ell} \in X_{\tau h}^{r,p}$
  of the stabilized transport problem \eqref{eq:17:fully_transport}.

\item \textbf{Break} if the \textbf{goal} is reached 
(here, i.e.  
$\|\concentration-\concentration_{\tau h}^{\ell}\|< \texttt{tol}$).

\item \textbf{Find the dual solution}
  $\dualz_{\tau h}^{\ell} \in X_{\tau h}^{r,q}, q > p,$
  of the dual transport problem 
%   (third Eq. in 
  \eqref{eq:?:DualProblems}.
%   ).

\item Evaluate the \textbf{localized DWR-weighted error indicators}
  $\eta_{\tau}^{\trans,\ell}$ and $\eta_{h}^{\trans,\ell}$ given by \eqref{eq:24:eta_transport} 
  for the transport problem and the \textbf{localized non-weighted Kelly error indicators}
  $\tilde{\eta}_{\sigma}^{\flow,\ell}$ and $\tilde{\eta}_{h}^{\flow,\ell}$ given
  by \eqref{eq:24:auxiliary_eta_stokes} for the flow problem.
   
\item \textbf{If} 
$|\tilde{\eta}_{\sigma h}^{\flow,\ell}|=|\tilde{\eta}_{\sigma}^{\flow,\ell}| + |\tilde{\eta}_{h}^{\flow,\ell}|
>
\varpi\,|\eta_{\tau}^{\trans,\ell}| + |\eta_{h}^{\trans,\ell}|
\,,\varpi \geq 1$: 

Refine the temporal and spatial meshes of the \textbf{flow} problem as follows:
\begin{enumerate}
\item[(i)] 
\textbf{Mark the slabs} $Q_{\tilde{n}}^{\flow,\ell}$, $\tilde{n}\in\{1,\dots,N^{\flow,\ell}\}$,
\textbf{for temporal refinement} if the corresponding 
$\tilde{\eta}_{\sigma}^{\flow,\tilde{n},\ell}$ is in 
the set of $\theta_\sigma^{\textnormal{top}}\,,0\leq\theta_\sigma^{\textnormal{top}}\leq 1\,,$ 
percent of the worst indicators.
\item[(ii)] 
\textbf{Mark the cells\,} $\tilde{K}^{\flow}\in\mathcal{T}_{h,n}^{\flow,\ell}$
\textbf{\,for spatial refinement} if the corresponding 
$\tilde{\eta}_{h}^{\flow,n,\ell}|_{\tilde{K}^{\flow}}$ 
is in the set of $\theta_{h,1}^\textnormal{\flow,top}$ or 
$\theta_{h,2}^\textnormal{\flow,top}$
(for a slab that is or is not marked for temporal refinement),
$0\leq\theta_{h,2}^\textnormal{\flow,top}
\leq
\theta_{h,1}^\textnormal{\flow,top}\leq 1\,,$ 
percent of the worst indicators,
\textbf{or}, respectively,
mark \textbf{for spatial coarsening} if 
$\tilde{\eta}_{h}^{\flow,n,\ell}|_{\tilde{K}^{\flow}}$  is in the 
set of $\theta_h^\textnormal{\flow,bottom}\,,
0\leq\theta_h^\textnormal{\flow,bottom}\leq 1\,,$ percent of the best indicators.
\item[(iii)] \textbf{Execute spatial adaptations} on all slabs 
of the flow problem under the use of mesh smoothing operators.
\item[(iv)] \textbf{Execute temporal refinement} on all slabs of the flow 
problem. 
\end{enumerate}
\textbf{Else}:
Do not refine the temporal and spatial meshes of the flow problem and continue with Step~7.
Refine the temporal and spatial meshes of the \textbf{transport} problem 
as follows:
\begin{enumerate}
\item[(i)] \textbf{If} 
$|\eta_{\tau}^{\trans,\ell}|>
\omega\,|\eta_{h}^{\trans,\ell}|\,,\omega \geq 1$: 
\textbf{Mark the slabs} $Q_{\tilde{n}}^{\trans,\ell}$, $\tilde{n}\in\{1,\dots,N^{\trans,\ell}\}$,
\textbf{for temporal refinement} if the corresponding 
$\eta_{\tau}^{\trans,\tilde{n},\ell}$ is in 
the set of $\theta_\tau^\mathrm{top}\,,0\leq\theta_\tau^\mathrm{top}\leq 1\,,$ 
percent of the worst indicators.
  
\item[(ii)] \textbf{Else if} 
$|\eta_{h}^{\trans,\ell}|>
\omega\,|\eta_{\tau}^{\trans,\ell}|$: 
\textbf{Mark the cells\,} $\tilde{K}^{\trans}\in\mathcal{T}_{h,n}^{\trans,\ell}$
\textbf{\,for spatial refinement} if the corresponding 
$\eta_{h}^{\trans,n,\ell}|_{\tilde{K}\trans,}$ 
is in the set of $\theta_{h,1}^{\trans,\mathrm{top}}$ or 
$\theta_{h,2}^{\trans,\mathrm{top}}$
(for a slab that is or is not marked for temporal refinement),
$0\leq\theta_{h,2}^{\trans,\mathrm{top}}
\leq
\theta_{h,1}^{\trans,\mathrm{top}}\leq 1\,,$ 
percent of the worst indicators,
\textbf{or}, respectively,
mark \textbf{for spatial coarsening} if 
$\eta_{h}^{\trans,n,\ell}|_{\tilde{K}\trans,}$  is in the 
set of $\theta_h^{\trans,\mathrm{bottom}}\,,
0\leq\theta_h^{\trans,\mathrm{bottom}}\leq 1\,,$ percent of the best 
indicators.

\item[(iii)] \textbf{Else}: 
\textbf{Mark the slabs} $Q_{\tilde{n}}^{\trans,\ell}$
\textbf{for temporal refinement} as well as \textbf{mark the cells\,} 
$\tilde{K}^{\trans}\in\mathcal{T}_{h,n}^{\trans,\ell}$ \textbf{for spatial coarsening} 
and \textbf{refinement} as described in Step~7(i)-(ii).% and Step~7(ii), respectively.
  
\item[(iv)] \textbf{Execute spatial adaptations} on all slabs
of the transport problem under the use of mesh smoothing operators.
  
\item[(v)] \textbf{Execute temporal refinement} on all slabs of the transport 
problem. 
\end{enumerate}

\item Increase $\ell$ to $\ell+1$ and return to Step~1.
\end{enumerate}
\vspace{-0.45cm}
\noindent\rule{\textwidth}{0.5pt}
\newpage
\begin{remark}
\label{rem:3:algorithm}
Let us remark some aspects about the multirate adaptive algorithm:
\begin{itemize}
\item For the spatial discretization of the flow problem we are using 
Taylor-Hood elements $Q_p/Q_{p-1}\,, p \geq 2$.
In order to ensure the conditions to the temporal meshes outlined in 
Sec.~\ref{sec:4:2:MultirateAspects}, we refine the transport meshes after the  
flow meshes such that the endpoints of the temporal mesh of the flow solver 
match with endpoints of the temporal mesh of the transport solver.

\item Within the framework of coupled problems, it is essential to know which 
equation contributes most to the overall error. 
For this purpose, the problem equilibration constant $\varpi$ 
(a value in the range of $ 1 \leq \varpi \leq 3$ is used in our 
numerical experiments) is introduced in Step~6 and 7. 

\item To ensure an equilibrated reduction of the temporal and spatial 
discretization error within the transport problem, the equilibration constant 
$\omega$ (a value in the range of $1.5  \leq \omega \leq 3.5$ is used in our
numerical examples) is introduced in Step~7.  

\item Within the Steps 2, 4 and 5 of the algorithm, the computed flow field
$\mathbf{v}_{\sigma h}$ of the flow problem is interpolated to the adaptively 
refined spatial and temporal triangulation of the current space-time transport 
slab as described in Sec.~\ref{sec:4:2:MultirateAspects}.
By means of the variational equation \eqref{eq:16:fully_stokes} along 
with the definition \eqref{eq:12:BsigmaFsigma} of the bilinear form 
$B_{\sigma}$ the fully discrete flow field $\mathbf{v}_{\sigma h}$ is weakly 
divergence free on the spatial mesh $\mathcal T_{h,n}^{\flow}$ of the flow problem. 
However, on the spatial mesh $\mathcal T_{h,n}^{\trans}$ of the transport problem this 
constraint is in general violated due to the different decomposition of $\Omega$. 
The lack of divergence-freeness might be the source of an approximation error. 
By the application of an additional Helmholtz or Stokes projection 
(cf. \cite[Rem.~2.1]{Brenner2014}) the divergence-free constraint
can be recovered on the spatial mesh of the transport quantity $u_{\tau h}$. 
This amounts to the solution of a problem of Stokes type for each of the fully 
discrete flow fields of the temporal flow mesh (cf.~Fig.~\ref{fig:2:InitializationTemporalMeshes}). 
In the numerical experiments of Sec.~\ref{sec:5:examples} we did not observe any 
problems without post-proccesing the velocity fields to ensure discrete 
divergence-freeness on $\mathcal T_{h,n}^{\trans}$.

\item Our simulation tools of the \texttt{DTM++} project are frontend solvers
for the \texttt{deal.II} library; cf. \cite{dealiiReference93}.
Technical details of the implementation are given in
\cite{Koecher2019,Bause2021}.
\end{itemize}
\end{remark}

\section{Numerical Examples}
\label{sec:5:examples}

In this section, we study robustness, stability as well as computational 
accuracy and efficiency of our cost reduced algorithm for coupled flow and 
transport problems.
The first example is an academic test problem with given analytical solutions. 
It serves to present the performance properties of the algorithm with regard 
to adaptive mesh refinement in space and time.
The second example is motivated by a problem of physical relevance in which we 
simulate a convection-dominated transport of a species through a channel with a 
constraint.
Finally, in a third part we modify this example to the case of a strongly 
convection-dominated transport in order to investigate the interaction
of stabilization combined with goal-oriented error control.

\subsection{Example 1 (Space-Time Adaptivity Studies for the Coupled Problem)}
\label{sec:5:1}

In a first numerical example, we study the algorithm introduced in Sec.~\ref{sec:4:3:algorithm} 
with regard to accuracy, reliability and efficiency reasons.
For this purpose, we consider a so-called effectivity index given as the ratio 
of the estimated error over the exact error, i.e.
\begin{equation}
\label{eq:30:Ieff}
\mathcal{I}_{\textnormal{eff}} = \left|\frac{\eta_{\tau}^{\trans}+\eta_{h}^{\trans}}{J(\concentration)
-J(\concentration_{\tau h})}\right|\,.
\end{equation}
Desirably, this index should be close to one. 
Moreover, it is essential to guarantee an equilibrated reduction
of the temporal and spatial discretization error ensured by means of well-balanced
error indicators $\eta_{\tau}^{\trans}\approx\eta_{h}^{\trans}$.

For the first test case, we investigate the coupled flow and transport problem
with the following setting.
Regarding the flow problem \eqref{eq:1:stokes_problem}, we choose the volume 
force term $\mathbf{f}$ in such a way that the exact solution $\mathbf{u}=\{\mathbf{v}, p\}$
is given by
\begin{equation}
\label{eq:26:BRdynamic}
\begin{array}{r@{\,}c@{\,}l}
\mathbf{v}(\bold x, t) &:=&
\left(
\begin{array}{c}
\sin(t)\sin^2(\pi x_1)\sin(\pi x_2)\cos(\pi x_2) \\
-\sin(t)\sin(\pi x_1)\cos(\pi x_1)\sin^2(\pi x_2)
\end{array}\right)
\,,
\\[2.5ex]
p(\bold x, t) &:=&
\sin(t)\sin(\pi x_1)\cos(\pi x_1)\sin(\pi x_2)\cos(\pi x_2)\,,
\end{array}
\end{equation}
with $\bold x = (x_1, x_2)^\top \in \mathbb{R}^2\,, t \in \mathbb{R}$ and a 
divergence-free flow field $\nabla \cdot \mathbf{v} = 0$.
The viscosity is set to $\nu=0.5$.
The problem is defined on $Q=\Omega\times I:=(0,1)^2\times (0,1]$ and the initial 
and boundary conditions are given as
\begin{displaymath}
\convection = 0 \,\textnormal{ on }\, \Sigma_0 = \Omega \times \{ 0 \}\,,
\quad 
\mathbf{v} = 0 \,\textnormal{ on }\, \Sigma = \partial\Omega \times (0,1)\,.
\end{displaymath}
This is a typical test problem for time-dependent incompressible flow and can 
be found, for instance, in \cite[Example 1]{Besier2012}. 
Regarding the convection-diffusion transport problem \eqref{eq:2:transport_problem} 
coupled with the flow problem via the flow field $\mathbf{v}_{\sigma h}$, we 
choose the force $g$ in such a way that the exact solution $u$ is given by
\begin{equation}
\label{eq:26:KB2}
\begin{array}{l@{\,}c@{\,}l}
u(\bold x, t) &:=&
u_1 \cdot u_2\,,\,\,
\bold x = (x_1, x_2)^\top \in \mathbb{R}^2 \textnormal{ and }
t \in \mathbb{R}\,,\\[.5ex]
u_1(\bold x, t) &:=&
  (1 + a \cdot ( (x_1 - m_1(t))^2 + (x_2 - m_2(t))^2  ) )^{-1}\,,\\[.5ex]
u_2(t) &:=& \nu_1(t) \cdot s \cdot \arctan( \nu_2(t) )\,,
\end{array}
\end{equation}
with $m_1(t) := \frac{1}{2} + \frac{1}{4} \cos(2 \pi t)$ and
$m_2(t) := \frac{1}{2} + \frac{1}{4} \sin(2 \pi t)$, and,
$\nu_1(\hat t) := -1$,
$\nu_2(\hat t) := 5 \pi \cdot (4 \hat t - 1)$,
for $\hat t \in [0, 0.5)$ and
$\nu_1(\hat t) := 1$,
$\nu_2(\hat t) := 5 \pi \cdot (4 (\hat t-0.5) - 1)$,
for $\hat t \in [0.5, 1)$, $\hat t = t - k$,
$k \in \mathbb{N}_0$, and,
scalars $a = 50$ and $s=-\frac{1}{3}$.
The problem is defined on $Q=\Omega\times I:=(0,1)^2\times (0,1]$ and the initial 
and nonhomogeneous Dirichlet boundary conditions are given by the exact 
solution \eqref{eq:26:KB2}.
We choose the diffusion coefficient $\varepsilon = 1$ and the reaction 
coefficient is set to $\alpha = 1$. 
Furthermore, no stabilization ($\delta_K = 0$) is used in this test case.
The solution mimics a counterclockwise rotating cone, cf.~\cite[Ch.~1.4.2]{Hartmann98}, 
and is modified by means of additionally changing the height and orientation of 
the cone over the period $T=1$. 
Precisely, the orientation of the cone switches from negative to positive
while passing $t=0.25$ and from positive to negative while passing $t=0.75$.
This poses a considerable challenge with regard to the adaptive mesh refinements 
in space and time for the transport problem, since the spatial refinements have 
to follow the current position of the cone and the temporal refinements should
detect the special dynamics close to the orientation changes of the cone.
With regard to the characteristic times of the two subproblems defined in 
\eqref{eq:3:characteristic-times}, the respective coefficients are chosen in 
such a way that there holds $t_{\textnormal{transport}} < t_{\textnormal{flow}}$.
Thus, the initial space-time meshes of the transport problem are finer 
compared to the initial meshes of the flow problem, cf. the first row of 
Table~\ref{table:1:L2final-KB2BR12-2232}.
The goal quantity for the transport problem is chosen to control the $L^2$-error
$e_N^-, e_N^-=u(\bold x, T)-u_{\tau h}(\bold x, T)$, at the final time point 
$T=1$, i.e.
\begin{equation}
\label{eq:27:Goal-L2final}
J_T(u)= \frac{(u(\bold x, T),e_N^-)}{\|e_N^-\|_{T}}\,,
\end{equation}
where $\|e_N^-\|_{T}$ denotes the $L^2$-norm at the final time point $T$.
Finally, as outlined in Rem.~\ref{rem:3:algorithm}, the tuning parameters of the 
algorithm are chosen in a way to balance automatically the potential misfit of 
the spatial and temporal errors as
\begin{displaymath}
\begin{array}{l@{\,}l@{\,}l@{\,}}
\theta_{h,1}^{\trans,\textnormal{top}} \geq \theta_{h,2}^{\trans,\textnormal{top}}  =  
\frac{1}{2} \cdot \frac{|\eta_h^{\trans}|}{|\eta_h^{\trans}| + |\eta_\tau^{\trans}|} \,,
& \hspace{0.8cm}
\theta_h^{\trans,\textnormal{bottom}} = 0.02\,, 
& \hspace{0.8cm}
\theta_\tau^\textnormal{top} = \frac{1}{2} \cdot \frac{\eta_\tau^{\trans}}{|\eta_h^{\trans}| + |\eta_\tau^{\trans}|}\,,
\\[1.5ex]
\theta_{h,1}^\textnormal{\flow,top} = \theta_{h,2}^\textnormal{\flow,top} = 0.38\,,
& \hspace{0.8cm} 
\theta_h^\textnormal{\flow,bottom} = 0.02\,,  
& \hspace{0.8cm}
\theta_{\sigma}^{\textnormal{top}} = 1.0\,,
\\[1.5ex]
\varpi = 1.0\,,
& \hspace{0.8cm}
\omega = 2.0\,.
& \hspace{0.8cm}
\end{array}
\end{displaymath}
We approximate the primal and dual transport solutions $u$ and $z$ by means of a 
cG($1$)-dG($0$) and a cG($2$)-dG($0$) method, respectively. The primal flow 
solution $\{\mathbf{v},p\}$ is approximated by using a \{cG(2)-dG(0),cG(1)-dG(0)\}
discretization.
The transport problem is adaptively refined in space and time using the DWR-based
error indicators \eqref{eq:24:eta_transport}.
% along with an approximated flow solution $\mathbf{v}_{\sigma h}$ on coarser global in time and adaptive in space refined meshes.
%
In this example, the adaptivity regarding the flow problem is initially restricted
to the spatial meshes along with a global refinement in time in order to
investigate the refinement behavior based on the non-weighted, auxiliary error 
indicators \eqref{eq:24:auxiliary_eta_stokes} obtained by means of the Kelly 
Error Estimator, cf. Sec.\ref{sec:4:1:EIandWeights}.
The following example given in Sec.~\ref{sec:5:2} then also deals with adaptive 
refinement in time.
For implementational simplicity, we set $\tilde{\eta}_{\sigma h}^{\flow}:=
||\mathbf{v}-\mathbf{v}_{\sigma h}^{2,0}||_{(0,T)\times\Omega}$ in Step~6 of the 
algorithm, i.e. the meshes of the flow problem are refined only if the global 
$L^2(L^2)$-error with respect to the flow field becomes larger than the error 
aimed to be controlled by the goal of the transport problem, here the $L^2$-error 
at the final time point, cf.~\eqref{eq:27:Goal-L2final}.
This is reasonable since as mentioned before we try to control the transport 
problem under the condition that the influence of the error in the flow problem 
stays comparatively small.

In Table~\ref{table:1:L2final-KB2BR12-2232}, we present the development of the 
total discretization error $J_T(e_N^-)=\|e_N^-\|_{T}$ for goal functional 
\eqref{eq:27:Goal-L2final} as well as the global $L^2(L^2)$-error 
$||\mathbf{v}-\mathbf{v}_{\sigma h}^{p,r}||$ for the flow field solution.
Additionally, the spatial and temporal error indicators $\eta_h^{\trans}$ and 
$\eta_{\tau}^{\trans}$ as well as the effectivity index $\mathcal{I}_{\mathrm{eff}}$
during an adaptive refinement process are displayed.
Here and in the following, $\ell$ denotes the refinement level or DWR loop, 
$N$ the number of slabs, $N_{K}^{\text{max}}$ the number of spatial cells on the 
finest mesh within the current loop, and $N_{\text{DoF}}^{\text{tot}}$ the total
space-time degrees of freedom of the flow or transport problem, respectively. 
We observe a very good estimation of the discretization error $J_T(e)$ identified
by effectivity indices close to one (cf. the last column of 
Table~\ref{table:1:L2final-KB2BR12-2232}).
Thus, with regard to accuracy the underlying algorithm performs very well.
Moreover, well-balanced error indicators $\eta_\tau^{\trans}$ and $\eta_h^{\trans}$ 
are obtained in the course of the refinement process (cf. columns ten and eleven 
of Table~\ref{table:1:L2final-KB2BR12-2232}).
Note the existing mismatch of these indicators at the beginning or, for instance,
in the DWR loops 4, 10 or 17, such that the refinement only takes place in time
here. 

%
%%%%%%%%%%%%%%%%%%%%%%%%%%%%%%%%%%%%%%%%%%%%%%%%%%%%%%%%%%%%%%%%%%%%%%%%%%%%%%%%
\begin{table}[t]
\begin{minipage}[c]{\linewidth}
\centering
\resizebox{0.93\linewidth}{!}{%
\begin{tabular}{c | rcrc | rrr | c | cc | c}
\hline
\hline
DWR & \multicolumn{4}{c|}{Flow} & \multicolumn{7}{c}{Transport}
\\
\hline
$\ell$ & $N^{\flow}$ & $N_K^{\flow,\text{max}}$ & $N_{\text{DoF}}^{\flow,\text{tot}}$ & 
$||\mathbf{v}-\mathbf{v}_{\sigma h}^{2,0}||$ &
$N^{\trans}$ & $N_K^{\trans,\text{max}}$ & $N_{\text{DoF}}^{\trans,\text{tot}}$ & $J_T(e^{1,0,2,0})$ & %\|e^{1,0,2,0}\|_{T}
${\eta}_h^{\trans}$ & ${\eta}_\tau^{\trans}$ &
$\mathcal{I}_{\textnormal{eff}}$
\\
\hline
 1 &  2 &    4 &     118 & 1.96e-02 &  10 &   16 &     250 & 2.14e-02 & -1.26e-03 & 4.67e-02 & 2.12\\ 
 2 &    &      &         & 1.96e-02 &  12 &   16 &     300 & 2.10e-02 &  5.06e-03 & 1.87e-03 & 0.33\\%%
 3 &    &      &         & 1.96e-02 &  15 &   40 &     699 & 1.88e-02 &  3.73e-03 & 4.07e-03 & 0.41\\ 
 4 &  4 &   16 &     748 & 4.58e-03 &  19 &  112 &    1975 & 1.05e-02 &  7.08e-04 & 2.97e-03 & 0.35\\ 
 5 &    &      &         & 4.58e-03 &  24 &  112 &    2512 & 6.32e-03 &  1.22e-03 & 2.02e-03 & 0.51\\ 
 6 &    &      &         & 4.58e-03 &  30 &  196 &    4374 & 4.29e-03 &  8.97e-04 & 3.72e-03 & 1.08\\
 7 &  8 &   64 &    5272 & 1.72e-03 &  38 &  196 &    5688 & 3.10e-03 &  1.23e-03 & 2.72e-03 & 1.27\\%%
 8 &    &      &         & 1.72e-03 &  48 &  196 &    7022 & 3.02e-03 &  1.25e-03 & 2.60e-03 & 1.27\\
 9 &    &      &         & 1.72e-03 &  60 &  196 &    8856 & 2.30e-03 &  1.64e-03 & 1.28e-03 & 1.27\\%%
10 &    &      &         & 1.72e-03 &  76 &  268 &   14436 & 2.22e-03 &  8.26e-04 & 1.80e-03 & 1.18\\
11 &    &      &         & 1.72e-03 &  96 &  268 &   18040 & 2.05e-03 &  8.99e-04 & 1.25e-03 & 1.05\\
12 &    &      &         & 1.72e-03 & 121 &  400 &   31421 & 1.47e-03 &  6.71e-04 & 9.64e-04 & 1.11\\
13 & 16 &  208 &   33072 & 1.01e-03 & 153 &  556 &   52213 & 1.22e-03 &  4.30e-04 & 9.60e-04 & 1.14\\ 
14 &    &      &         & 1.01e-03 & 194 &  556 &   64340 & 1.18e-03 &  4.42e-04 & 6.89e-04 & 0.96\\  
15 &    &      &         & 1.01e-03 & 246 & 1060 &  121064 & 8.30e-04 &  3.72e-04 & 4.35e-04 & 0.97\\
16 & 32 &  688 &  215980 & 5.65e-04 & 312 & 1372 &  197706 & 6.79e-04 &  2.43e-04 & 4.14e-04 & 0.97\\ 
17 &    &      &         & 5.65e-04 & 396 & 1744 &  320712 & 4.50e-04 &  1.81e-04 & 4.05e-04 & 1.30\\ 
18 & 64 & 2176 & 1354656 & 3.04e-04 & 502 & 1744 &  415834 & 4.03e-04 &  2.02e-04 & 2.14e-04 & 1.03\\
19 &    &      &         & 3.04e-04 & 637 & 2860 &  707155 & 3.56e-04 &  1.37e-04 & 2.16e-04 & 1.00\\ 
20 & 64 & 2176 & 1354656 & 3.04e-04 & 808 & 3484 & 1146746 & 2.74e-04 &  9.77e-05 & 1.76e-04 & 1.00\\ 
\hline
\end{tabular}
}

% \vskip-1ex
\caption{Adaptive refinement in the transport problem (based on the DWR method)
including effectivity indices for goal functional \eqref{eq:27:Goal-L2final},
with $\varepsilon = 1$, $\delta_0=0$, and $\omega = 2.0$ using a flow solution 
$\mathbf{v}_{\sigma h}^{2,0}$ corresponding to a cG(2)-dG(0) approximation on a 
global refined mesh in time and adaptive refined mesh in space (based on the 
Kelly Error Estimator).
$e^{1,0,2,0}$ corresponds to the adaptive solution approximation 
$u_{\tau h}^{1,0}$ in cG(1)-dG(0) and dual solution approximation
$z_{\tau h}^{2,0}$ in cG(2)-dG(0).}

\label{table:1:L2final-KB2BR12-2232}
\end{minipage}
\end{table}
%%%%%%%%%%%%%%%%%%%%%%%%%%%%%%%%%%%%%%%%%%%%%%%%%%%%%%%%%%%%%%%%%%%%%%%%%%%%%%%%
%% End Table 6 Adaptive Refinement KB2BR21-2232

\begin{figure}[h!]
%%%%%%%%%%%%%%%%%%%%%%%%%%%%%%%%%%%%%%%%%%%%%%%%%%%%%%%%%%%%%%%%%%%%%%%%%%%%%%%%
%% 1
\begin{minipage}{.8\linewidth}
\centering
\begin{tikzpicture}
\begin{axis}[%
width=4.in,
height=1.2in,
scale only axis,
xlabel={t},
ylabel={\textcolor{navyblue}{$\tau_K(I_n^{7})$}, \textcolor{HSUred}{$\sigma_K(I_n^{\textnormal{F},7})$}},
xmin=0.0,
xmax=1.0,
ymin=0.,
ymax=0.15,
yticklabels={0.0,0, 0.05 , 0.1},
]
\addplot [
color=navyblue,
solid,
line width=1.5pt,
mark=*,
mark size = 1.5,
only marks,
mark options={solid,navyblue}
]
table[row sep=crcr]{
0.05 0.05 \\
0.1 0.05 \\
0.125 0.025 \\
0.15 0.025 \\
0.2 0.05 \\
0.225 0.025 \\
0.2375 0.0125 \\
0.25 0.0125 \\
0.2625 0.0125 \\
0.275 0.0125 \\
0.3 0.025 \\
0.35 0.05 \\
0.4 0.05 \\
0.45 0.05 \\
0.5 0.05 \\
0.55 0.05 \\
0.6 0.05 \\
0.625 0.025 \\
0.65 0.025 \\
0.7 0.05 \\
0.725 0.025 \\
0.7375 0.0125 \\
0.75 0.0125 \\
0.7625 0.0125 \\
0.775 0.0125 \\
0.8 0.025 \\
0.85 0.05 \\
0.9 0.05 \\
0.925 0.025 \\
0.95 0.025 \\
0.9625 0.0125 \\
0.975 0.0125 \\
0.98125 0.00625 \\
0.9875 0.00625 \\
0.990625 0.003125 \\
0.99375 0.003125 \\
0.996875 0.003125 \\
1 0.003125 \\
};
\addplot [
color=HSUred,
solid,
line width=1.5pt,
mark=*,
mark size = 1.5,
only marks,
mark options={fill=HSUred}
]
table[row sep=crcr]{
0.125 0.125 \\
0.25 0.125 \\
0.375 0.125 \\
0.5 0.125 \\
0.625 0.125 \\
0.75 0.125 \\
0.875 0.125 \\
1 0.125 \\
};
\end{axis}
\end{tikzpicture}
\end{minipage}

\begin{minipage}{.8\linewidth}
\centering
\begin{tikzpicture}
\begin{axis}[%
width=4.in,
height=1.2in,
scale only axis,
/pgf/number format/.cd, 1000 sep={},
xlabel={t},
ylabel={\textcolor{navyblue}{$\tau_K(I_n^{18})$}, \textcolor{HSUred}{$\sigma_K(I_n^{\textnormal{F},18})$}},
xmin=0.0,
xmax=1.0,
ymax=0.018,
yticklabels={0.0,0, \phantom{-}0.5 , 1, 1.5 },
]

\addplot [
color=navyblue,
solid,
line width=1.0pt,
mark=*,
mark size = 1.5,
only marks,
mark options={solid,navyblue}
]
table[row sep=crcr]{
0.0015625 0.0015625 \\
0.003125 0.0015625 \\
0.0046875 0.0015625 \\
0.00625 0.0015625 \\
0.009375 0.003125 \\
0.0125 0.003125 \\
0.0140625 0.0015625 \\
0.015625 0.0015625 \\
0.0171875 0.0015625 \\
0.01875 0.0015625 \\
0.021875 0.003125 \\
0.025 0.003125 \\
0.028125 0.003125 \\
0.03125 0.003125 \\
0.034375 0.003125 \\
0.0375 0.003125 \\
0.040625 0.003125 \\
0.04375 0.003125 \\
0.046875 0.003125 \\
0.05 0.003125 \\
0.053125 0.003125 \\
0.05625 0.003125 \\
0.059375 0.003125 \\
0.0625 0.003125 \\
0.065625 0.003125 \\
0.06875 0.003125 \\
0.071875 0.003125 \\
0.075 0.003125 \\
0.078125 0.003125 \\
0.08125 0.003125 \\
0.084375 0.003125 \\
0.0875 0.003125 \\
0.090625 0.003125 \\
0.09375 0.003125 \\
0.096875 0.003125 \\
0.1 0.003125 \\
0.103125 0.003125 \\
0.10625 0.003125 \\
0.109375 0.003125 \\
0.1125 0.003125 \\
0.115625 0.003125 \\
0.11875 0.003125 \\
0.121875 0.003125 \\
0.125 0.003125 \\
0.128125 0.003125 \\
0.13125 0.003125 \\
0.134375 0.003125 \\
0.1375 0.003125 \\
0.140625 0.003125 \\
0.14375 0.003125 \\
0.146875 0.003125 \\
0.15 0.003125 \\
0.153125 0.003125 \\
0.15625 0.003125 \\
0.159375 0.003125 \\
0.1625 0.003125 \\
0.164063 0.0015625 \\
0.165625 0.0015625 \\
0.16875 0.003125 \\
0.171875 0.003125 \\
0.175 0.003125 \\
0.178125 0.003125 \\
0.18125 0.003125 \\
0.184375 0.003125 \\
0.1875 0.003125 \\
0.190625 0.003125 \\
0.19375 0.003125 \\
0.196875 0.003125 \\
0.2 0.003125 \\
0.203125 0.003125 \\
0.20625 0.003125 \\
0.209375 0.003125 \\
0.2125 0.003125 \\
0.214063 0.0015625 \\
0.215625 0.0015625 \\
0.217188 0.0015625 \\
0.21875 0.0015625 \\
0.220312 0.0015625 \\
0.221875 0.0015625 \\
0.223438 0.0015625 \\
0.225 0.0015625 \\
0.226562 0.0015625 \\
0.228125 0.0015625 \\
0.229688 0.0015625 \\
0.23125 0.0015625 \\
0.232813 0.0015625 \\
0.234375 0.0015625 \\
0.235156 0.00078125 \\
0.235937 0.00078125 \\
0.236719 0.00078125 \\
0.2375 0.00078125 \\
0.238281 0.00078125 \\
0.239062 0.00078125 \\
0.239844 0.00078125 \\
0.240625 0.00078125 \\
0.241406 0.00078125 \\
0.242188 0.00078125 \\
0.242969 0.00078125 \\
0.24375 0.00078125 \\
0.244531 0.00078125 \\
0.245312 0.00078125 \\
0.246094 0.00078125 \\
0.246875 0.00078125 \\
0.247656 0.00078125 \\
0.248438 0.00078125 \\
0.249219 0.00078125 \\
0.25 0.00078125 \\
0.250781 0.00078125 \\
0.251563 0.00078125 \\
0.252344 0.00078125 \\
0.253125 0.00078125 \\
0.253906 0.00078125 \\
0.254687 0.00078125 \\
0.255469 0.00078125 \\
0.25625 0.00078125 \\
0.257031 0.00078125 \\
0.257812 0.00078125 \\
0.258594 0.00078125 \\
0.259375 0.00078125 \\
0.260156 0.00078125 \\
0.260938 0.00078125 \\
0.2625 0.0015625 \\
0.264062 0.0015625 \\
0.265625 0.0015625 \\
0.267188 0.0015625 \\
0.26875 0.0015625 \\
0.270313 0.0015625 \\
0.271875 0.0015625 \\
0.273438 0.0015625 \\
0.275 0.0015625 \\
0.276563 0.0015625 \\
0.278125 0.0015625 \\
0.279687 0.0015625 \\
0.28125 0.0015625 \\
0.282813 0.0015625 \\
0.284375 0.0015625 \\
0.2875 0.003125 \\
0.290625 0.003125 \\
0.29375 0.003125 \\
0.296875 0.003125 \\
0.3 0.003125 \\
0.303125 0.003125 \\
0.30625 0.003125 \\
0.309375 0.003125 \\
0.3125 0.003125 \\
0.315625 0.003125 \\
0.31875 0.003125 \\
0.321875 0.003125 \\
0.325 0.003125 \\
0.328125 0.003125 \\
0.33125 0.003125 \\
0.334375 0.003125 \\
0.3375 0.003125 \\
0.340625 0.003125 \\
0.34375 0.003125 \\
0.346875 0.003125 \\
0.35 0.003125 \\
0.353125 0.003125 \\
0.35625 0.003125 \\
0.359375 0.003125 \\
0.3625 0.003125 \\
0.365625 0.003125 \\
0.36875 0.003125 \\
0.371875 0.003125 \\
0.375 0.003125 \\
0.378125 0.003125 \\
0.38125 0.003125 \\
0.384375 0.003125 \\
0.3875 0.003125 \\
0.390625 0.003125 \\
0.39375 0.003125 \\
0.396875 0.003125 \\
0.4 0.003125 \\
0.403125 0.003125 \\
0.40625 0.003125 \\
0.409375 0.003125 \\
0.4125 0.003125 \\
0.415625 0.003125 \\
0.41875 0.003125 \\
0.421875 0.003125 \\
0.425 0.003125 \\
0.428125 0.003125 \\
0.43125 0.003125 \\
0.434375 0.003125 \\
0.4375 0.003125 \\
0.440625 0.003125 \\
0.44375 0.003125 \\
0.446875 0.003125 \\
0.45 0.003125 \\
0.453125 0.003125 \\
0.45625 0.003125 \\
0.459375 0.003125 \\
0.4625 0.003125 \\
0.465625 0.003125 \\
0.46875 0.003125 \\
0.471875 0.003125 \\
0.475 0.003125 \\
0.478125 0.003125 \\
0.48125 0.003125 \\
0.484375 0.003125 \\
0.4875 0.003125 \\
0.490625 0.003125 \\
0.49375 0.003125 \\
0.496875 0.003125 \\
0.5 0.003125 \\
0.503125 0.003125 \\
0.50625 0.003125 \\
0.509375 0.003125 \\
0.5125 0.003125 \\
0.515625 0.003125 \\
0.51875 0.003125 \\
0.521875 0.003125 \\
0.525 0.003125 \\
0.528125 0.003125 \\
0.53125 0.003125 \\
0.534375 0.003125 \\
0.5375 0.003125 \\
0.540625 0.003125 \\
0.54375 0.003125 \\
0.546875 0.003125 \\
0.55 0.003125 \\
0.553125 0.003125 \\
0.55625 0.003125 \\
0.559375 0.003125 \\
0.5625 0.003125 \\
0.565625 0.003125 \\
0.56875 0.003125 \\
0.571875 0.003125 \\
0.575 0.003125 \\
0.578125 0.003125 \\
0.58125 0.003125 \\
0.584375 0.003125 \\
0.5875 0.003125 \\
0.590625 0.003125 \\
0.59375 0.003125 \\
0.596875 0.003125 \\
0.6 0.003125 \\
0.603125 0.003125 \\
0.60625 0.003125 \\
0.609375 0.003125 \\
0.6125 0.003125 \\
0.615625 0.003125 \\
0.61875 0.003125 \\
0.621875 0.003125 \\
0.625 0.003125 \\
0.628125 0.003125 \\
0.63125 0.003125 \\
0.634375 0.003125 \\
0.6375 0.003125 \\
0.640625 0.003125 \\
0.64375 0.003125 \\
0.646875 0.003125 \\
0.65 0.003125 \\
0.653125 0.003125 \\
0.65625 0.003125 \\
0.659375 0.003125 \\
0.6625 0.003125 \\
0.664063 0.0015625 \\
0.665625 0.0015625 \\
0.66875 0.003125 \\
0.671875 0.003125 \\
0.675 0.003125 \\
0.678125 0.003125 \\
0.68125 0.003125 \\
0.684375 0.003125 \\
0.6875 0.003125 \\
0.690625 0.003125 \\
0.69375 0.003125 \\
0.696875 0.003125 \\
0.7 0.003125 \\
0.703125 0.003125 \\
0.70625 0.003125 \\
0.709375 0.003125 \\
0.7125 0.003125 \\
0.714063 0.0015625 \\
0.715625 0.0015625 \\
0.717188 0.0015625 \\
0.71875 0.0015625 \\
0.720313 0.0015625 \\
0.721875 0.0015625 \\
0.723438 0.0015625 \\
0.725 0.0015625 \\
0.726563 0.0015625 \\
0.728125 0.0015625 \\
0.729688 0.0015625 \\
0.73125 0.0015625 \\
0.732813 0.0015625 \\
0.733594 0.00078125 \\
0.734375 0.00078125 \\
0.735156 0.00078125 \\
0.735938 0.00078125 \\
0.736719 0.00078125 \\
0.7375 0.00078125 \\
0.738281 0.00078125 \\
0.739063 0.00078125 \\
0.739844 0.00078125 \\
0.740625 0.00078125 \\
0.741406 0.00078125 \\
0.742188 0.00078125 \\
0.742969 0.00078125 \\
0.74375 0.00078125 \\
0.744531 0.00078125 \\
0.745313 0.00078125 \\
0.746094 0.00078125 \\
0.746875 0.00078125 \\
0.747656 0.00078125 \\
0.748437 0.00078125 \\
0.749219 0.00078125 \\
0.75 0.00078125 \\
0.750781 0.00078125 \\
0.751563 0.00078125 \\
0.752344 0.00078125 \\
0.753125 0.00078125 \\
0.753906 0.00078125 \\
0.754687 0.00078125 \\
0.755469 0.00078125 \\
0.75625 0.00078125 \\
0.757031 0.00078125 \\
0.757812 0.00078125 \\
0.758594 0.00078125 \\
0.759375 0.00078125 \\
0.760156 0.00078125 \\
0.760937 0.00078125 \\
0.7625 0.0015625 \\
0.764062 0.0015625 \\
0.765625 0.0015625 \\
0.767188 0.0015625 \\
0.76875 0.0015625 \\
0.770313 0.0015625 \\
0.771875 0.0015625 \\
0.773438 0.0015625 \\
0.775 0.0015625 \\
0.776563 0.0015625 \\
0.778125 0.0015625 \\
0.779687 0.0015625 \\
0.78125 0.0015625 \\
0.782813 0.0015625 \\
0.784375 0.0015625 \\
0.785938 0.0015625 \\
0.7875 0.0015625 \\
0.790625 0.003125 \\
0.79375 0.003125 \\
0.795313 0.0015625 \\
0.796875 0.0015625 \\
0.798438 0.0015625 \\
0.8 0.0015625 \\
0.801563 0.0015625 \\
0.803125 0.0015625 \\
0.804688 0.0015625 \\
0.80625 0.0015625 \\
0.807813 0.0015625 \\
0.809375 0.0015625 \\
0.8125 0.003125 \\
0.815625 0.003125 \\
0.81875 0.003125 \\
0.821875 0.003125 \\
0.825 0.003125 \\
0.828125 0.003125 \\
0.83125 0.003125 \\
0.834375 0.003125 \\
0.8375 0.003125 \\
0.840625 0.003125 \\
0.84375 0.003125 \\
0.846875 0.003125 \\
0.85 0.003125 \\
0.853125 0.003125 \\
0.85625 0.003125 \\
0.859375 0.003125 \\
0.8625 0.003125 \\
0.865625 0.003125 \\
0.86875 0.003125 \\
0.871875 0.003125 \\
0.875 0.003125 \\
0.878125 0.003125 \\
0.88125 0.003125 \\
0.884375 0.003125 \\
0.8875 0.003125 \\
0.890625 0.003125 \\
0.89375 0.003125 \\
0.896875 0.003125 \\
0.9 0.003125 \\
0.903125 0.003125 \\
0.90625 0.003125 \\
0.909375 0.003125 \\
0.9125 0.003125 \\
0.914063 0.0015625 \\
0.915625 0.0015625 \\
0.917188 0.0015625 \\
0.91875 0.0015625 \\
0.920313 0.0015625 \\
0.921875 0.0015625 \\
0.923438 0.0015625 \\
0.925 0.0015625 \\
0.926563 0.0015625 \\
0.928125 0.0015625 \\
0.929688 0.0015625 \\
0.93125 0.0015625 \\
0.932813 0.0015625 \\
0.934375 0.0015625 \\
0.935937 0.0015625 \\
0.9375 0.0015625 \\
0.939063 0.0015625 \\
0.940625 0.0015625 \\
0.942187 0.0015625 \\
0.94375 0.0015625 \\
0.945312 0.0015625 \\
0.946875 0.0015625 \\
0.948437 0.0015625 \\
0.95 0.0015625 \\
0.951562 0.0015625 \\
0.953125 0.0015625 \\
0.954687 0.0015625 \\
0.95625 0.0015625 \\
0.957812 0.0015625 \\
0.959375 0.0015625 \\
0.960937 0.0015625 \\
0.9625 0.0015625 \\
0.964062 0.0015625 \\
0.964844 0.00078125 \\
0.965625 0.00078125 \\
0.966406 0.00078125 \\
0.967187 0.00078125 \\
0.967969 0.00078125 \\
0.96875 0.00078125 \\
0.969531 0.00078125 \\
0.970313 0.00078125 \\
0.971094 0.00078125 \\
0.971875 0.00078125 \\
0.972656 0.00078125 \\
0.973437 0.00078125 \\
0.974219 0.00078125 \\
0.975 0.00078125 \\
0.975781 0.00078125 \\
0.976562 0.00078125 \\
0.977344 0.00078125 \\
0.978125 0.00078125 \\
0.978516 0.000390625 \\
0.978906 0.000390625 \\
0.979687 0.00078125 \\
0.980469 0.00078125 \\
0.98125 0.00078125 \\
0.982031 0.00078125 \\
0.982422 0.000390625 \\
0.982812 0.000390625 \\
0.983203 0.000390625 \\
0.983594 0.000390625 \\
0.983984 0.000390625 \\
0.984375 0.000390625 \\
0.984766 0.000390625 \\
0.985156 0.000390625 \\
0.985547 0.000390625 \\
0.985938 0.000390625 \\
0.986328 0.000390625 \\
0.986719 0.000390625 \\
0.987109 0.000390625 \\
0.9875 0.000390625 \\
0.987891 0.000390625 \\
0.988281 0.000390625 \\
0.988672 0.000390625 \\
0.989063 0.000390625 \\
0.989453 0.000390625 \\
0.989844 0.000390625 \\
0.990234 0.000390625 \\
0.990625 0.000390625 \\
0.991016 0.000390625 \\
0.991406 0.000390625 \\
0.991797 0.000390625 \\
0.992188 0.000390625 \\
0.992578 0.000390625 \\
0.992969 0.000390625 \\
0.993359 0.000390625 \\
0.99375 0.000390625 \\
0.994141 0.000390625 \\
0.994531 0.000390625 \\
0.994922 0.000390625 \\
0.995313 0.000390625 \\
0.995508 0.000195313 \\
0.995703 0.000195312 \\
0.995898 0.000195312 \\
0.996094 0.000195313 \\
0.996289 0.000195313 \\
0.996484 0.000195312 \\
0.99668 0.000195312 \\
0.996875 0.000195313 \\
0.99707 0.000195313 \\
0.997266 0.000195312 \\
0.997461 0.000195312 \\
0.997656 0.000195313 \\
0.997852 0.000195313 \\
0.998047 0.000195313 \\
0.998242 0.000195313 \\
0.998437 0.000195312 \\
0.998633 0.000195312 \\
0.998828 0.000195313 \\
0.999023 0.000195313 \\
0.999219 0.000195313 \\
0.999414 0.000195313 \\
0.999609 0.000195312 \\
0.999805 0.000195312 \\
0.999902 9.76563e-05 \\
0.999951 4.88281e-05 \\
1.0 4.88281e-05 \\
};

\addplot [
color=HSUred,
solid,
line width=1.5pt,
mark=*,
mark size = 1.,
only marks,
mark options={fill=HSUred}
]
table[row sep=crcr]{
0.015625 0.015625 \\
0.03125 0.015625 \\
0.046875 0.015625 \\
0.0625 0.015625 \\
0.078125 0.015625 \\
0.09375 0.015625 \\
0.109375 0.015625 \\
0.125 0.015625 \\
0.140625 0.015625 \\
0.15625 0.015625 \\
0.171875 0.015625 \\
0.1875 0.015625 \\
0.203125 0.015625 \\
0.21875 0.015625 \\
0.234375 0.015625 \\
0.25 0.015625 \\
0.265625 0.015625 \\
0.28125 0.015625 \\
0.296875 0.015625 \\
0.3125 0.015625 \\
0.328125 0.015625 \\
0.34375 0.015625 \\
0.359375 0.015625 \\
0.375 0.015625 \\
0.390625 0.015625 \\
0.40625 0.015625 \\
0.421875 0.015625 \\
0.4375 0.015625 \\
0.453125 0.015625 \\
0.46875 0.015625 \\
0.484375 0.015625 \\
0.5 0.015625 \\
0.515625 0.015625 \\
0.53125 0.015625 \\
0.546875 0.015625 \\
0.5625 0.015625 \\
0.578125 0.015625 \\
0.59375 0.015625 \\
0.609375 0.015625 \\
0.625 0.015625 \\
0.640625 0.015625 \\
0.65625 0.015625 \\
0.671875 0.015625 \\
0.6875 0.015625 \\
0.703125 0.015625 \\
0.71875 0.015625 \\
0.734375 0.015625 \\
0.75 0.015625 \\
0.765625 0.015625 \\
0.78125 0.015625 \\
0.796875 0.015625 \\
0.8125 0.015625 \\
0.828125 0.015625 \\
0.84375 0.015625 \\
0.859375 0.015625 \\
0.875 0.015625 \\
0.890625 0.015625 \\
0.90625 0.015625 \\
0.921875 0.015625 \\
0.9375 0.015625 \\
0.953125 0.015625 \\
0.96875 0.015625 \\
0.984375 0.015625 \\
1.0 0.015625 \\
};
\end{axis}
\end{tikzpicture}
\end{minipage}

\caption{Distribution of the temporal step size $\tau_K$ of the transport
(adaptive, based on the DWR method) and $\sigma_K$ of the flow problem (global)
over the time interval $I=(0,T]$, exemplary after 7 
and 18 DWR-loops, corresponding to Table~\ref{table:1:L2final-KB2BR12-2232}.}
\label{fig:5:L2final-DistributionTauSigmaKB2BR12-1020}
\end{figure}
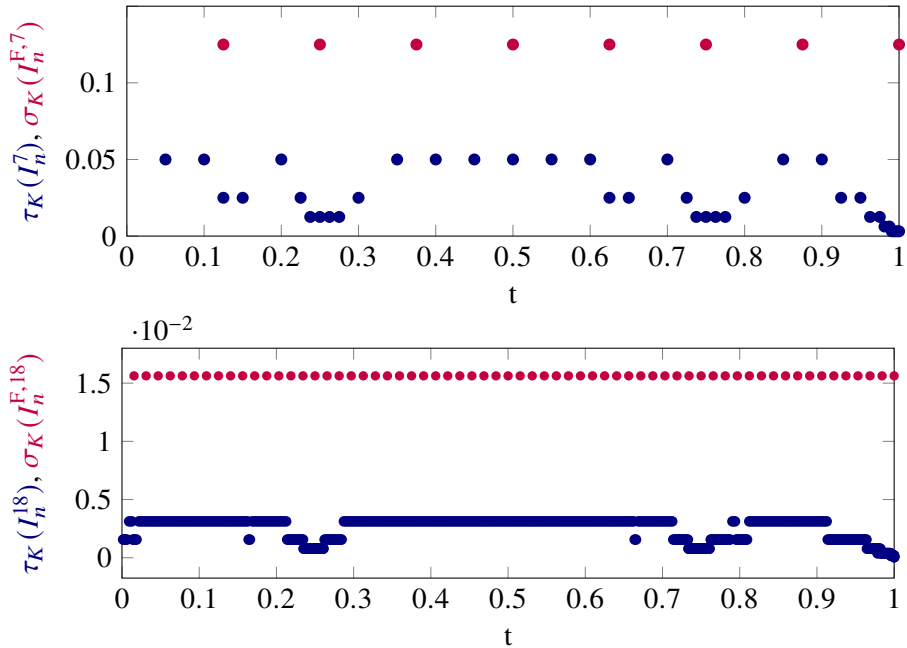

In Fig.~\ref{fig:5:L2final-DistributionTauSigmaKB2BR12-1020}, we 
visualize the distribution of the adaptively determined time cell lengths 
$\tau_K$ of $\mathcal{T}_{\tau,n}$, used for the transport problem, as well as 
the distribution of the globally determined time cell lengths 
$\sigma_K$ of $\mathcal{T}_{\sigma,n}$, used for the flow problem, over 
the whole time interval $I$ for selected DWR loops corresponding to 
Table~\ref{table:1:L2final-KB2BR12-2232}.
We point out that the time steps for the transport problem become significantly
smaller when the cone is changing its orientation ($t=0.25$ and $t=0.75$) as 
well as reaching the final time point $T=1$, while the time steps for the flow 
problem stay comparatively large in the course of the refinement process.
Away from these time points, the temporal mesh of the transport problem is almost 
equally decomposed.
This behavior is desirably since the underlying goal functional \eqref{eq:27:Goal-L2final} 
aims to control the $L^2$-error at the final time point. 
Moreover, the algorithm is able to identify specific dynamics in time arising close
to the time points where the cone is changing its orientation ($t=0.25$ and $t=0.75$)
automatically which indicates the potential of our multirate approach regarding
different characteristic time scales of the subproblems.

%
%%%%%%%%%%%%%%%%%%%%%%%%%%%%%%%%%%%%%%%%%%%%%%%%%%%%%%%%%%%%%%%%%%%%%%%%%%%%%%%%
%% Begin Figure 9 Timesteps DWR-loop Transport and Stokes Flow L2final
%%
%%%%%%%%%%%%%%%%%%%%%%%%%%%%%%%%%%%%%%%%%%%%%%%%%%%%%%%%%%%%%%%%%%%%%%%%%%%%%%%%
%% End Figure 9 Timesteps DWR loop Transport and Stokes Flow L2final

\begin{figure}%[H]
\centering
\subfloat[$t=0.00$]{
  \centering
\begin{minipage}{.23\linewidth}
\centering
\includegraphics[width=3.5cm]{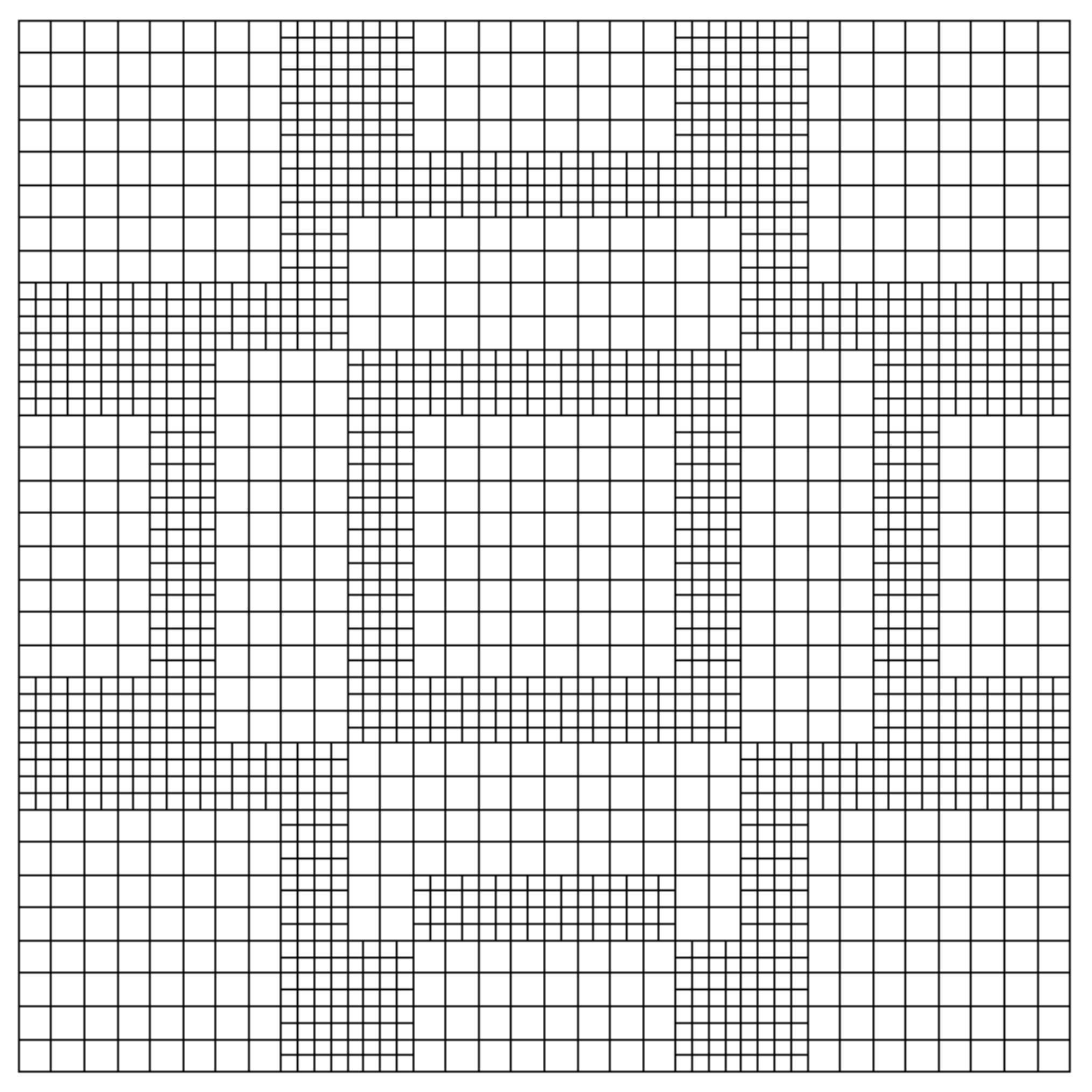}
\end{minipage}
\label{fig:6:Stokes00}
}
\subfloat[$t=0.25$]{
  \centering
\begin{minipage}{.23\linewidth}
\centering
\includegraphics[width=3.5cm]{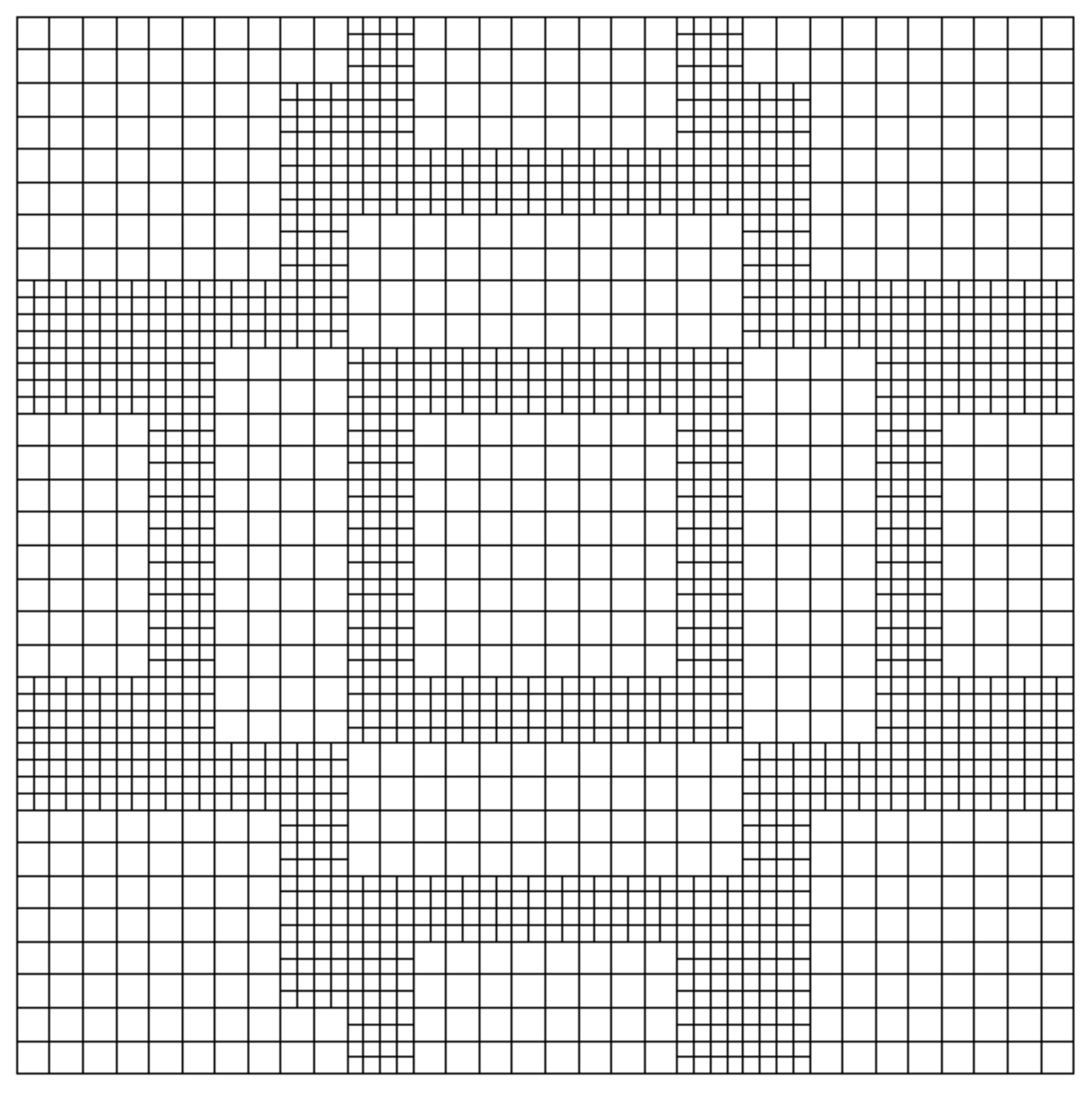}
\end{minipage}
\label{fig:6:Stokes25}
}
\subfloat[$t=0.75$]{
  \centering
\begin{minipage}{.23\linewidth}
\centering
\includegraphics[width=3.5cm]{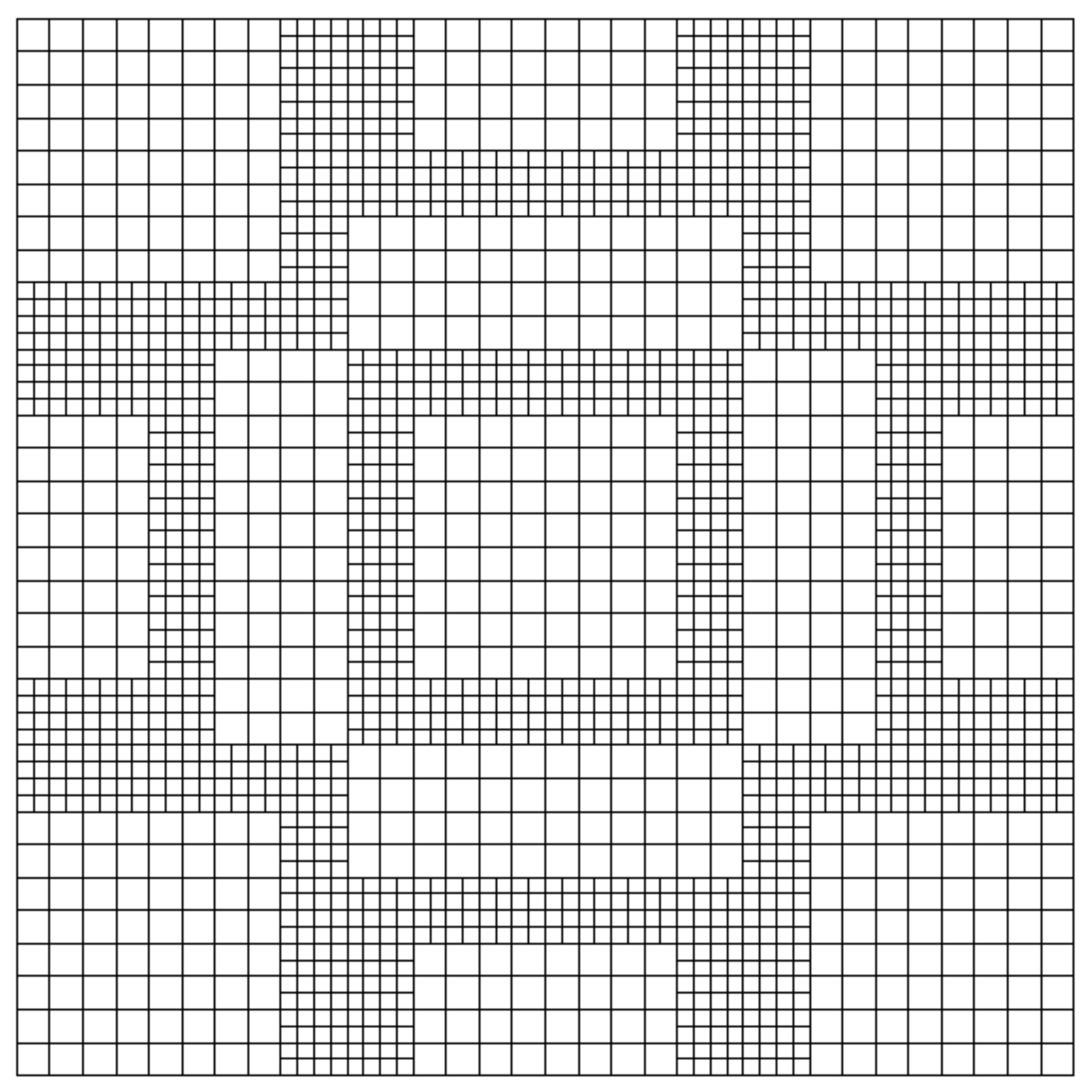}
\end{minipage}
\label{fig:6:Stokes75}
}
\subfloat[$t=1.00$]{
  \centering
\begin{minipage}{.23\linewidth}
\centering
\includegraphics[width=3.5cm]{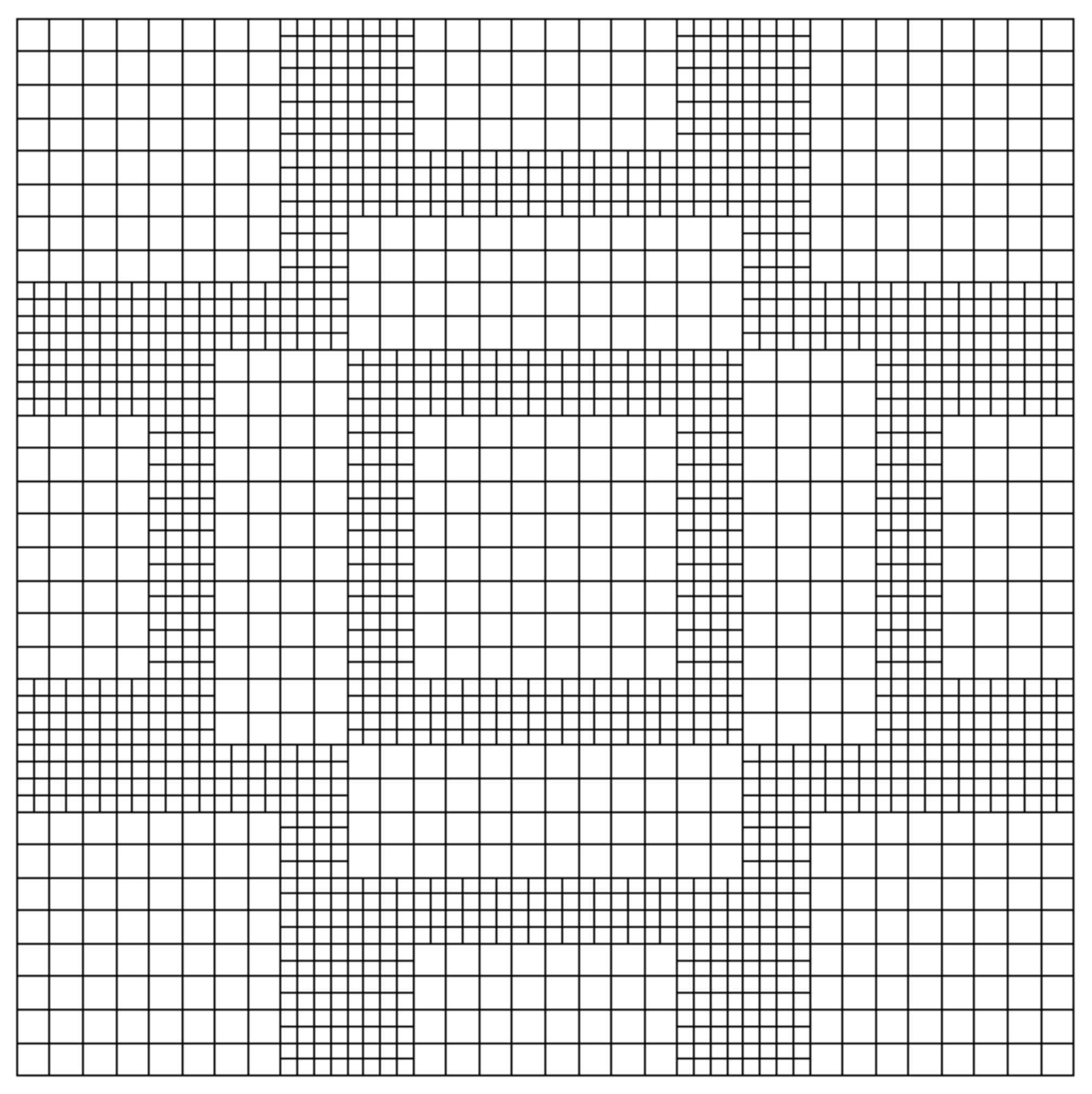}
\end{minipage}
\label{fig:6:Stokes100}
}

\subfloat[$t=0.00$]{
  \centering
\begin{minipage}{.23\linewidth}
\centering
\includegraphics[width=3.5cm]{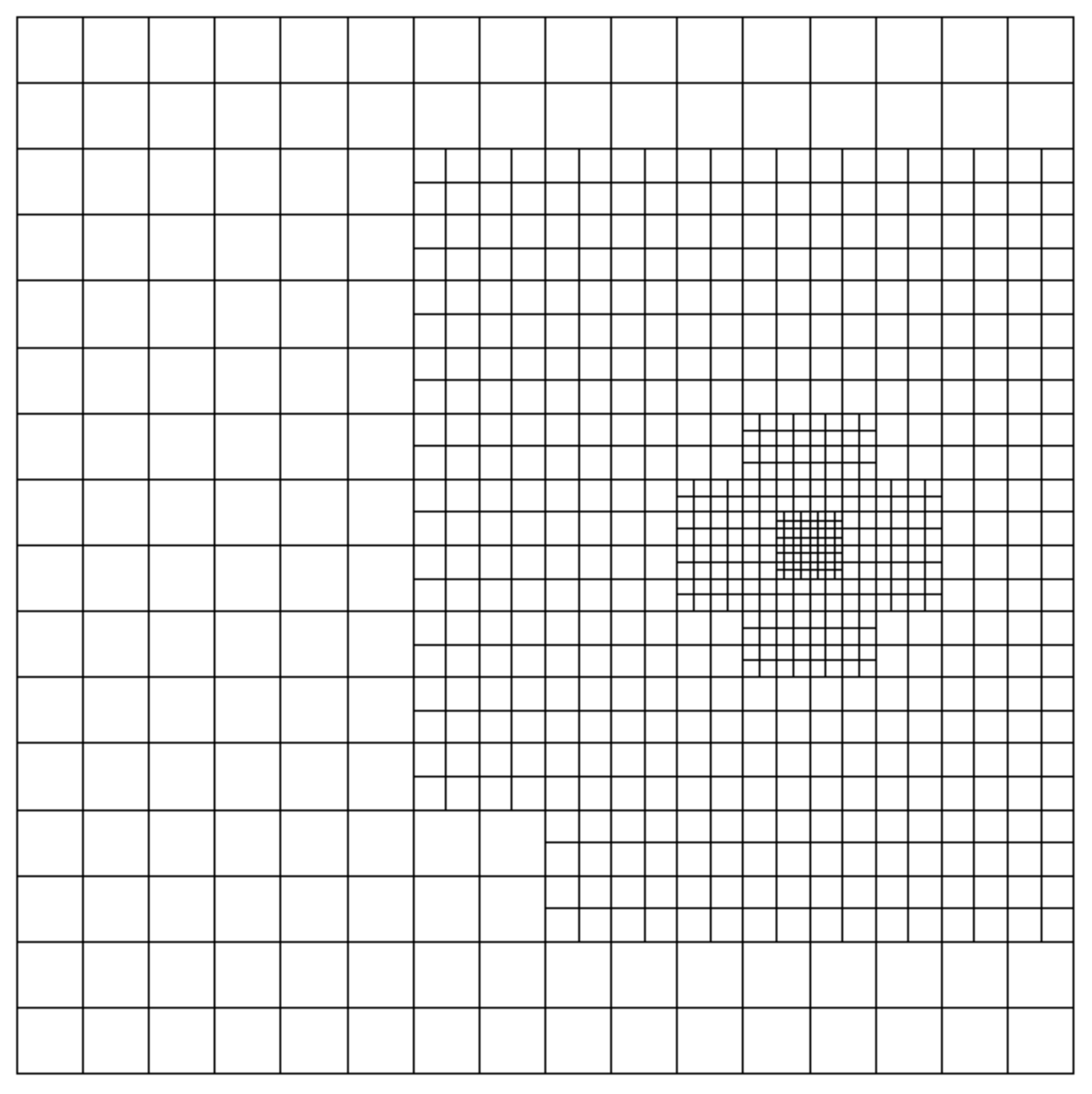}
\end{minipage}
\label{fig:6:Transport00}
}
\subfloat[$t=0.25$]{
  \centering
\begin{minipage}{.23\linewidth}
\centering
\includegraphics[width=3.5cm]{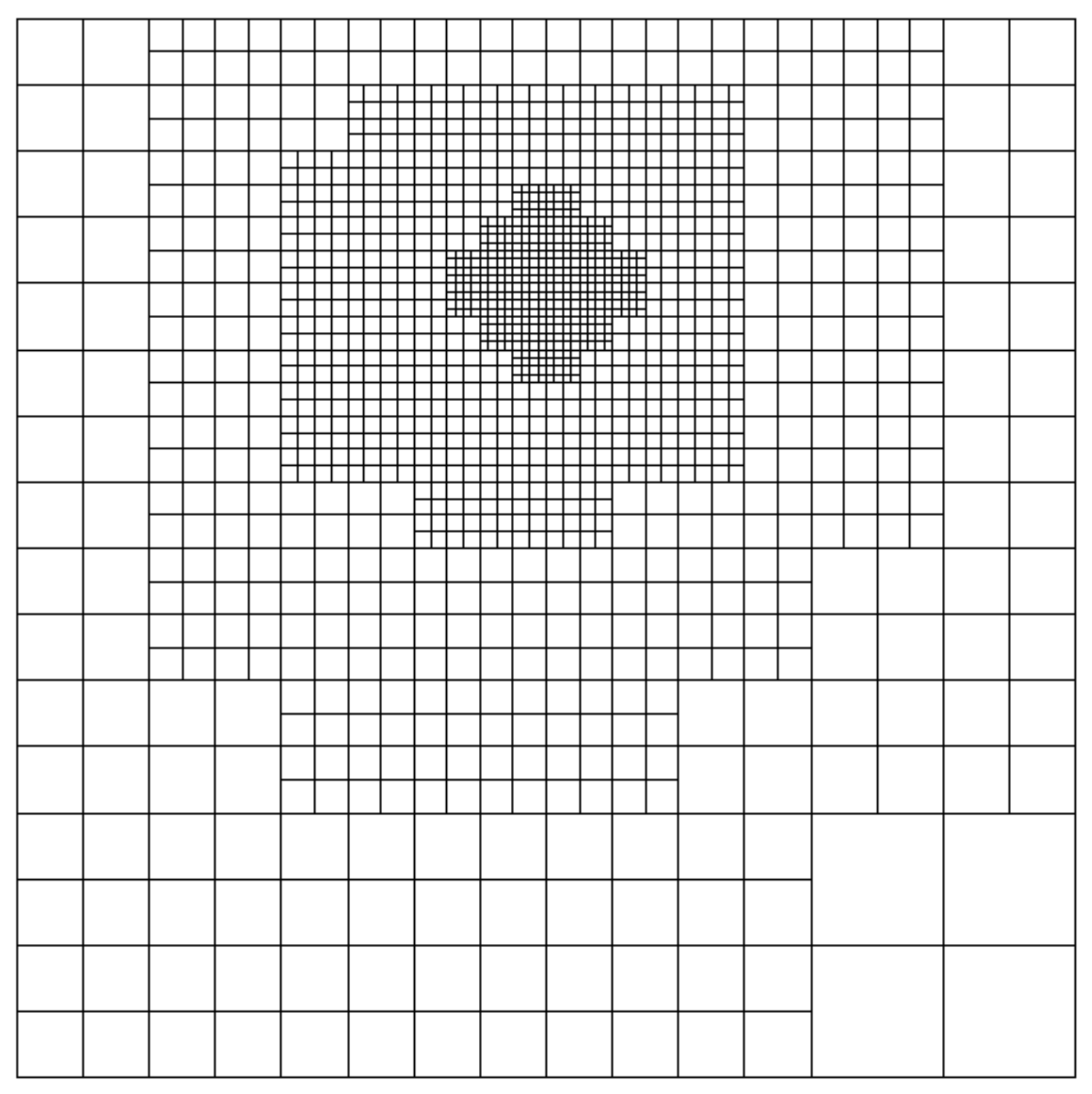}
\end{minipage}
\label{fig:6:Transport25}
}
\subfloat[$t=0.75$]{
  \centering
\begin{minipage}{.23\linewidth}
\centering
\includegraphics[width=3.5cm]{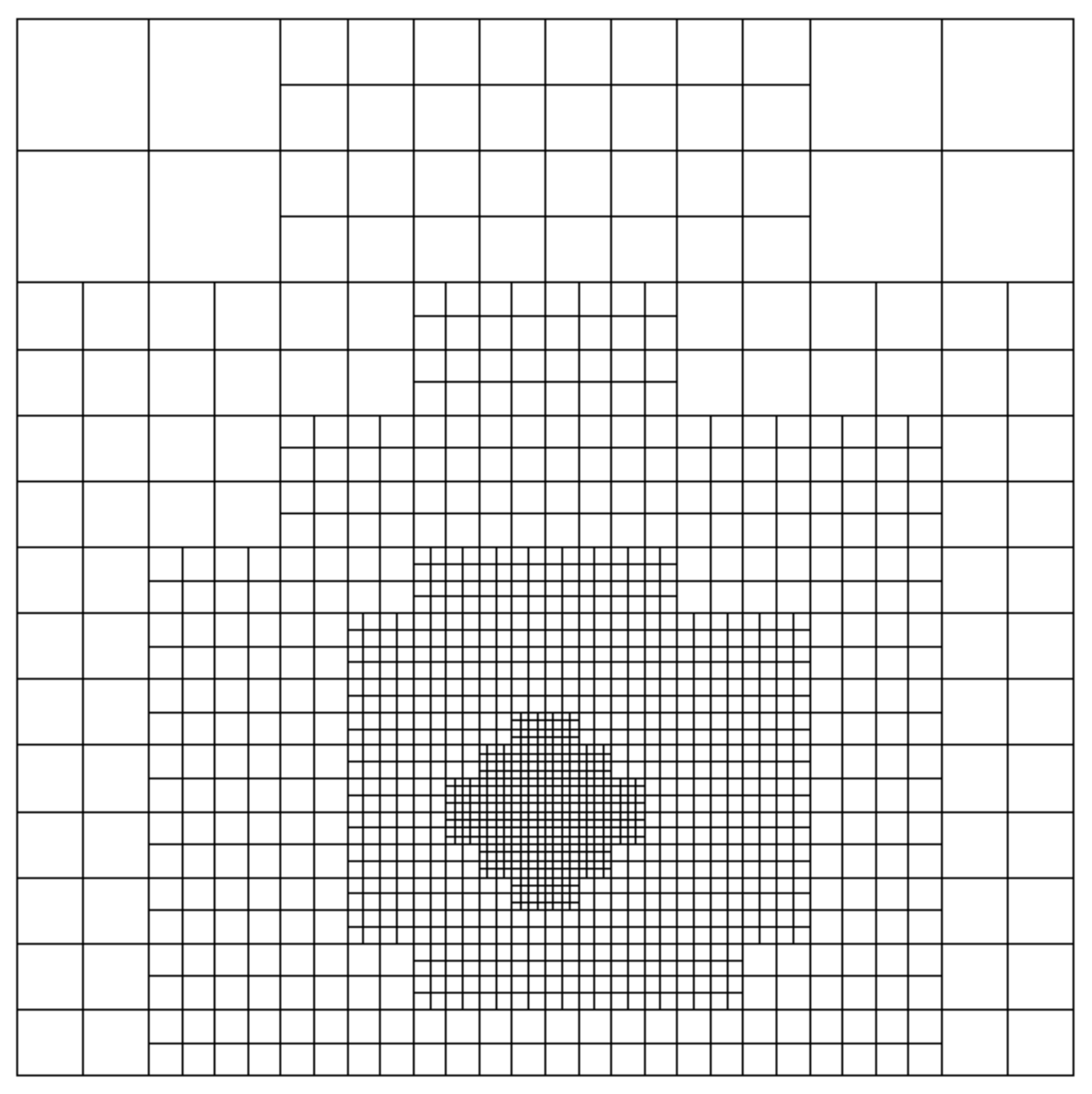}
\end{minipage}
\label{fig:6:Transport75}
}
\subfloat[$t=1.00$]{
  \centering
\begin{minipage}{.23\linewidth}
\centering
\includegraphics[width=3.5cm]{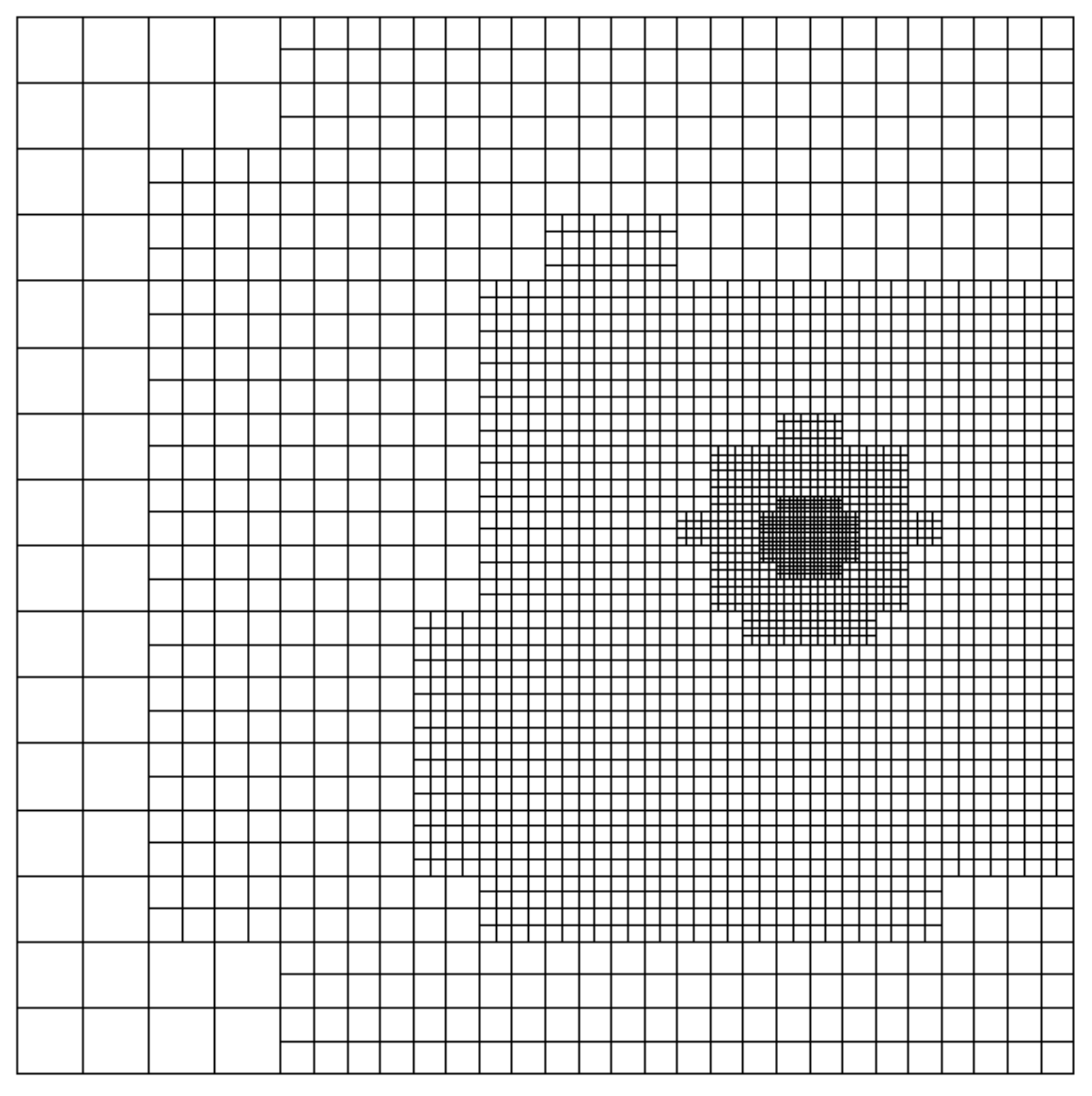}
\end{minipage}
\label{fig:6:Transport100}
}

\caption{Comparison of adaptive spatial meshes at selected specific time points 
for the flow (based on the Kelly Error Estimator) (top, (a)--(d)) and transport 
(based on the DWR method) (bottom, (e)--(h)) problem, respectively, 
corresponding to the final loop in Table~\ref{table:1:L2final-KB2BR12-2232}.
}
\label{fig:6:SpatialMeshesTransportStokes}
\end{figure}

Finally, in Fig.~\ref{fig:6:SpatialMeshesTransportStokes} we present some adaptive
spatial meshes at selected time points corresponding to the final loop in 
Table~\ref{table:1:L2final-KB2BR12-2232} for the flow and transport problem.
Considering the spatial meshes with regard to the transport problem, we note 
that the local refinements take place at the current position of the cone, i.e. 
the adaptivity runs synchronously to the rotation of the cone.
Moreover, the total number and distribution of the respective spatial cells is 
almost equal, although the refinement is slightly stronger at the final time 
point in accordance to the underlying local in time acting goal functional 
\eqref{eq:27:Goal-L2final}.
Regarding the spatial meshes of the flow problem obtained by using non-weighted 
error indicators by means of a Kelly Error Estimator, we observe an equal number
and distribution of the spatial cells over the whole time located to the course 
of the stream lines of the flow field solution, cf.~\eqref{eq:26:BRdynamic}.
This observation is in good agreement with the results obtained in 
\cite[Sec. 4.7.2]{BruchhaeuserS09} using so-called heuristic error indicators,
cf.~\cite[Sec. 4.6]{BruchhaeuserS09} for further details of this approach.
All in all, the algorithm provides very efficient spatial meshes with regard to 
the underlying goal functional, additionally taking into account the dynamics
in time.

%%%%%%%%%%%%%%%%%%%%%%%%%%%%%%%%%%%%%%%%%%%%%%%%%%%%%%%%%%%%%%%%%%%%%%%%%%%%%%%%
%% Begin Figure 10 Spatial Meshes Transport and Stokes
%%%%%%%%%%%%%%%%%%%%%%%%%%%%%%%%%%%%%%%%%%%%%%%%%%%%%%%%%%%%%%%%%%%%%%%%%%%%%%%%
%% End Figure 10 Spatial Meshes Transport and Stokes

%%%%%%%%%%%%%%%%%%%%%%%%%%%%%%%%%%%%%%%%%%%%%%%%%%%%%%%%%%%%%%%%%%%%%%%%%%%%%%%%
%% END EXAMPLE L2FINAL

\subsection{Example 2 (Transport in a Constricted Channel)}
\label{sec:5:2}

In this example, we simulate a convection-dominated transport of a species 
through a channel with a constraint. Here, we investigate the algorithm 
presented in Sec.~\ref{sec:4:3:algorithm} with regard to fully adaptive 
space-time refinements of the transport and flow problem obtained by means of 
the weighted, DWR-based error indicators \eqref{eq:24:eta_transport} and 
non-weighted, auxiliary Kelly error indicators \eqref{eq:24:auxiliary_eta_stokes},
respectively. The following results may be compared to Example 3 in \cite{Bruchhaeuser2022},
where a naive, fixed in advance adaptive refinement strategy without using 
specific error indicators was used for the flow problem. 
Hence, we study the coupled flow and transport problem with the following setting,
where we refer to \cite{Bruchhaeuser2022} for further details.

\begin{figure}%[H]
\centering

\includegraphics[width=.44\linewidth]{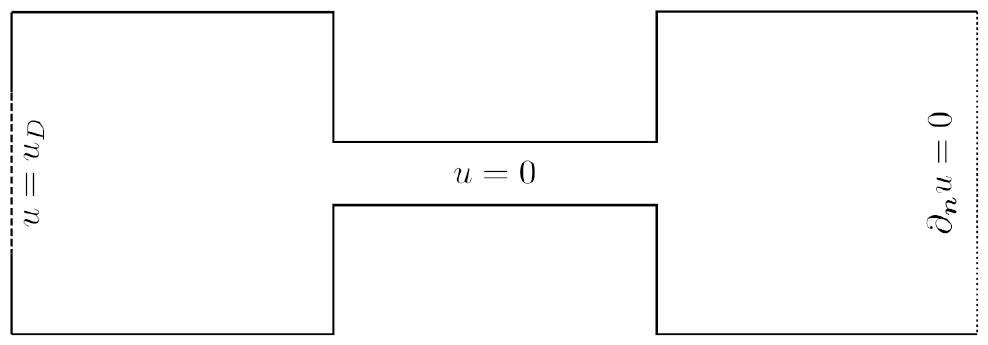}
~~
\includegraphics[width=.44\linewidth]{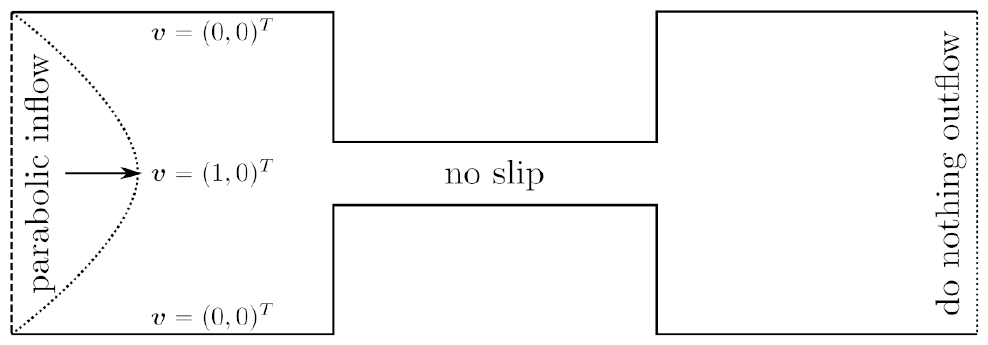}

\caption{Boundary colorization for the transport problem (left)
and the coupled flow problem (right) for Sec.~\ref{sec:5:2}.}

\label{fig:7:boundary}
\end{figure}

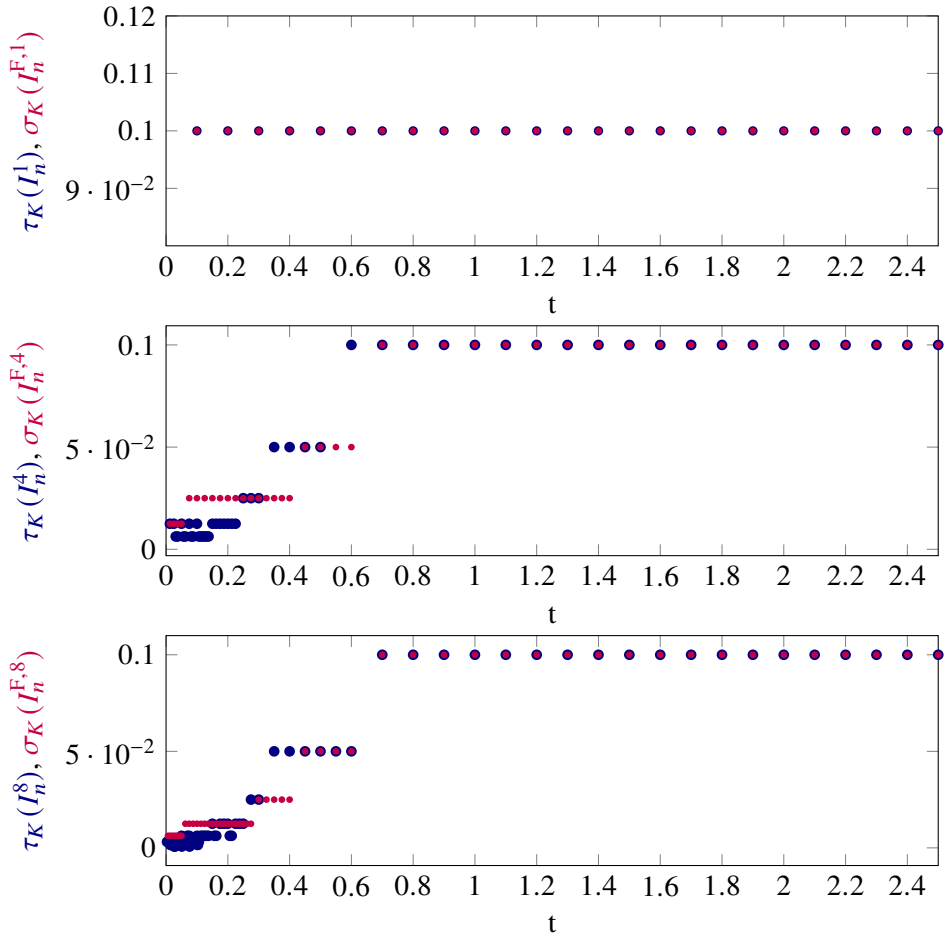
\begin{figure}[hbt!]
\begin{minipage}{.8\linewidth}
\centering
\begin{tikzpicture}
\begin{axis}[%
width=4.in,
height=1.2in,
scale only axis,
xlabel={t},
ylabel={\textcolor{navyblue}{$\tau_K(I_n^{1})$}, 
\textcolor{HSUred}{$\sigma_K(I_n^{\textnormal{F},1})$}},
xmin=0.0,
xmax=2.5,
]
\addplot [
color=navyblue,
solid,
line width=1.5pt,
mark=*,
mark size = 1.,
only marks,
mark options={solid,navyblue}
]
table[row sep=crcr]{
0.1 0.1 \\
0.2 0.1 \\
0.3 0.1 \\
0.4 0.1 \\
0.5 0.1 \\
0.6 0.1 \\
0.7 0.1 \\
0.8 0.1 \\
0.9 0.1 \\
1 0.1 \\
1.1 0.1 \\
1.2 0.1 \\
1.3 0.1 \\
1.4 0.1 \\
1.5 0.1 \\
1.6 0.1 \\
1.7 0.1 \\
1.8 0.1 \\
1.9 0.1 \\
2 0.1 \\
2.1 0.1 \\
2.2 0.1 \\
2.3 0.1 \\
2.4 0.1 \\
2.5 0.1 \\
};
\addplot [
color=HSUred,
solid,
line width=0.5pt,
mark=*,
mark size = 1.,
only marks,
mark options={fill=HSUred}
]
table[row sep=crcr]{
0.1 0.1 \\
0.2 0.1 \\
0.3 0.1 \\
0.4 0.1 \\
0.5 0.1 \\
0.6 0.1 \\
0.7 0.1 \\
0.8 0.1 \\
0.9 0.1 \\
1 0.1 \\
1.1 0.1 \\
1.2 0.1 \\
1.3 0.1 \\
1.4 0.1 \\
1.5 0.1 \\
1.6 0.1 \\
1.7 0.1 \\
1.8 0.1 \\
1.9 0.1 \\
2 0.1 \\
2.1 0.1 \\
2.2 0.1 \\
2.3 0.1 \\
2.4 0.1 \\
2.5 0.1 \\
};
\end{axis}
\end{tikzpicture}
\end{minipage}

\begin{minipage}{.8\linewidth}
\centering
\begin{tikzpicture}
\begin{axis}[%
width=4.in,
height=1.2in,
scale only axis,
/pgf/number format/.cd, 1000 sep={},
xlabel={t},
ylabel={\textcolor{navyblue}{$\tau_K(I_n^{4})$}, 
\textcolor{HSUred}{$\sigma_K(I_n^{\textnormal{F},4})$}},
xmin=0.0,
xmax=2.5,
]

\addplot [
color=navyblue,
solid,
line width=1.0pt,
mark=*,
mark size = 1.5,
only marks,
mark options={solid,navyblue}
]
table[row sep=crcr]{
0.0125 0.0125 \\
0.025 0.0125 \\
0.03125 0.00625 \\
0.0375 0.00625 \\
0.05 0.0125 \\
0.05625 0.00625 \\
0.0625 0.00625 \\
0.075 0.0125 \\
0.08125 0.00625 \\
0.0875 0.00625 \\
0.1 0.0125 \\
0.10625 0.00625 \\
0.1125 0.00625 \\
0.11875 0.00625 \\
0.125 0.00625 \\
0.13125 0.00625 \\
0.1375 0.00625 \\
0.15 0.0125 \\
0.1625 0.0125 \\
0.175 0.0125 \\
0.1875 0.0125 \\
0.2 0.0125 \\
0.2125 0.0125 \\
0.225 0.0125 \\
0.25 0.025 \\
0.275 0.025 \\
0.3 0.025 \\
0.35 0.05 \\
0.4 0.05 \\
0.45 0.05 \\
0.5 0.05 \\
0.6 0.1 \\
0.7 0.1 \\
0.8 0.1 \\
0.9 0.1 \\
1 0.1 \\
1.1 0.1 \\
1.2 0.1 \\
1.3 0.1 \\
1.4 0.1 \\
1.5 0.1 \\
1.6 0.1 \\
1.7 0.1 \\
1.8 0.1 \\
1.9 0.1 \\
2 0.1 \\
2.1 0.1 \\
2.2 0.1 \\
2.3 0.1 \\
2.4 0.1 \\
2.5 0.1 \\
};

\addplot [
color=HSUred,
solid,
line width=0.5pt,
mark=*,
mark size = 1.,
only marks,
mark options={fill=HSUred}
]
table[row sep=crcr]{
0.0125 0.0125 \\
0.025 0.0125 \\
0.0375 0.0125 \\
0.05 0.0125 \\
0.075 0.025 \\
0.1  0.025 \\
0.125 0.025 \\
0.15 0.025 \\
0.175 0.025 \\
0.2 0.025 \\
0.225 0.025 \\
0.25 0.025 \\
0.275 0.025 \\
0.3 0.025 \\
0.325 0.025 \\
0.35 0.025 \\
0.375 0.025 \\
0.4 0.025 \\
0.45 0.05 \\
0.5 0.05 \\
0.55 0.05 \\
0.6 0.05 \\
0.7 0.1 \\
0.8 0.1 \\
0.9 0.1 \\
1 0.1 \\
1.1 0.1 \\
1.2 0.1 \\
1.3 0.1 \\
1.4 0.1 \\
1.5 0.1 \\
1.6 0.1 \\
1.7 0.1 \\
1.8 0.1 \\
1.9 0.1 \\
2 0.1 \\
2.1 0.1 \\
2.2 0.1 \\
2.3 0.1 \\
2.4 0.1 \\
2.5 0.1 \\
};
\end{axis}
\end{tikzpicture}
\end{minipage}
\begin{minipage}{.8\linewidth}
\centering
\begin{tikzpicture}
\begin{axis}[%
width=4.in,
height=1.2in,
scale only axis,
/pgf/number format/.cd, 1000 sep={},
xlabel={t},
ylabel={\textcolor{navyblue}{$\tau_K(I_n^{8})$}, 
\textcolor{HSUred}{$\sigma_K(I_n^{\textnormal{F},8})$}},
xmin=0.0,
xmax=2.5,
]

\addplot [
color=navyblue,
solid,
line width=1.0pt,
mark=*,
mark size = 1.5,
only marks,
mark options={solid,navyblue}
]
table[row sep=crcr]{
0.003125 0.003125 \\
0.00625 0.003125 \\
0.009375 0.003125 \\
0.0125 0.003125 \\
0.0140625 0.0015625 \\
0.015625 0.0015625 \\
0.0171875 0.0015625 \\
0.01875 0.0015625 \\
0.0203125 0.0015625 \\
0.021875 0.0015625 \\
0.0234375 0.0015625 \\
0.025 0.0015625 \\
0.0257813 0.00078125 \\
0.0265625 0.00078125 \\
0.028125 0.0015625 \\
0.0289063 0.00078125 \\
0.0296875 0.00078125 \\
0.03125 0.0015625 \\
0.0328125 0.0015625 \\
0.034375 0.0015625 \\
0.0375 0.003125 \\
0.040625 0.003125 \\
0.04375 0.003125 \\
0.05 0.00625 \\
0.0507812 0.00078125 \\
0.0515625 0.00078125 \\
0.053125 0.0015625 \\
0.05625 0.003125 \\
0.059375 0.003125 \\
0.0625 0.003125 \\
0.06875 0.00625 \\
0.075 0.00625 \\
0.0757813 0.00078125 \\
0.0765625 0.00078125 \\
0.078125 0.0015625 \\
0.08125 0.003125 \\
0.084375 0.003125 \\
0.0875 0.003125 \\
0.090625 0.003125 \\
0.09375 0.003125 \\
0.1 0.00625 \\
0.101562 0.0015625 \\
0.103125 0.0015625 \\
0.10625 0.003125 \\
0.1125 0.00625 \\
0.11875 0.00625 \\
0.125 0.00625 \\
0.13125 0.00625 \\
0.1375 0.00625 \\
0.15 0.0125 \\
0.15625 0.00625 \\
0.1625 0.00625 \\
0.175 0.0125 \\
0.1875 0.0125 \\
0.2 0.0125 \\
0.20625 0.00625 \\
0.2125 0.00625 \\
0.225 0.0125 \\
0.2375 0.0125 \\
0.25 0.0125 \\
0.275 0.025 \\
0.3 0.025 \\
0.35 0.05 \\
0.4 0.05 \\
0.45 0.05 \\
0.5 0.05 \\
0.55 0.05 \\
0.6 0.05 \\
0.7 0.1 \\
0.8 0.1 \\
0.9 0.1 \\
1 0.1 \\
1.1 0.1 \\
1.2 0.1 \\
1.3 0.1 \\
1.4 0.1 \\
1.5 0.1 \\
1.6 0.1 \\
1.7 0.1 \\
1.8 0.1 \\
1.9 0.1 \\
2 0.1 \\
2.1 0.1 \\
2.2 0.1 \\
2.3 0.1 \\
2.4 0.1 \\
2.5 0.1 \\
};

\addplot [
color=HSUred,
solid,
line width=0.5pt,
mark=*,
mark size = 1.,
only marks,
mark options={fill=HSUred}
]
table[row sep=crcr]{
0.00625 0.00625 \\
0.0125 0.00625 \\
0.01875 0.00625 \\
0.025 0.00625 \\
0.03125 0.00625 \\
0.0375 0.00625 \\
0.04375 0.00625 \\
0.05 0.00625 \\
0.0625 0.0125 \\
0.075 0.0125 \\
0.0875 0.0125 \\
0.1 0.0125 \\
0.1125 0.0125 \\
0.125 0.0125 \\
0.1375 0.0125 \\
0.15 0.0125 \\
0.1625 0.0125 \\
0.175 0.0125 \\
0.1875 0.0125 \\
0.2 0.0125 \\
0.2125 0.0125 \\
0.225 0.0125 \\
0.2375 0.0125 \\
0.25 0.0125 \\
0.2625 0.0125 \\
0.275 0.0125 \\
0.3 0.025 \\
0.325 0.025 \\
0.35 0.025 \\
0.375 0.025 \\
0.4 0.025 \\
0.45 0.05 \\
0.5 0.05 \\
0.55 0.05 \\
0.6 0.05 \\
0.7 0.1 \\
0.8 0.1 \\
0.9 0.1 \\
1 0.1 \\
1.1 0.1 \\
1.2 0.1 \\
1.3 0.1 \\
1.4 0.1 \\
1.5 0.1 \\
1.6 0.1 \\
1.7 0.1 \\
1.8 0.1 \\
1.9 0.1 \\
2 0.1 \\
2.1 0.1 \\
2.2 0.1 \\
2.3 0.1 \\
2.4 0.1 \\
2.5 0.1 \\
};
\end{axis}
\end{tikzpicture}
\end{minipage}

\caption{Distribution of the temporal step size $\tau_K$ of the transport problem
(based on the DWR method) and $\sigma_K$ of the Stokes flow problem (based on the 
Kelly Error Estimator) over the time interval $I=(0,T]$ for the initial (1) and
after 5 and 8 DWR-loops.}
\label{fig:9:DistributionTauSigma-stabilizedTransport}
\end{figure}

The domain and its boundary colorization are presented by Fig.~\ref{fig:7:boundary}, 
cf. also Rem.~\ref{rem:1:boundary}. 
Precisely, the spatial domain is composed of two unit squares and a constraint 
in the middle which restricts the channel height by a factor of 5.
More precisely, $\Omega = (-1,0)\times(-0.5,0.5) \cup (0,1)\times(-0.1,0.1) \cup
(1,2)\times(-0.5,0.5)$ with an initial cell diameter of $h=\sqrt{2 \cdot 0.025^2}$.
The time domain is set to $I=(0,2.5)$.
With regard to the characteristic times of the two subproblems defined in 
\eqref{eq:3:characteristic-times}, the coefficients are chosen in such a way 
that the convective part of $t_{\textnormal{transport}}$ becomes dominant towards
the diffusive and reactive part and thus there holds $t_{\textnormal{transport}}
\approx t_{\textnormal{flow}}$.
Thus, the time domain $I$ is here discretized using the same initial 
$\sigma=\tau=0.1$ for the flow and transport problem within the first 
loop $\ell=1$, cf. the first plot in 
Fig.~\ref{fig:9:DistributionTauSigma-stabilizedTransport}. 
On the left boundary $\Gamma_{\textnormal{in}}$ a time-dependent inflow profile in
the positive $x_1$-direction is prescribed for the flow field $\mathbf{v}_D$ 
given by 
\begin{equation}
\label{eq:30:insta-inflow-condition}
\mathbf{v}_D(\boldsymbol{x},t) =
\begin{cases}
 \frac{\arctan(t)}{\pi/2}\cdot(1-4x_2^2,0)^\top & \textnormal{ for } 0 \leq t \leq 0.1\,,\\
 (1,0)^\top & \textnormal{ for } 0.1 < t \leq T\,.
\end{cases}
\end{equation}
Moreover, for the transport problem, the Dirichlet boundary function value is 
homogeneous on $\Gamma_D$ except for the line $(-1,-1) \times (-0.4,0.4)$ and 
time $0 \leq t \leq 0.1$ where the constant value
\begin{displaymath}
u(\boldsymbol{x},t)=1
\end{displaymath}
is prescribed on the solution. 
The diffusion coefficient has the constant and small value of $\varepsilon=10^{-4}$ 
and the reaction coefficient is chosen $\alpha=0.1$. 
The local SUPG stabilization coefficient is here set to $\delta_K = 
\delta_0 \cdot h_K$, $\delta_0=0$, i.e. a vanishing stabilization here.
The initial value function $u_0=0$ as well as the forcing term
$g=0$ are homogeneous.
The viscosity is set to $\nu=1$.
The goal functional is
\begin{equation}
\label{eq:31:goal-mean}
J(u)= \frac{1}{T\cdot |\Omega|}
\displaystyle\int_I\int_\Omega u(\bold x, t)\, \mathrm{d}\bold{x}\mathrm{d}t\,.
\end{equation}
Finally, the tuning parameters with regard to the adaptive refinement process are chosen
here as
\begin{displaymath}
\begin{array}{l@{\,}l@{\,}l@{\,}}
\theta_{h,1}^{\trans,\textnormal{top}} \geq \theta_{h,2}^{\trans,\textnormal{top}} 
=  \frac{1}{2} \cdot
\min\left\{ \frac{|\eta_h^{\trans}|}{|\eta_h^{\trans}| + |\eta_\tau^{\trans}|} \,, 1\right\},
& \hspace{0.2cm}
\theta_h^{\trans,\textnormal{bottom}} = 0.02\,, 
& \hspace{0.2cm}
\theta_\tau^\textnormal{top} = \frac{1}{2} \cdot
\min\left\{ \frac{|\eta_\tau^{\trans}|}{|\eta_h^{\trans}| + |\eta_\tau^{\trans}|} \,, 1\right\},
\\[1.5ex]
\theta_{h,1}^\textnormal{\flow,top} = \theta_{h,2}^\textnormal{\flow,top} = 1.0\,,
& \hspace{0.2cm} 
\theta_h^\textnormal{\flow,bottom} = 0.0\,,  
& \hspace{0.2cm}
\theta_{\sigma}^{\textnormal{top}} = 1.0 \textnormal{ (for } 0 \leq t \leq 0.2)\,,
\\[1.5ex]
\varpi = 1.0\,,
& \hspace{0.2cm}
\omega = 3.0\,.
& \hspace{0.2cm}
\theta_{\sigma}^{\textnormal{top}} = 0.0 \textnormal{ (for } 0.2 < t \leq 2.5)\,.
\end{array}
\end{displaymath}
We approximate the primal and dual transport solutions $u$ and $z$ by means of 
a cG(1)-dG(0) and a cG(2)-dG(0) method, respectively, and the primal flow 
solution $\mathbf{u}=\{\mathbf{v},p\}$ by means of a \{cG(2)-dG(0),cG(1)-dG(0)\}
discretization.

In Fig.~\ref{fig:9:DistributionTauSigma-stabilizedTransport}, we 
visualize the distribution of the adaptively determined time cell lengths $\tau_K$ 
and $\sigma_K$ used for the transport and  flow problem, respectively,
over the whole time interval $I$ for different DWR refinement loops.
We observe an adaptive refinement in time at the beginning, consistent with the 
restriction in time of the inflow boundary conditions given above. 
The closer we get to the final time point $T$ the coarser the temporal mesh is 
chosen.
This behavior nicely brings out the feature of automatically controlled mesh
refinement within the underlying algorithm regarding dynamics in time, note that
the goal functional \eqref{eq:31:goal-mean} acts global in time
here.
These observations are in good agreement to the results obtained for the fixed
refinement strategy mentioned above in \cite[Ex.~3]{Bruchhaeuser2022}. 
This validates the underlying approach using non-weighted, auxiliary error indicators
for the flow problem in order to reduce numerical costs significantly. 
%%%%%%%%%%%%%%%%%%%%%%%%%%%%%%%%%%%%%%%%%%%%%%%%%%%%%%%%%%%%%%%%%%%%%%%%%%%%%%%%
%% Begin Figure 20 Timesteps DWR loop Transport and  Flow Ex. 7 Stabilized 
%%
%%%%%%%%%%%%%%%%%%%%%%%%%%%%%%%%%%%%%%%%%%%%%%%%%%%%%%%%%%%%%%%%%%%%%%%%%%%%%%%%
%% End Figure 20 Timesteps DWR loop Transport and Stokes Flow Ex. 7 Stabilized T

\subsection{Example 3 (Convection-Dominated Transport with Stabilization)}
\label{sec:5:3}

In a final step, we investigate our multirate approach in view of focusing on the 
interaction of stabilization techniques combined with goal-oriented error control.
For this purpose, we modify Example~2 to the case of a strongly convection-dominated 
transport problem by increasing the P\'{e}clet number by two orders of magnitude.
In this case a solely application of adaptive mesh refinement is no longer 
sufficient to capture strong gradients and avoid spurious and non-physical 
oscillations. 
Then, the transport problem additionally has to be stabilized.
Here, we compare a non-stabilized solution ($\delta_0=0$) with the case
of a SUPG stabilized solution ($\delta_0\neq0$) for the transport problem. 
This final investigation is summarized in the following setting.

We study the coupled flow and transport problem given by \eqref{eq:1:stokes_problem},
\eqref{eq:2:transport_problem} with the same setting as outlined in 
Example 2., except for the following.
The diffusion coefficient has the constant and small value of 
\begin{displaymath}
\varepsilon=10^{-6}\,.
\end{displaymath}
The transport problem is stabilized using SUPG stabilization. Therefore, the 
local SUPG stabilization parameter is set to
\begin{displaymath}
\delta_K = \delta_0 \cdot h_K\,, \delta_0=0.1\,.
\end{displaymath}
The goal functional is given by \eqref{eq:31:goal-mean}.
Finally, the tuning parameters are chosen here as
\begin{displaymath}
\begin{array}{l@{\,}l@{\,}l@{\,}}
\theta_{h,1}^{\trans,\textnormal{top}} \geq \theta_{h,2}^{\trans,\textnormal{top}}
=  \frac{1}{2} \cdot
\min\left\{ \frac{|\eta_h^{\trans}|}{|\eta_h^{\trans}| + |\eta_\tau^{\trans}|} \,, 1\right\},
& \hspace{0.2cm}
\theta_h^{\trans,\textnormal{bottom}} = 0.02\,, 
& \hspace{0.2cm}
\theta_\tau^\textnormal{top} = \frac{1}{2} \cdot
\min\left\{ \frac{|\eta_\tau^{\trans}|}{|\eta_h^{\trans}| + |\eta_\tau^{\trans}|} \,, 1\right\},
\\[1.5ex]
\theta_{h,1}^\textnormal{\flow,top} = \theta_{h,2}^\textnormal{\flow,top} = 0.33\,,
& \hspace{0.2cm} 
\theta_h^\textnormal{\flow,bottom} = 0.02\,,  
& \hspace{0.2cm}
\theta_{\sigma}^{\textnormal{top}} = 0.2\,,
\\[1.5ex]
\varpi = 1.0\,,
& \hspace{0.2cm}
\omega = 3.0\,.
& \hspace{0.2cm}
\end{array}
\end{displaymath}
%
%%%%%%%%%%%%%%%%%%%%%%%%%%%%%%%%%%%%%%%%%%%%%%%%%%%%%%%%%%%%%%%%%%%%%%%%%%%%%%%%
%% End Example 6 Channel time-dependent
%
We approximate the primal and dual transport solution $u$ and $z$ by means of a 
cG(1)-dG(0) and cG(2)-dG(0) method, respectively, on adaptively refined meshes
in space and time based on weighted error indicators based on the DWR
method, given by Eq.~\eqref{eq:24:eta_transport}. 
However, the flow solution $\mathbf{u}=\{\mathbf{v},p\}$ is approximated
with a \{cG(2)-dG(0),cG(1)-dG(0)\} discretization on adaptively refined meshes 
in space and time based on auxiliary, non-weighted error indicators based on the 
Kelly Error Estimator, given by Eq.~\eqref{eq:24:auxiliary_eta_stokes}.

\begin{figure}[hbt!]
\centering

\includegraphics[width=.48\linewidth]{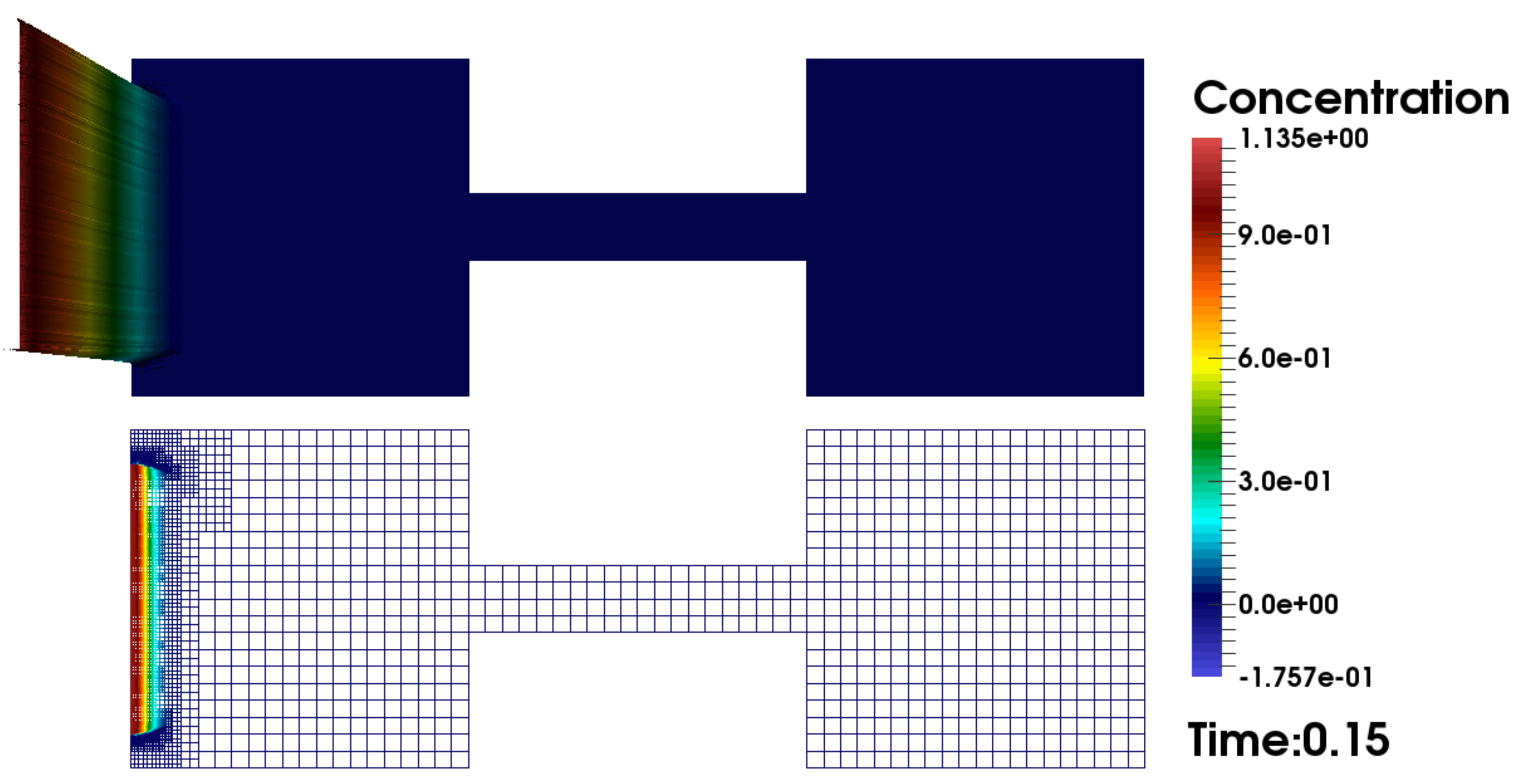}
~~
\includegraphics[width=.48\linewidth]{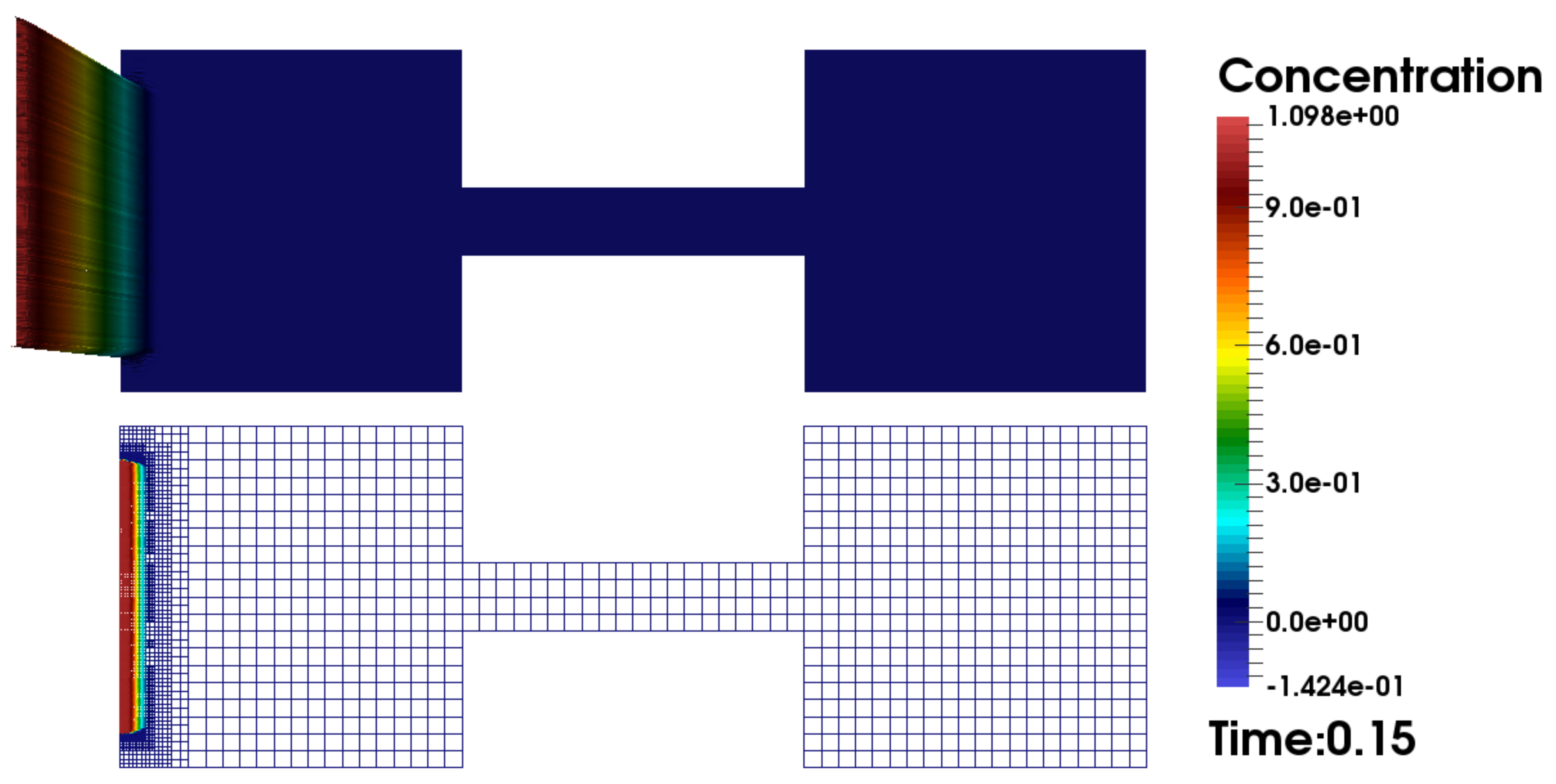}

\includegraphics[width=.48\linewidth]{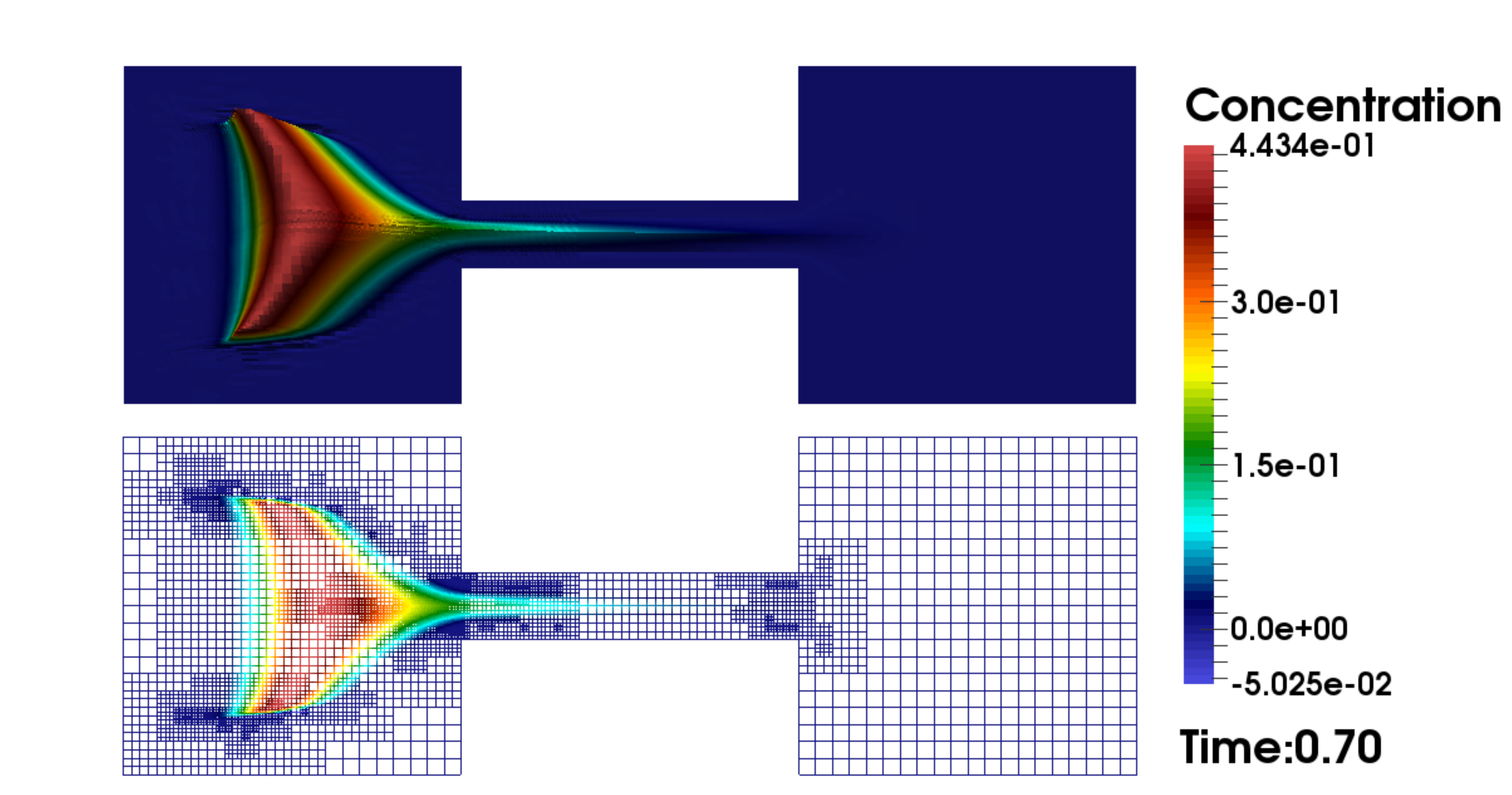}
~~
\includegraphics[width=.48\linewidth]{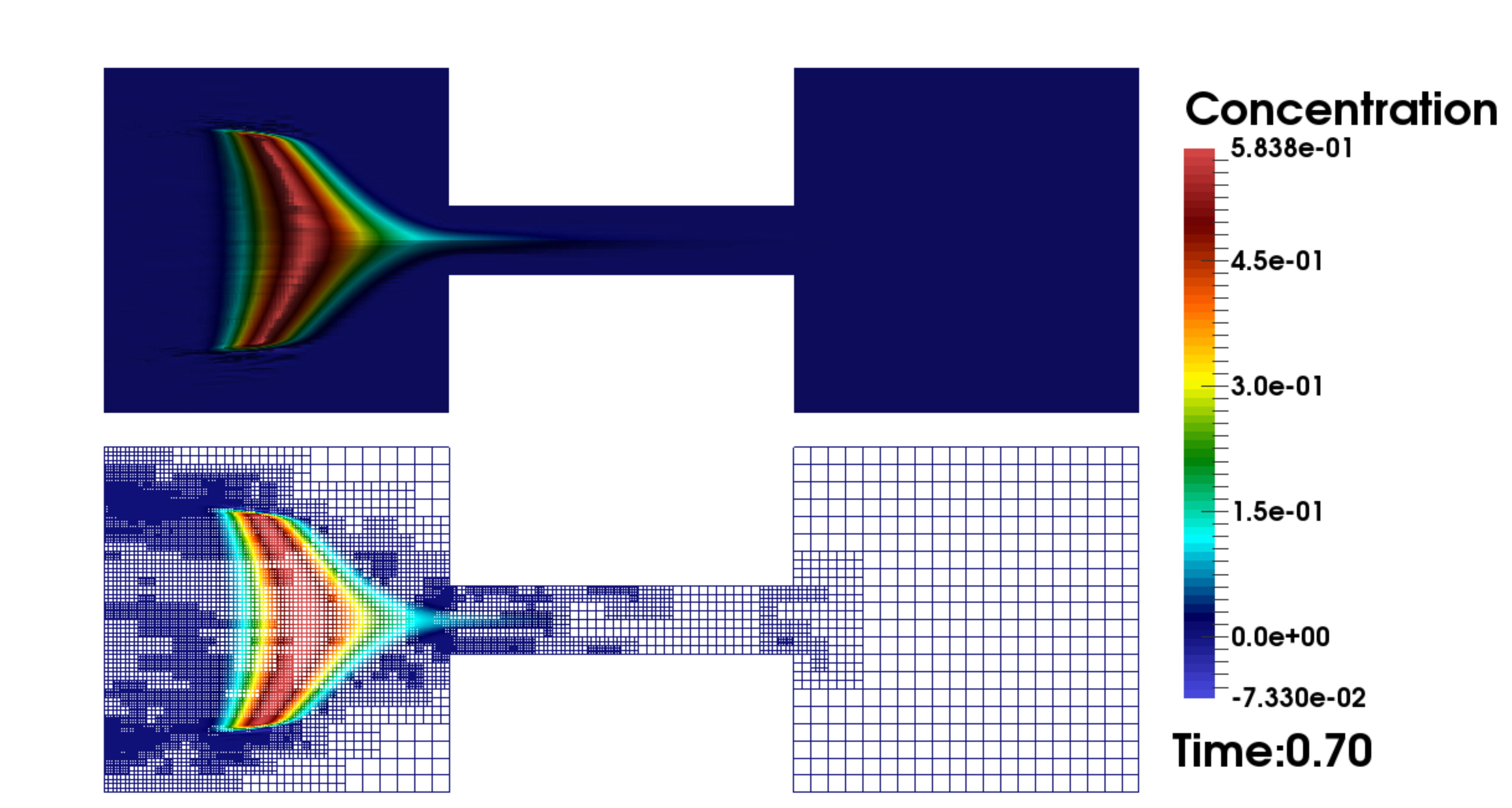}

\includegraphics[width=.48\linewidth]{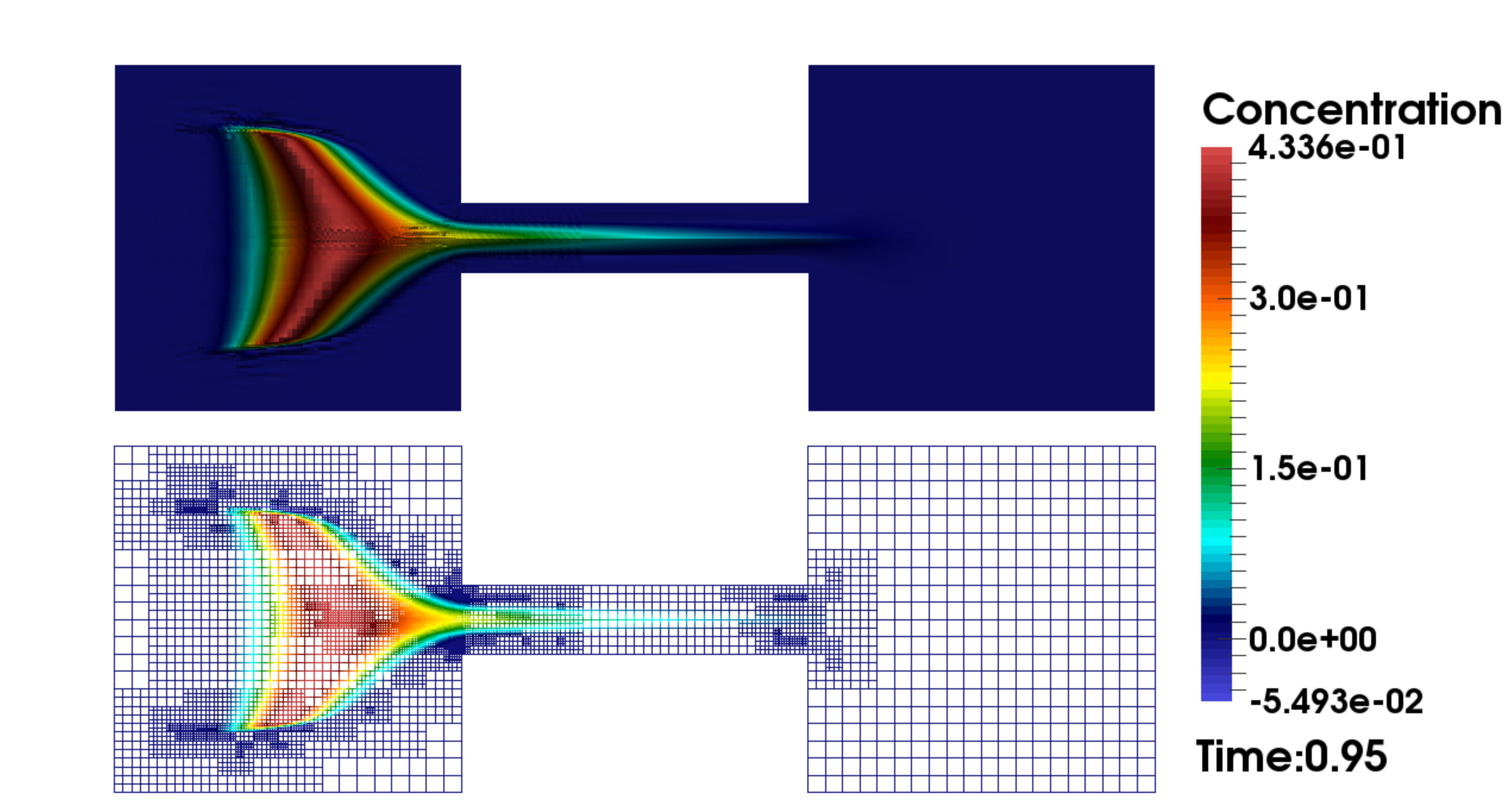}
~~
\includegraphics[width=.48\linewidth]{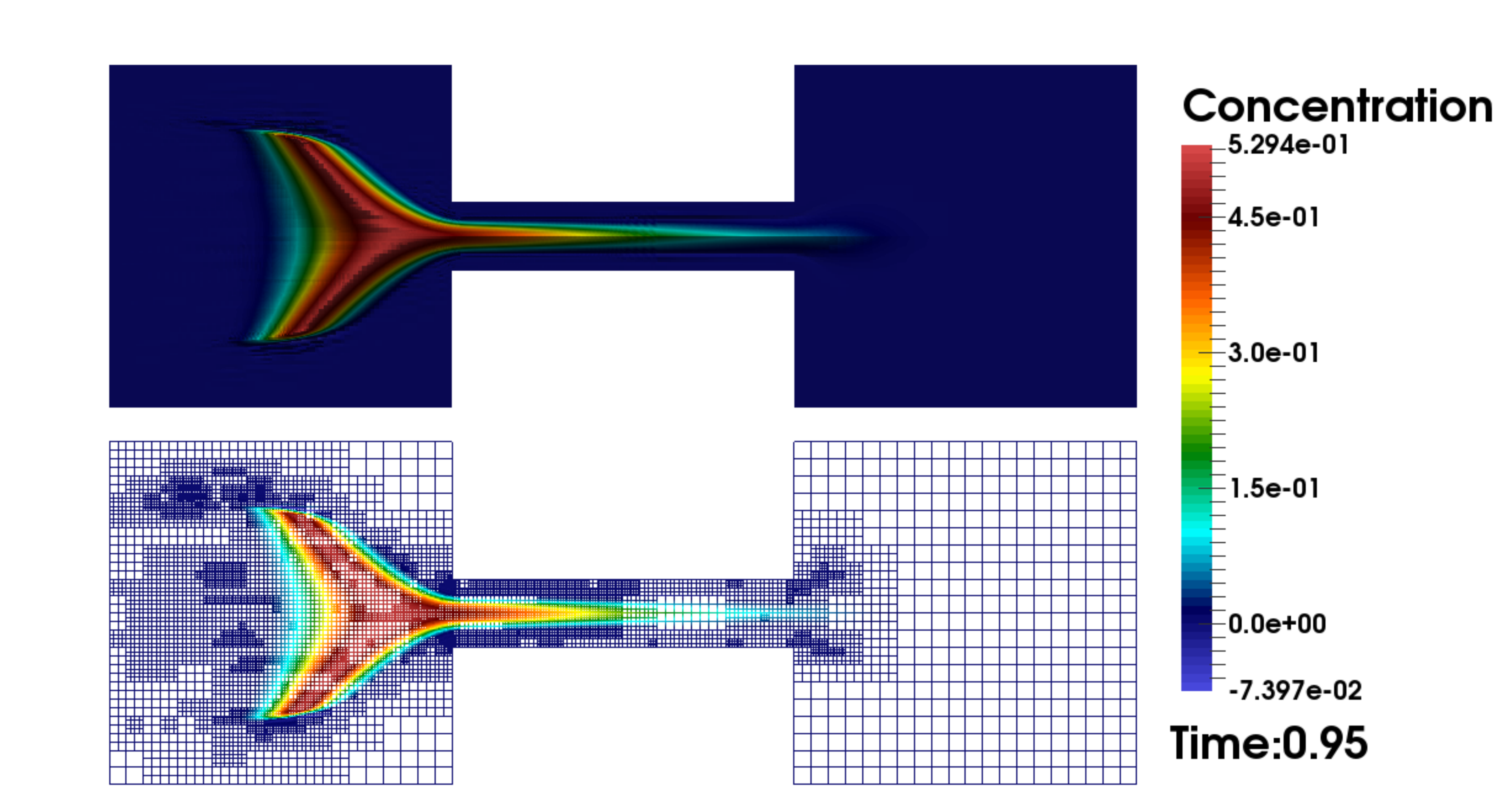}

\includegraphics[width=.48\linewidth]{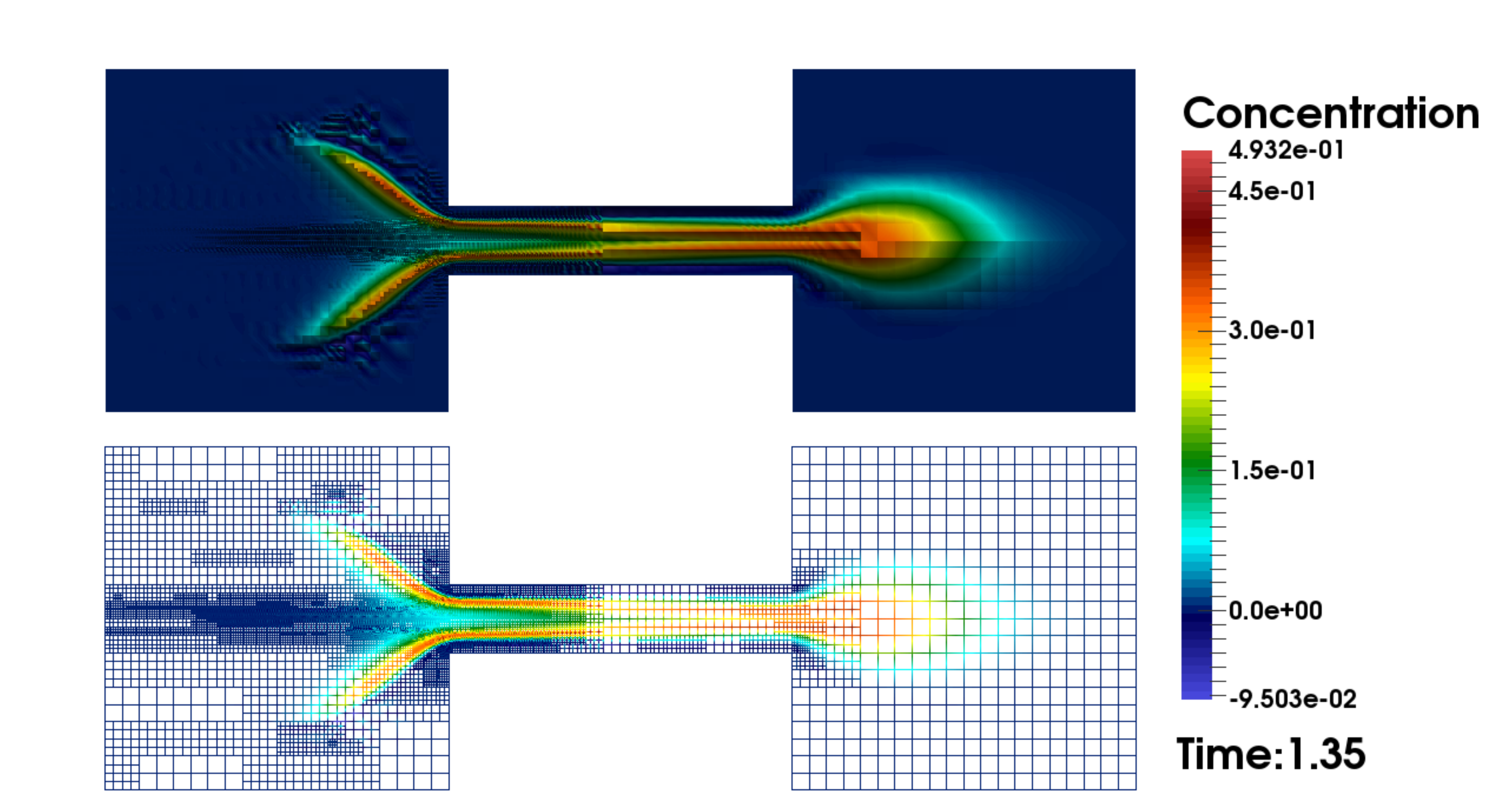}
~~
\includegraphics[width=.48\linewidth]{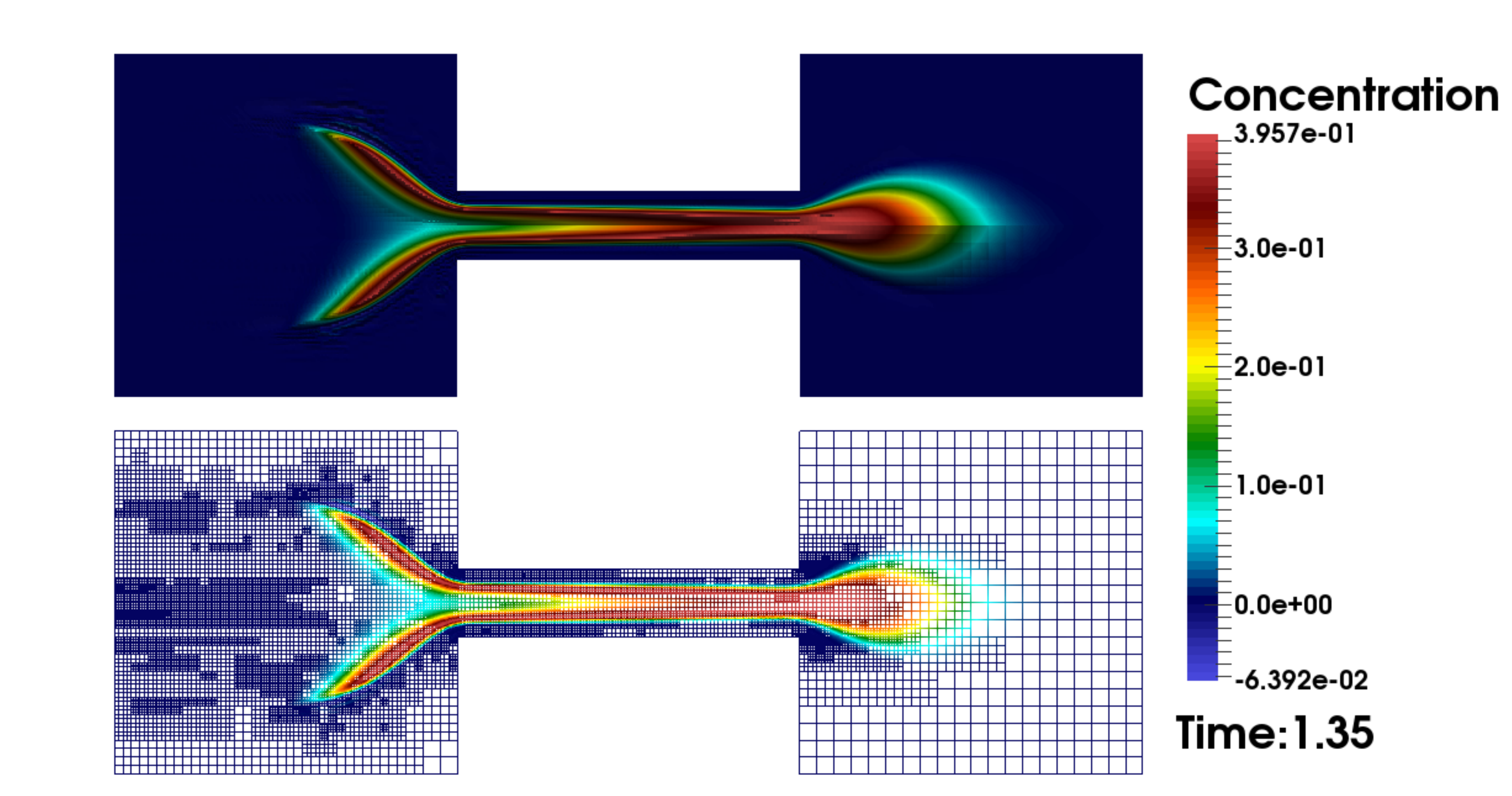}

\includegraphics[width=.48\linewidth]{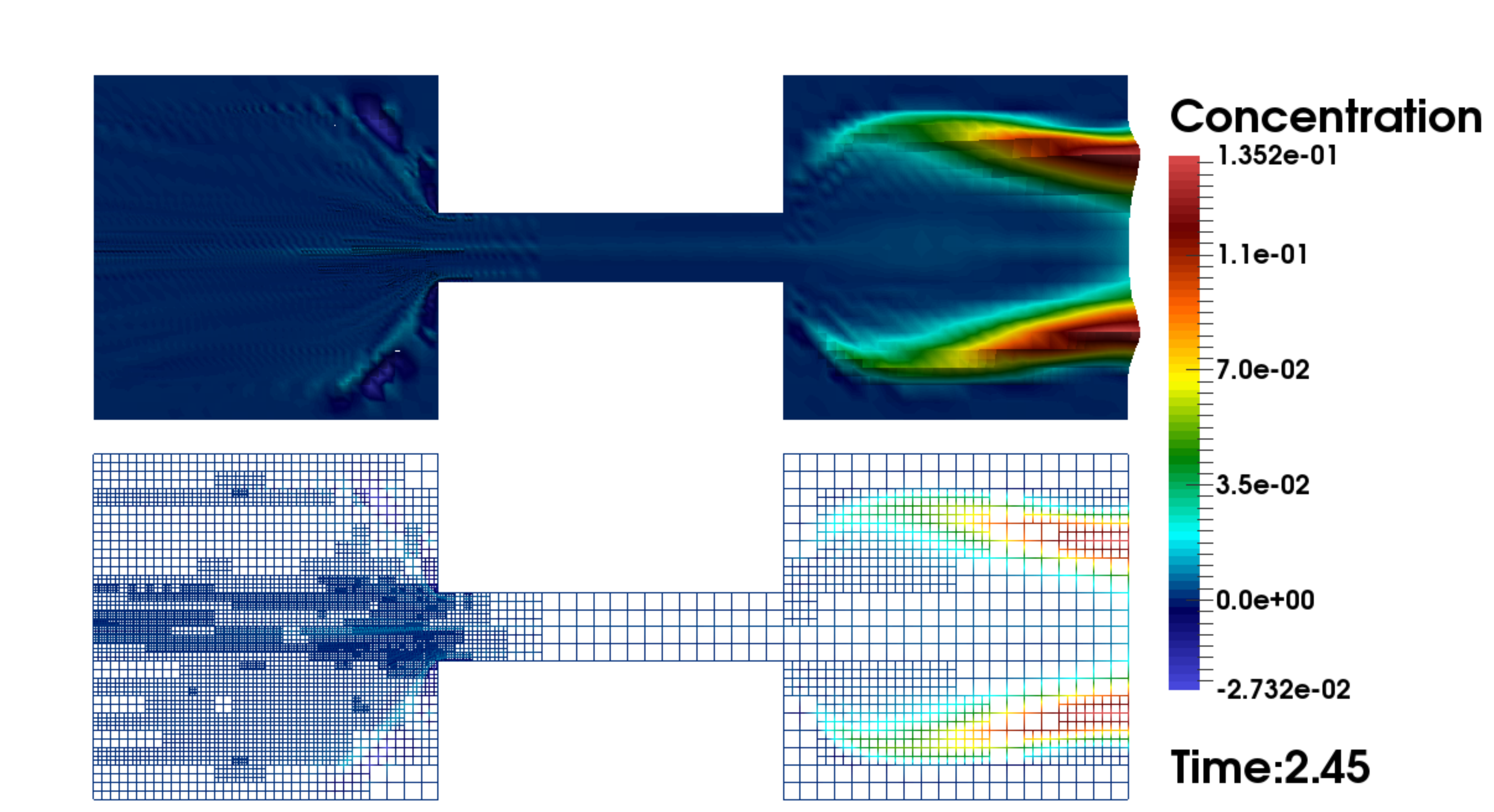}
~~
\includegraphics[width=.48\linewidth]{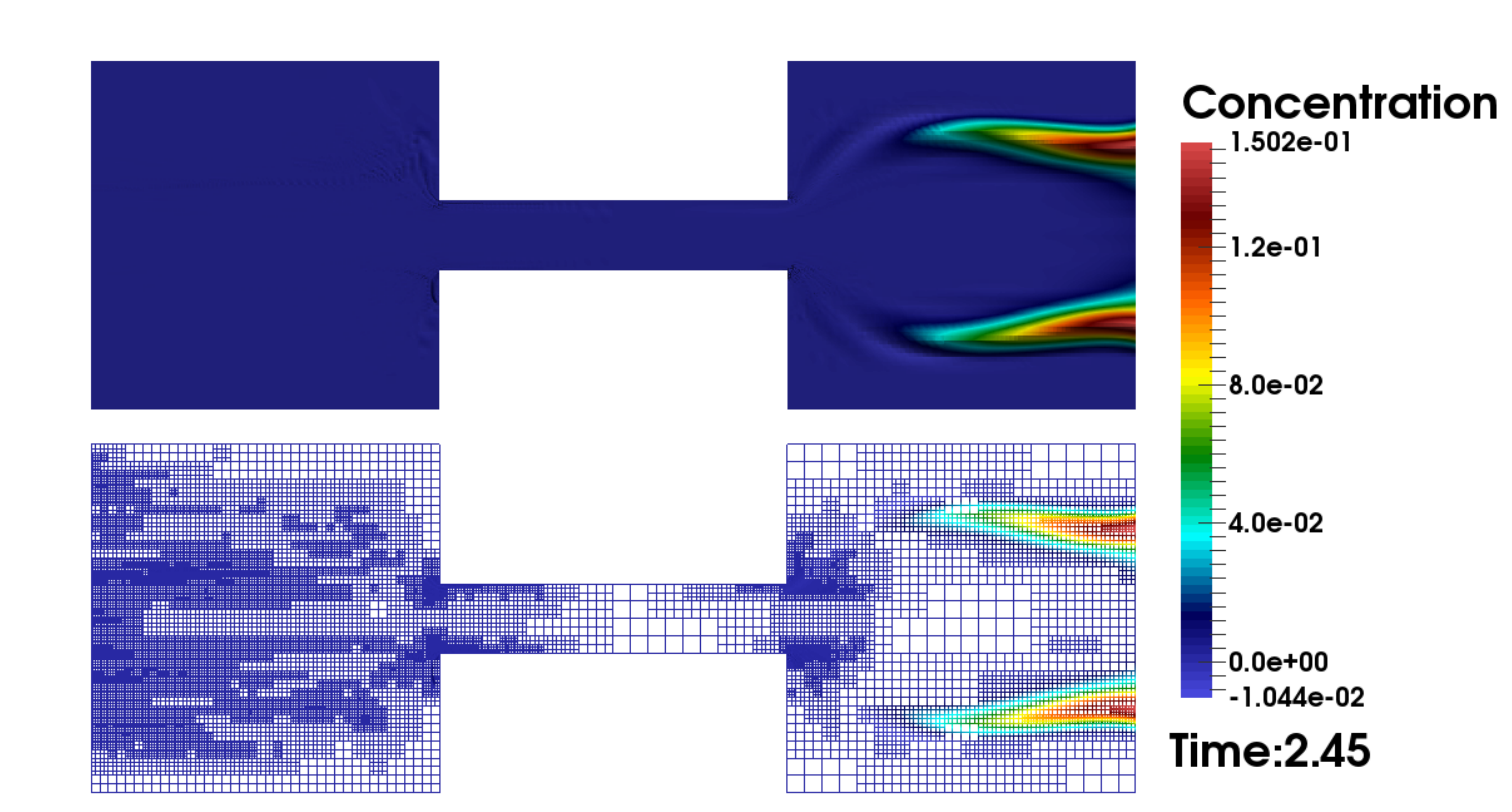}

\caption{Comparison of transport solution profiles and related 
adaptively refined spatial meshes (based on the DWR method)
without stabilization $\delta_0=0.0$ (left) 
and with SUPG stabilization $\delta_0=0.1$ (right) for $\varepsilon=10^{-6}$ 
at different time points corresponding to loop 
$\ell=8$ for Example~3.}
\label{fig:18:ComparisonDifferentTimes}
\end{figure}

\begin{figure}[hbt!]
\centering

\subfloat[No Stabilization ($\delta_0=0.0$), 29843 spatial DoFs]{
  \centering
\begin{minipage}{.95\linewidth}
\centering
\includegraphics[width=.95\linewidth]{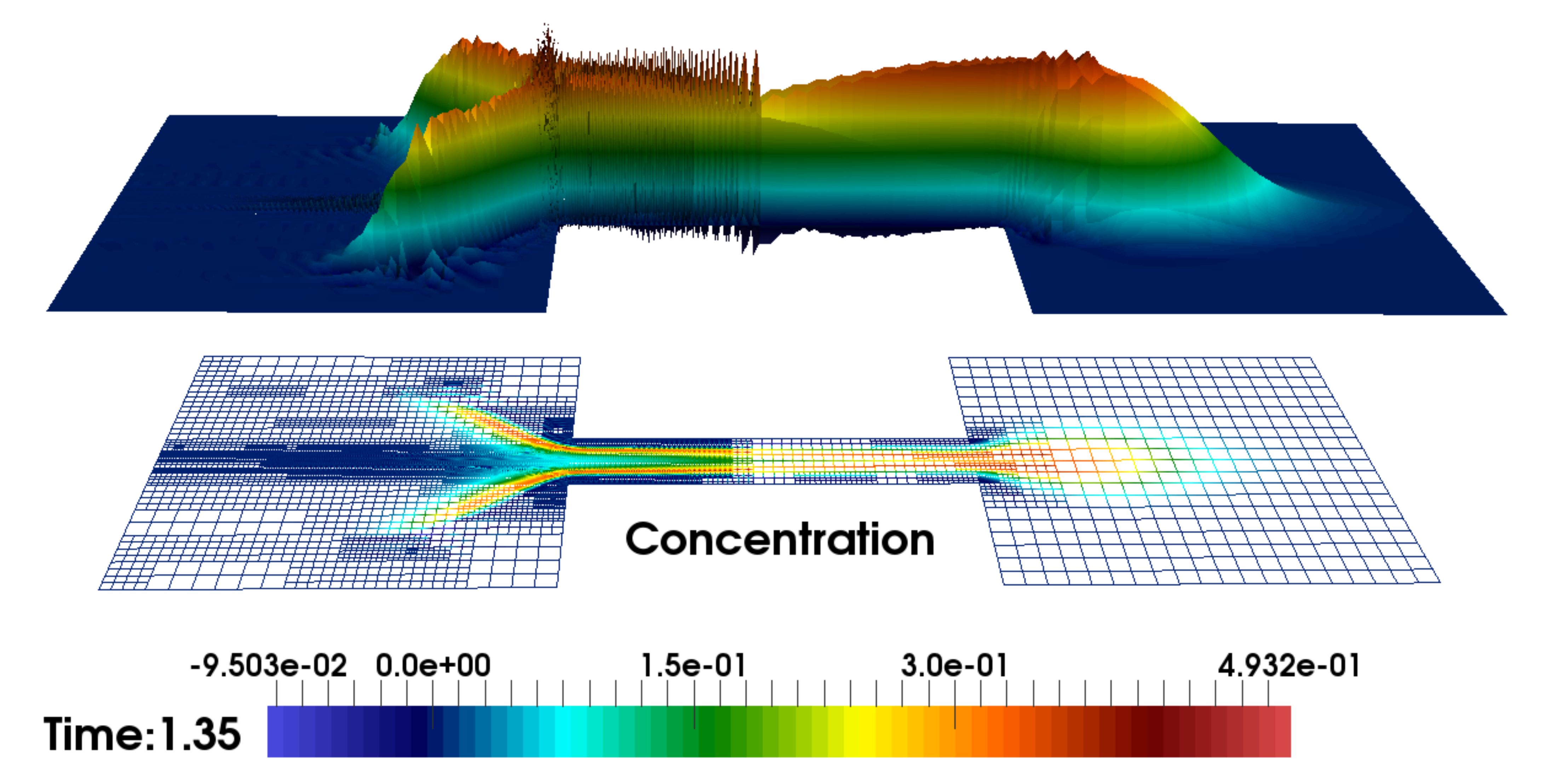}
\end{minipage}
\label{fig:?:NoStabil}
}

\subfloat[SUPG Stabilization ($\delta_0=0.1$), 27471 spatial DoFs]{
  \centering
\begin{minipage}{.95\linewidth}
\centering
\includegraphics[width=.95\linewidth]{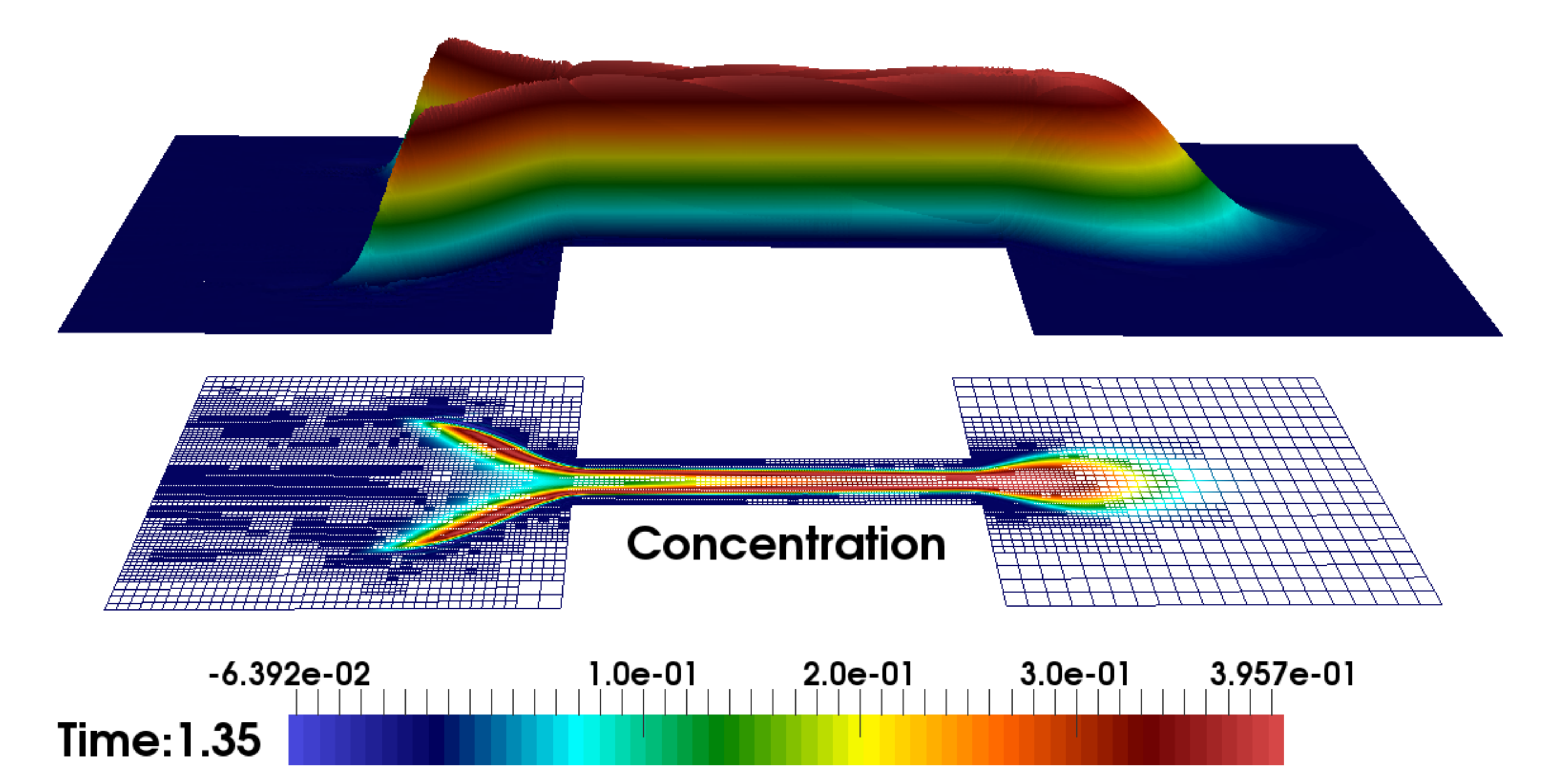}
\end{minipage}
\label{fig:?:Stabil}
}

\caption{Comparison of side profiles for the transport solution and related 
adaptively refined spatial meshes (based on the DWR method) without stabilization 
$\delta_0=0.0$ (top) and with SUPG stabilization $\delta_0=0.1$ (bottom) for 
$\varepsilon=10^{-6}$ corresponding to loop $\ell=8$ for 
Example~3, exemplary at time 
point $t=1.35$.}
\label{fig:19:ComparisonFixedTime135}
\end{figure}

In Fig.~\ref{fig:18:ComparisonDifferentTimes}, we compare the 
solution profiles and corresponding adaptive spatial meshes of the primal 
transport solution $u_{\tau h}^{1,0}$ for a non-stabilized and stabilized case, 
respectively, at selected time points within the final DWR loop $\ell=8$.
It becomes clear that in the strongly convection-dominated case a solely adaptive 
mesh refinement without stabilization is no longer sufficient to resolve the 
arising layers of the transported species within the channel, especially regarding
the solution profiles in the course of time; cf. the blurred solution profiles 
in the course of the transported species on the left part of 
Fig.~\ref{fig:18:ComparisonDifferentTimes}.
This becomes even clearer considering the exemplary side profile of the primal
transport solution at time $t=1.35$ given by the upper part of 
Fig.~\ref{fig:19:ComparisonFixedTime135} that is strongly perturbed
by means of spurious and non-physical oscillations, especially in the part of the
constriction within the channel.
Without additional stabilization techniques the underlying algorithm is not able
to capture the strong gradients and resolve the layers and sharp moving fronts
of the underlying transported species.
In contrast, regarding the stabilized solution profiles given by the right and 
lower part of Fig.~\ref{fig:18:ComparisonDifferentTimes} and 
Fig.~\ref{fig:19:ComparisonFixedTime135}, respectively, the 
solution profile fronts are resolved in a visibly more accurate way along with a 
significantly reduction of the spurious oscillations.
Moreover, regarding the underlying spatial meshes, we point out that in the 
stabilized case the adaptive refinement is located close to the whole solution
front of the underlying transported species within the channel. 
More precisely, we observe local refinements located in the left unit square 
corresponding to the wing-like fronts and behind them, within the constriction
corresponding to the course of the solution profile as well as at the exit of 
the restriction, in particular at the corners of the exit, corresponding to the 
head of the solution profile.

In contrast to that, in the non-stabilized case most of the local refinement 
takes place at the beginning of the constriction where most of the oscillations 
are visible.
In comparison, the remaining parts along the solution front are less refined, 
in particular, regarding the wing-like fronts in the left unit square.
This is obvious since the goal functional \eqref{eq:31:goal-mean}
acts global in space and time and thus those parts of the error indicators are 
weighted stronger involving larger errors in the respective quantity.

Due to the additional stabilization along with a significantly reduction of the 
oscillations, the algorithm in the stabilized case is capable to distribute the 
refinement more evenly to the regions belonging to strong gradients of the 
underlying solution profile.
In summary, the stabilized solution shows a significantly improvement with regard 
to resolving layers and sharp moving fronts along with efficient underlying 
spatial meshes, even though some slight perturbations located at the course of 
the layers are still visibly.

%%%%%%%%%%%%%%%%%%%%%%%%%%%%%%%%%%%%%%%%%%%%%%%%%%%%%%%%%%%%%%%%%%%%%%%%%%%%%%%%
%%
%%%%%%%%%%%%%%%%%%%%%%%%%%%%%%%%%%%%%%%%%%%%%%%%%%%%%%%%%%%%%%%%%%%%%%%%%%%%%%%%
%%

%%%%%%%%%%%%%%%%%%%%%%%%%%%%%%%%%%%%%%%%%%%%%%%%%%%%%%%%%%%%%%%%%%%%%%%%%%%%%%%%
%%
%%%%%%%%%%%%%%%%%%%%%%%%%%%%%%%%%%%%%%%%%%%%%%%%%%%%%%%%%%%%%%%%%%%%%%%%%%%%%%%%
%%

\begin{figure}[hbt!]
\centering
\includegraphics[width=.8\linewidth]{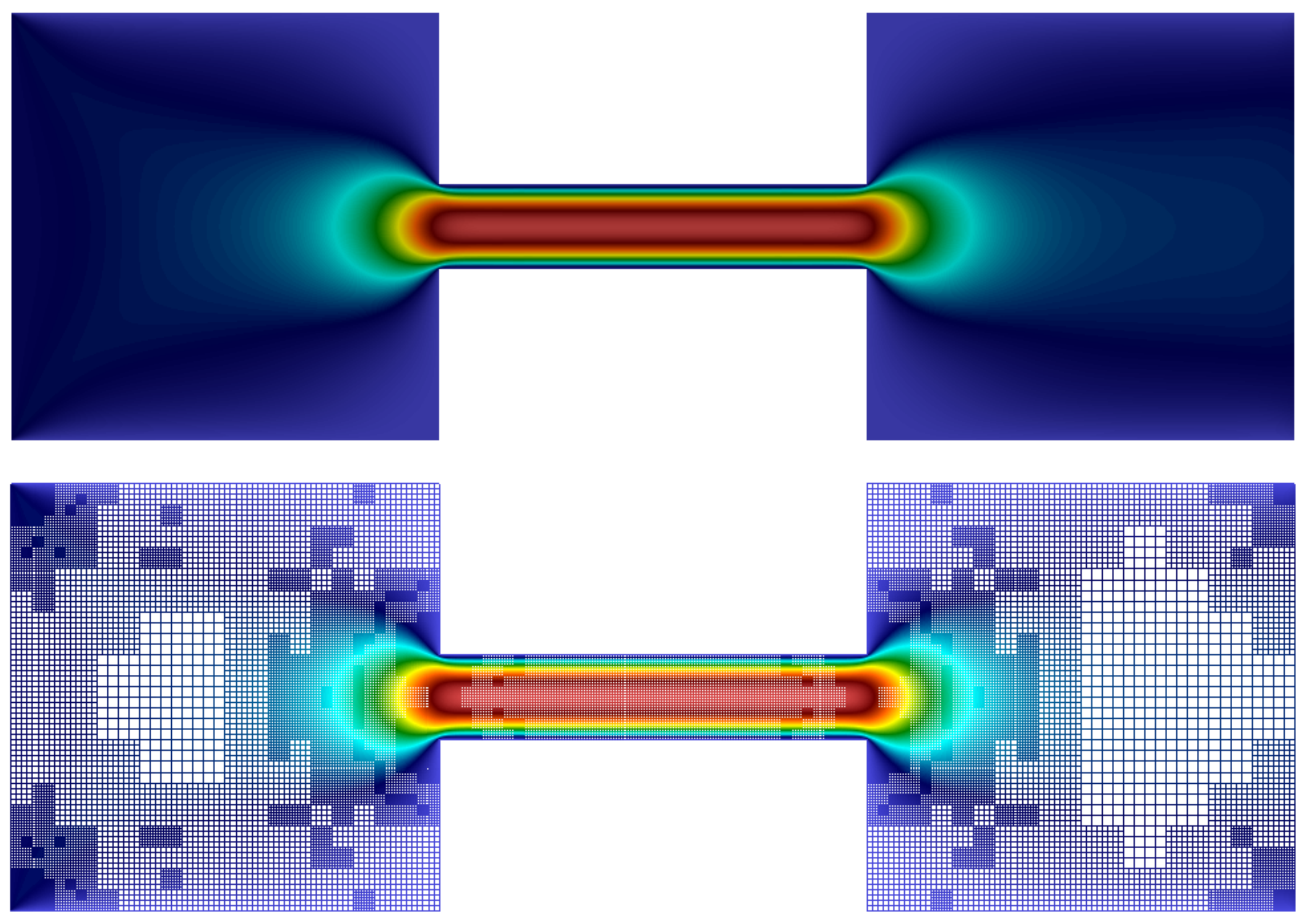}

\caption{Convection $\mathbf{v}_{\sigma h}$ solution of the Stokes flow problem 
and related adaptively refined spatial mesh (based on the Kelly Error Estimator) 
with $Q_2$-$Q_1$ finite elements and 48186 spatial DoFs, corresponding to loop $\ell=8$ for 
Example~3, exemplary at time 
point $t=0.1$.
On the left boundary a time-dependent inflow profile in the positive $x_1$-direction
is prescribed for the convection $\mathbf{v}_D$ given by 
Eq.~\eqref{eq:30:insta-inflow-condition}.}
\label{fig:21:AdaptiveMeshStokes}
\end{figure}

Finally, for the sake of completeness, we present in Fig.~\ref{fig:21:AdaptiveMeshStokes}
the solution profile and corresponding adaptive spatial mesh of the primal 
flow field solution $\mathbf{v}_{\sigma h}$ based on the Kelly Error 
Estimator, exemplary at time $t=0.1$ corresponding to the time-dependence of the
inflow boundary condition given by \eqref{eq:30:insta-inflow-condition} 
within the final loop $\ell=8$.
The adaptive spatial refinement is located to the spreading of the convection 
flow field $\mathbf{v}_{\sigma h}$, cf. the upper plot of 
Fig.~\ref{fig:21:AdaptiveMeshStokes}. 
Moreover, the spatial mesh is visibly more refined close to the corners of the 
entrance and exit of the channels' constriction consistent with occurring 
challenges arising in such regions of the underlying meshes, 
cf., e.g., \cite{Schmich2008,Endtmayer2017}.
%

%%%%%%%%%%%%%%%%%%%%%%%%%%%%%%%%%%%%%%%%%%%%%%%%%%%%%%%%%%%%%%%%%%%%%%%%%%%%%%%%
%%
%%%%%%%%%%%%%%%%%%%%%%%%%%%%%%%%%%%%%%%%%%%%%%%%%%%%%%%%%%%%%%%%%%%%%%%%%%%%%%%%

\section{Summary} 
\label{sec:6:conclusion}

In this work, we presented a cost-efficient space-time adaptive algorithm for 
our multirate approach combined with stabilization and goal-oriented error control
based on the DWR method applied to coupled flow and transport.
Different adaptive time step sizes on different temporal meshes initialized with 
the help of characteristic times were used for the two subproblems.
Furthermore, the transport problem was assumed to be convection-dominated and 
thus stabilized using the residual based SUPG method which puts an additional
facet of complexity on the algorithmic design.
Both subproblems are discretized using a discontinuous Galerkin method dG($r$) 
with an arbitrary polynomial degree $r \geq 0$ in time and a continuous Galerkin
method cG($p$) with an arbitrary polynomial degree $p \geq 1$ in space.
Goal-oriented a posteriori error representation based on the DWR method were 
derived for the flow as well as for transport problem. 
These error representations are splitted into quantities in space and time such that
their localized forms serve as error indicators for the adaptive mesh refinement 
process in space and time.
To reduce numerical costs significantly, in the numerical experiments such 
weighted error indicators were used only within the adaptivity for the transport
problem, where the adaptivity for the flow problem was achieved by using auxiliary,
non-weighted error indicators based on the Kelly Error Estimator that avoids an 
explicit computation of a dual flow problem.
Nevertheless, considering duality for the flow problem within the numerical 
examples is work that remains to be done and will be interesting to be compared 
to the results obtained in the present work.

The practical realization as well as some implementational aspects regarding 
the specifications within the multirate approach were demonstrated. 
In numerical experiments, the algorithm was studied with regard to accuracy, 
efficiency and reliability reasons by investigating a typical benchmark as well 
as a problem of practical interest. 
Spurious oscillations that typically arise in numerical approximations of 
convection-dominated problems could be reduced significantly. 
High-efficient adaptively refined meshes in space and time were obtained for 
both subproblems, where using auxiliary, non-weighted error indicators for the 
flow problem had any negative impact on the adaptive refinement process for 
the transport problem as 
effectivity indices close to one and well-balanced error indicators in space 
and time were obtained.
Moreover, the numerical results give a first hint of the potential regarding a 
temporal mesh that is adapted to the dynamics of the active or fast components 
compared to a more standard fixed-time strategy where the whole mesh
is refined due to these fast components.
Finally, the here presented approach for coupled flow and transport is fairly 
general and can be easily adopted to other multi-physics systems coupling phenomena 
that are characterized by strongly differing time scales.

%%%%%%%%%%%%%%%%%%%%%%%%%%%%%%%%%%%%%%%%%%%%%%%%%%%%%%%%%%%%%%%%%%%%%%%%%%%%%%%%
%% End own script

\section*{Acknowledgement}
We acknowledge U. K\"ocher for his support in the design and 
implementation of the underlying software \texttt{dwr-stokes-condiffrea}; 
cf. the software project \texttt{DTM++.Project/dwr} \cite{Koecher2019}.

% \section*{Funding}
% The authors have no relevant financial or non-financial interests to disclose.

\section*{Appendix}

In the Appendix, we give some detailed definitions and remarks regarding
the main results derived in Se.~\ref{sec:3:error}. 

\subsection*{\textbf{Galerkin Orthogonality for Temporal and Spatial Error of Transport
Problem}}
For the temporal error $e = \concentration-\concentration_{\tau}$
we get the following Galerkin orthogonality by subtracting 
Eq.~\eqref{eq:13:semi_transport} from Eq.~\eqref{eq:8:weak_transport} 
\begin{equation}
\label{eq:?:Galerkin_orthogonality_time}
\begin{array}{l}
\displaystyle
\sum_{n=1}^{N^\ell}  \sum_{K_n \in \mathcal{T}_{\tau,n}}
\int_{K_n}
\big\{
(
\partial_{t} e,\varphi_{\tau})
+ a(e, \convection_{\sigma})(\varphi_{\tau})
\big\}
\mathrm{d} t
\\[1.5ex]
=
\displaystyle\sum_{t_F \in \mathcal{F}_\tau}
([u_\tau]_{t_F},\varphi_{\tau}(t_F^+))
\displaystyle
- \sum_{n=1}^{N^\ell}\sum_{K_n \in \mathcal{T}_{\tau,n}}\int_{K_n}
\big(
(\convection-\convection_{\sigma}) \cdot \nabla \concentration,
\varphi_{\tau}
\big)
\mathrm{d} t\,,
\end{array}
\end{equation}
with a non-vanishing right-hand side term depending on the the temporal
error in the approximation of the flow field.
For the spatial error $e = \concentration_{\tau}-\concentration_{\tau h}$
we get the following Galerkin orthogonality by subtracting 
Eq.~\eqref{eq:17:fully_transport} from Eq.~\eqref{eq:13:semi_transport} 
\begin{equation}
\label{eq:?:Galerkin_orthogonality_space}
\begin{array}{l}
\displaystyle
\sum_{n=1}^{N^\ell}  \sum_{K_n \in \mathcal{T}_{\tau,n}}
\int_{K_n}
\big\{
(
\partial_{t} e,\varphi_{\tau h})
+ a(e; \convection_{\sigma h})(\varphi_{\tau h})
\big\}
\mathrm{d} t
\\
+  \displaystyle\sum_{t_F \in \mathcal{F}_\tau}
(
[e]_{t_F},\varphi_{\tau h}(t_F^+)
+ (e(0^+),\varphi_{\tau h}(0^+))
\\
=
\displaystyle
S_A(\concentration_{\tau h}; \convection_{\sigma h})(\varphi_{\tau h})
- \sum_{n=1}^{N^\ell}\sum_{K_n \in \mathcal{T}_{\tau,n}}\int_{K_n}
\big(
(\convection_{\sigma}-\convection_{\sigma h}) \cdot \nabla \concentration_{\tau},
\varphi_{\tau h}
\big)
\mathrm{d} t\,,
\end{array}
\end{equation}
with a non-vanishing right-hand side term depending on the stabilization and the
spatial error in the approximation of the flow field.

\subsection*{\textbf{Transport: Dual Problems and Residuals}}

Within the context of the DWR philosophy, the dual problems are generally given as 
optimality or stationary conditions regarding the underlying Lagrangian functionals.
More precisely, considering the directional derivatives of the Lagrangian 
functionals \eqref{eq:20:Lagrangian_transport}, also known as G\^{a}teaux 
derivatives (cf., e.g., \cite{Besier2012}), with respect to their first argument, 
i.e.\
\begin{displaymath}
\mathcal{L}^{\prime}_{\concentration}(\concentration,\dualz;\convection)(\varphi) :=
\lim_{t\neq0,t\rightarrow 0}
t^{-1}\big\{\mathcal{L}(\concentration + t\varphi,\dualz;\convection)
-\mathcal{L}(\concentration,\dualz;\convection)\big\},
\quad \varphi \in \mathcal{X}\,,
\end{displaymath}
leads to the following dual transport problems: 
Find the continuous dual solution $\dualz \in \mathcal{X}$, the semi-discrete
dual solution $\dualz_{\tau} \in \mathcal{X}_{\tau}^{r}$ and the fully
discrete dual solution $\dualz_{\tau h} \in \mathcal{X}_{\tau h}^{r,p}$,
respectively, such that
\begin{equation}
\label{eq:?:DualProblems}
\begin{array}{r@{\,}c@{\,}l@{\,}l@{\,}}
A^{\prime}(\concentration; \convection)(\varphi,\dualz)
&=&
J^{\prime}(\concentration)(\varphi)
&
\quad\forall \varphi \in \mathcal{X}\,,
\\[1.ex]
A_{\tau}^{\prime}(\concentration_{\tau}; \convection_{\sigma})(\varphi_{\tau},\dualz_{\tau})
& = &
J^{\prime}(\concentration_{\tau})(\varphi_{\tau})
&
\quad\forall \varphi_{\tau}\in \mathcal{X}_{\tau}^{r}\,,
\\[1.ex]
A_{S}^{\prime}(\concentration_{\tau h}; \convection_{\sigma h})(\varphi_{\tau h},\dualz_{\tau h})
& = &
J^{\prime}(\concentration_{\tau h})(\varphi_{\tau h})
& 
\quad\forall \varphi_{\tau h}\in \mathcal{X}_{\tau h}^{r,p}\,,
\end{array}
\end{equation}
where we refer to our works \cite{Bause2021,Bruchhaeuser2022a} for a detailed 
description of the adjoint bilinear forms $A^{\prime},A_{\tau}^{\prime},A_S^{\prime}$ 
as well as the dual right hand side term $J^\prime$. 

The primal and dual residuals based on the continuous and semi-discrete in time
schemes are defined by means of the G\^{a}teaux derivatives of the Lagrangian 
functionals in the following way:
\begin{displaymath} 
\begin{array}{r@{\;}c@{\;}l@{\;}c@{\;}l@{\;}}
\rho(u;\convection)(\varphi) & := &  
\mathcal{L}_{\dualz}^{\prime}(\concentration,\dualz;\convection)(\varphi)
& = & G(\varphi)-A(\concentration;\convection)(\varphi)\,,
\\[1.5ex]
\rho^{\ast}(\concentration,\dualz;\convection)(\varphi) & := & 
\mathcal{L}_{\concentration}^{\prime}(\concentration,\dualz;\convection)(\varphi)
& = & J^{\prime}(\concentration)(\varphi)
-A^{\prime}(\concentration;\convection)(\varphi,\dualz)\,,
\\[1.5ex]
\rho_{\tau}(\concentration;\convection_{\sigma})(\varphi) & := &
\mathcal{L}_{\tau,\dualz}^{\prime}(\concentration,\dualz;\convection_{\sigma})(\varphi)
& = & G_{\tau}(\varphi)-A_{\tau}(\concentration;\convection_{\sigma})(\varphi)\,,
\\[1.5ex]
\rho_{\tau}^{\ast}(\concentration,\dualz;\convection_{\sigma})(\varphi) & := &
\mathcal{L}_{\tau,\concentration}^{\prime}(\concentration,\dualz;\convection_{\sigma})(\varphi)
& = & J^{\prime}(\concentration)(\varphi)
-A_{\tau}^{\prime}(\concentration;\convection_{\sigma})(\varphi,\dualz)\,.
\end{array}
\end{displaymath}

\subsection*{\textbf{Flow: Dual Problems and Residuals}}

For the sake of completeness, considering the directional derivatives of the Lagrangian 
functionals \eqref{eq:18:Lagrangian_stokes}, also known as G\^{a}teaux 
derivatives (cf., e.g., \cite{Besier2012}), with respect to their first argument, 
i.e.\
\begin{displaymath}
\mathcal{L}^{\prime}_{\mathbf{u}}(\mathbf{u},\mathbf{z})(\boldsymbol{\varphi}) :=
\lim_{t\neq0,t\rightarrow 0}
t^{-1}\big\{\mathcal{L}(\mathbf{u} + t\boldsymbol{\varphi},\mathbf{z})
-\mathcal{L}(\mathbf{u},\mathbf{z})\big\},
\quad \boldsymbol{\varphi} \in \mathcal{Y}\,,
\end{displaymath}
leads to the following dual flow problems, although these problems were not used 
within the underlying cost reduced approach here: 
Find the continuous dual flow solution $\mathbf{z}=\{\mathbf{w},q\} \in \mathcal{Y}$, 
the semi-discrete dual flow solution $\mathbf{z}_{\sigma}=\{\mathbf{w}_{\sigma},q_{\sigma}\} \in 
\mathcal{Y}_{\sigma}^{r}$ and the fully discrete dual flow solution 
$\mathbf{z}_{\sigma h}\{\mathbf{w}_{\sigma h},q_{\sigma h}\} \in \mathcal{Y}_{\sigma h}^{r,p}$,
respectively, such that
\begin{equation}
\label{eq:?:dual_stokes}
\begin{array}{r@{\,}c@{\,}l@{\,}l@{\,}}
B^{\prime}(\mathbf{u})(\boldsymbol{\varphi},\mathbf{z}) 
&=&
J^{\prime}(\mathbf{u})(\boldsymbol{\varphi}) 
&
\quad\forall  \boldsymbol{\varphi} = \{\boldsymbol{\psi},\chi\} \in \mathcal{Y}\,,
\\[1.ex]
B_{\sigma}^{\prime}(\mathbf{u}_{\sigma})(\boldsymbol{\varphi}_{\sigma},\mathbf{z}_{\sigma})
& = &
J^{\prime}(\mathbf{u}_{\sigma})(\boldsymbol{\varphi}_{\sigma})
&
\quad\forall \boldsymbol{\varphi}_{\sigma} = \{\boldsymbol{\psi}_{\sigma},\chi_{\sigma}\}
\in \mathcal{Y}_{\sigma}^{r}\,, 
\\[1.ex]
B_{\sigma h}^{\prime}(\mathbf{u}_{\sigma h})(\boldsymbol{\varphi}_{\sigma h},\mathbf{z}_{\sigma h})
& = &
J^{\prime}(\mathbf{u}_{\sigma h})(\boldsymbol{\varphi}_{\sigma h})
& \quad\forall \boldsymbol{\varphi}_{\sigma h} = \{\boldsymbol{\psi}_{\sigma h},\chi_{\sigma h}\}\in \mathcal{Y}_{\sigma h}^{r,p}\,,
\end{array}
\end{equation}
where we refer to \cite{Bruchhaeuser2022a} for a detailed description of the
adjoint bilinear forms $B^{\prime},B_{\sigma}^{\prime},B_{\sigma h}^{\prime}$ 
as well as the dual right hand side term $J^\prime$. 

Finally, the primal and dual residuals based on the semi-discrete in time schemes 
are defined by means of the G\^{a}teaux derivatives of the Lagrangian functionals
in the following way:
\begin{displaymath} 
\begin{array}{r@{\;}c@{\;}l@{\;}c@{\;}l@{\;}}
\rho_{\sigma}(\mathbf{u})(\boldsymbol{\varphi}) 
& := &
\mathcal{L}_{\sigma,\mathbf{z}}^{\prime}(\mathbf{u},\mathbf{z})(\boldsymbol{\varphi})
& = & F_{\sigma}(\boldsymbol{\psi})-B_{\sigma}(\mathbf{u})(\boldsymbol{\varphi})
\,,
\\[1.5ex]
\rho_{\sigma}^{\ast}(\mathbf{u},\mathbf{z})(\boldsymbol{\varphi}) 
& := & 
\mathcal{L}_{\sigma,\mathbf{u}}^{\prime}(\mathbf{u},\mathbf{z})(\boldsymbol{\varphi})
& = & J^{\prime}(\mathbf{u})(\boldsymbol{\varphi})-B_{\sigma}^{\prime}(\mathbf{u})(\boldsymbol{\varphi},\mathbf{z})
\,.
\end{array}
\end{displaymath}

\end{document}